W. B. Vasantha Kandasamy
Florentin Smarandache
N. Suresh Babu
R.S. Selvaraj

RANK DISTANCE BICODES AND
THEIR GENERALIZATION

2010

# RANK DISTANCE BICODES AND THEIR GENERALIZATION

W. B. Vasantha Kandasamy
Florentin Smarandache
N. Suresh Babu
R.S. Selvaraj

# CONTENTS

# PREFACE

In this book the authors introduce the new notion of rank distance bicodes and generalize this concept to Rank Distance n-codes (RD n-codes), n, greater than or equal to three. This definition leads to several classes of new RD bicodes like semi circulant rank bicodes of type I and II, semicyclic circulant rank bicode, circulant rank bicodes, bidivisible bicode and so on. It is important to mention that these new classes of codes will not only multitask simultaneously but also they will be best suited to the present computerised era. Apart from this, these codes are best suited in cryptography.

This book has four chapters. In chapter one we just recall the notion of RD codes, MRD codes, circulant rank codes and constant rank codes and describe their properties. In chapter two we introduce few new classes of codes and study some of their properties. In this chapter we introduce the notion of fuzzy RD codes and fuzzy RD bicodes. Rank distance m-codes are introduced in chapter three and the property of m-covering radius is analysed. Chapter four indicates some applications of these new classes of codes.

Our thanks are due to Dr. K. Kandasamy for proof-reading this book. We also acknowledge our gratitude to Kama and Meena for their help with corrections and layout.

W.B.VASANTHA KANDASAMY
FLORENTIN SMARANDACHE
N. SURESH BABU
R.S.SELVARAJ

Chapter One

# BASIC PROPERTIES OF
# RANK DISTANCE CODES

In this chapter we recall the basic definitions and properties of Rank Distance codes (RD codes). The Rank Distance (RD) codes are special type of codes endowed with rank metric introduced by Gabidulin [24, 27]. The rank metric introduced by Gabidulin is an ideal metric for it has the capability of handling varied error patterns efficiently.

    The significance of this new metric is that it recognizes the linear dependence between different symbols of the alphabet. Hence a code equipped with the rank metric detects and corrects more error patterns compared to those codes with other metric. In 1985 Gabidulin has studied a particular class of codes equipped with rank metric called Maximum Rank Distance (MRD) codes. Throughout this book $V^n$ denotes a linear space of dimension n over the Galois field $GF(2^N)$, $N > 1$. By fixing a basis for $V^n$ over $GF(2^N)$, we can represent any element $x \in V^n$ as an n-tuple $(x_1, x_2, \ldots, x_n)$ where $x_i \in GF(2^N)$.

Again, $GF(2^N)$ can be considered as a linear space of dimension 'N' over $GF(2)$. Hence an element $x_i \in GF(2^N)$ has a representation as a N-tuple $(\alpha_{i1}, \alpha_{i2}, \ldots, \alpha_{in})$ over $GF(2)$ with respect to some fixed basis. Hence associated with each $x \in V^n$, $(n \leq N)$ there is a matrix,

$$m(x) = \begin{bmatrix} \alpha_{11} & \cdots & \alpha_{1n} \\ \alpha_{21} & \cdots & \alpha_{2n} \\ \vdots & & \vdots \\ \alpha_{N1} & \cdots & \alpha_{Nn} \end{bmatrix}^T$$

where the $i^{th}$ column represents the $i^{th}$ coordinate '$x_i$' of x over $GF(2)$. If we assume $x_i \in GF(2^N)$ has a representation as a N-tuple $x_i = \alpha_{1i}u_1 + \alpha_{2i}u_2 + \ldots + \alpha_{Ni}u_N$ where $u_1, u_2, \ldots, u_N$ is some fixed basis of the field $F_{q^N}$ regarded as a vector space over $F_q$. If $A_N^n$ denote the ensemble of all $N \times n$ matrices over $F_q$ and if $A: V^n \to A_N^n$ is a bijection defined by the rule, for any vector x $= (x_1, \ldots, x_n) \in F_{q^N}^n$, the associated matrix denoted by,

$$A(x) = \begin{bmatrix} a_{11} & a_{12} & \cdots & a_{1n} \\ a_{21} & a_{22} & \cdots & a_{2n} \\ \vdots & \vdots & & \vdots \\ a_{N1} & a_{N2} & \cdots & a_{Nn} \end{bmatrix}$$

where the $i^{th}$ column represents the $i^{th}$ coordinate $x_i$ of x over $F_q$.

Now we proceed on to define the rank of an element $x \in V^n$.

**DEFINITION 1.1:** *The rank of a vector $x \in F_{q^N}^n$ is the rank of the associated matrix A(x). Let r(x) denote the rank of the vector x $\in F_{q^N}^n$, over $F_q$.*

By the usual properties of the rank of a matrix, it is easy to prove the following inequalities.

1. $r(x) \geq 0$ for every $x \in V^n$.
2. $r(x) = 0$ if and only if $x = 0$.
3. $r(x+y) \leq r(x) + r(y)$ for every $x, y \in V^n$.
4. $r(ax) = |a| r(x)$ for every $a \in GF(2)$ and $x \in V^n$.

Thus the function $x \to r(x)$ defines a norm on $V^n$, ($x \to r(x)$ defines a norm on $F_{q^N}^n$) and is called the rank norm. In this book we denote the rank norm by $r(x)$ or $wt(x)$ or by $\|x\|$. The rank norm induces a metric called rank metric (or rank distance) on $F_{q^N}^n$.

**DEFINITION 1.2:** *The metric induced by the rank norm is defined as the rank metric on $V^n$ ($F_{q^N}^n$) and it is denoted by $d_R$. The rank distance between $x, y \in V^n$ is the rank of their difference : $d_R(x, y) = r(x - y)$. The vector space $V^n$ ($F_{q^N}^n$) over $F_{q^N}$ equipped with the rank metric $d_R$ is defined as a rank distance space.*

**DEFINITION 1.3:** *A linear space $V^n$ over $GF(2^N)$, $N > 1$ of dimension $n$ such that $n \leq N$, equipped with the rank metric is defined as a rank space or rank distance space.*

Now we proceed on to recall the definition of Rank Distance (RD) codes.

**DEFINITION 1.4:** *A Rank Distance code (RD-code) of length $n$ over $GF(2^N)$ is a subset of the rank space $V^n$ over $GF(2^N)$. A linear $[n, k]$ RD code is a linear subspace of dimension $k$ in the rank space $V^n$. By $C[n, k]$ we denote a linear $[n, k]$ RD-code.*

**DEFINITION 1.5:** *A generator matrix of a linear $[n, k]$ RD code $C$ is a $k \times n$ matrix over $GF(2^N)$ whose rows form a basis for $C$. A generator matrix $G$ of a linear RD code $C[n, k]$ can be brought into the form $G = [I_k, A_{k, n-k}]$ where $I_k$ is the identity matrix and $A_{k, n-k}$ is some matrix over $GF(2^n)$. This form is called the standard form.*

**DEFINITION 1.6:** *Let G be a generator matrix of the linear RD code C[n, k], then a matrix H of order $(n - k) \times n$ over $GF(2^N)$ such that $GH^T = (0)$ is called a parity check matrix of C[n, k]. Suppose C is a linear [n, k] RD code with G and H as its generator and parity check matrices respectively, then C has two representations*
      *1. C is the row space of G and*
      *2. C is the solution space of H.*

We shall illustrate this situation by an example.

***Example 1.1:*** Let

$$G = \begin{bmatrix} \alpha_{11} & \alpha_{12} & ... & \alpha_{15} \\ \alpha_{21} & \alpha_{22} & ... & \alpha_{25} \\ \alpha_{31} & \alpha_{32} & ... & \alpha_{35} \end{bmatrix}_{3 \times 5}$$

be a generator matrix of the linear [5, 3] Rank Distance code C; over $GF(2^5)$, here $(\alpha_{11}, \alpha_{12}, \ldots, \alpha_{15})$, $(\alpha_{21}, \alpha_{22}, \ldots, \alpha_{25})$, $(\alpha_{31}, \alpha_{32}, \ldots, \alpha_{35})$ forms a basis of C; $\alpha_{ij} \in GF(2^5)$; $1 \le i \le 3$ and $1 \le j \le 5$.

Now as in case of linear codes with Hamming metric, we in case of Rank Distance codes have the concept of minimum distance. We just recall the definition.

**DEFINITION 1.7:** *Let C be a rank distance code, the minimum - rank distance d is defined by $d = \min\{d_R(x, y) \mid x, y \in C, x \ne y\}$. In other words, $d = \min\{r(x - y) \mid x, y \in C, x \ne y\}$. i.e., $d = \min\{r(x) \mid x \in C \text{ and } x \ne 0\}$.*

*If an RD code C has the minimum – rank distance d then it can correct all errors $e \in F_{q^N}^n$ with rank*

$$r(e) = \left\lfloor \frac{d-1}{2} \right\rfloor.$$

*Let C denote an [n, k] RD – code over $F_{q^N}$. A generator matrix G of C is a $k \times n$ matrix with entries from $F_{q^N}^n$ whose rows form* *a basis for C. Then an $(n - k) \times n$ matrix H with entries from $F_{q^N}^n$ such that $GH^T = (0)$ is called the parity check matrix of C.*

Result (singleton - style bound) The minimum - rank distance d of any linear [n, k] RD code $C \subseteq F_{q^N}^n$ satisfies the following bound: $d \leq n - k + 1$.

Now based on this, the notion of Maximum Rank Distance; MRD codes were defined in [24, 27].

**DEFINITION 1.8:** *An [n, k, d] RD code C is called Maximum Rank Distance (MRD) code if the singleton – style bound is reached; i.e., if $d = n - k + 1$.*

Now we just briefly recall the construction of MRD code.

Let $[s] = q^s$ for any integer s. Let $g_1, \ldots, g_n$ be any set of elements in $F_{q^N}$ that are linearly independent over $F_q$.

A generator matrix G of an MRD code C is defined by

$$G = \begin{bmatrix} g_1 & g_2 & \cdots & g_n \\ g_1^{[1]} & g_2^{[1]} & \cdots & g_n^{[1]} \\ g_1^{[2]} & g_2^{[2]} & \cdots & g_n^{[2]} \\ \vdots & \vdots & & \vdots \\ g_1^{[k-1]} & g_2^{[k-1]} & \cdots & g_n^{[k-1]} \end{bmatrix}.$$

It can be shown that the code C given by the above generator matrix G has the rank distance $d = n - k + 1$.

Any matrix of the above form is called a Frobenius matrix with generating vector $g_c = (g_1, g_2, \ldots, g_n)$.

Gabidulin has proved the following theorem.

**THEOREM 1.1:** *Let C[n, k] be a linear (n, k, d) MRD-code with d = 2t + 1. Then C[n, k] corrects all errors of rank atmost t and detects all errors of rank greater than t.*

Circulant Rank codes were defined by [61].

Consider the Galois field $GF(2^N)$ where $N > 1$. An element $\alpha \in GF(2^N)$ can be denoted by a N-tuple $(a_0, a_1, \ldots, a_{N-1})$ as well as by a polynomial $a_0 + a_1x + \ldots + a_{N-1}x^{N-1}$ over $GF(2)$.

**DEFINITION 1.9:** *The circulant transpose ($T_c$) of a vector $\alpha = (a_0, a_1, \ldots, a_{N-1}) \in GF(2^N)$ is defined as $\alpha^{T_c} = (a_0, a_1, \ldots, a_{N-1})$. If $\alpha \in GF(2^N)$ has polynomial representation $a_0 + a_1x + \ldots + a_{N-1}x^{N-1}$ in $\dfrac{[GF(2)](x)}{(x^N + 1)}$ then by $\alpha_i$, we denote the vector corresponding to the polynomial $[(a_0 + a_1x + \ldots + a_{N-1}x^{N-1}).x^i]$ (mod $x^N + 1$), for $i = 0$ to $N - 1$. (Note $\alpha_0 = \alpha$).*

**DEFINITION 1.10:** *Let $f: GF(2^N) \to [GF(2^N)]^N$ be defined as $f(\alpha) = (\alpha_0^{T_c}, \alpha_1^{T_c}, \ldots, \alpha_{N-1}^{T_c})$; we call $f(\alpha)$ as the 'word' generated by $\alpha$.*

MacWilliams F.J. and Sloane N.J.A., [61] defined circulant matrix associated with a vector in $GF(2^N)$ as follows.

**DEFINITION 1.11:** *A matrix of the form*

$$\begin{bmatrix} a_0 & a_1 & \ldots & a_{N-1} \\ a_{N-1} & a_0 & \ldots & a_{N-2} \\ \vdots & \vdots & & \vdots \\ a_1 & a_2 & \ldots & a_0 \end{bmatrix}$$

*is called the circulant matrix associated with the vector $(a_0, a_1, \ldots, a_{N-1}) \in GF(2^N)$. Thus with each $\alpha \in GF(2^N)$ we can associate a circulant matrix whose $i^{th}$ column represents $\alpha_i^{T_c}$, $i = 0, 1, 2, \ldots, N - 1$. f is nothing but a mapping of $GF(2^N)$ on to* *the algebra of all N × N circulant matrices over GF(2). Denote the space f(GF($2^N$)) by $V^N$.*

We define norm of a word $v \in V^N$ as follows:

**DEFINITION 1.12:** *The 'norm' of a word $v \in V^N$ is defined as the 'rank' of v over GF(2) (By considering it as a circulant matrix over GF(2)).*

We denote the 'norm' of v by r(v). We just prove the following theorem.

**THEOREM 1.2:** *Suppose $\alpha \in GF(2^N)$ has the polynomial representation g(x) over GF(2) such that the gcd(g(x), $x^N + 1$) has degree N – k, where $0 \le k \le N$. Then the 'norm' of the word generated by $\alpha$ is 'k'.*

*Proof:* We know the norm of the word generated by $\alpha$ is the rank of the circulant matrix $(\alpha_0^{T_C}, \alpha_1^{T_C}, ..., \alpha_{N-1}^{T_C})$, where $\alpha_i^{T_C}$ represents the polynomial $[x^i g(x)]$ (mod $x^N + 1$) over GF(2).

Suppose the gcd(g(x), $x^N + 1$) is a polynomial of degree 'N – k', $(0 \le k \le N – 1)$. To prove that the word generated by '$\alpha$' has rank 'k'. It is enough to prove that the space generated by the N-polynomials g(x)(mod $x^N + 1$), [x·g(x)](mod $x^N + 1$), …, [$x^{N-1}$·g(x)] (mod $x^N + 1$) has dimension 'k'. We will prove that the set of k-polynomials g(x)(mod $x^N + 1$), [x·g(x)] (mod $x^N + 1$), …, [$x^{N-1}$·g(x)] (mod $x^N + 1$) forms a basis for this space.

If possible, let $a_0(g(x)) + a_1(x \cdot g(x)) + a_2(x^2 \cdot g(x)) + ... + a_{k-1}(x^{k-1} \cdot g(x)) \equiv 0 (\mod (x^N + 1))$, where $a_i \in GF(2)$. This implies $x^N + 1$ divides $(a_0 + a_1 x + ... + a_{k-1} x^{k-1}) \cdot g(x)$.

Now if g(x) = h(x) a(x) where h(x) is the gcd(g(x), $x^N + 1$), then (a(x), $x^N + 1$) = 1. Thus $x^N + 1$ divides $(a_0 + a_1 x + ... + a_{k-1} x^{k-1}) \cdot g(x)$ implies that the quotient

$$\frac{(x^N + 1)}{h(x)}$$

divides $(a_0 + a_1 x + ... + a_{k-1} x^{k-1}) \cdot a(x)$. That is

$$\left[\frac{(x^N+1)}{h(x)}\right]$$

divides $(a_0 + a_1x + \ldots a_{k-1}x^{k-1})$ which is a contradiction as

$$\left[\frac{(x^N+1)}{h(x)}\right]$$

has degree k where as the polynomial $(a_0 + a_1x + \ldots + a_{k-1}x^{k-1})$ has degree atmost $k - 1$.

Hence the polynomials $g(x)\bmod(x^N+1)$, $[x \cdot g(x)]\bmod(x^N+1)$, ..., $[x^{k-1} \cdot g(x)]\bmod(x^N+1)$ are linearly independent over GF(2). We will prove that the polynomials $g(x)\bmod(x^N+1)$, $[x \cdot g(x)]\bmod(x^N+1)$, ..., $[x^{k-1} \cdot g(x)]\bmod(x^N+1)$ generate the space. For this, it is enough to prove that $x^i \cdot g(x)$ is a linear combination of these polynomials for $k \leq i \leq N - 1$.

Let $x^N + 1 = h(x)b(x)$, where $b(x) = b_0 + b_1x + \ldots + b_kx^k$ (Note that here $b_0 = b_k = 1$, since $b(x)$ divides $x^N + 1$).

Also, we have $g(x) = h(x) \cdot a(x)$. Thus

$$x^N + 1 = \frac{(g(x) \cdot b(x))}{a(x)}$$

i.e.,

$$\frac{g(x)(b_0 + b_1x + \ldots + b_kx^k)}{a(x)} = 0 \bmod(x^N + 1),$$

that is

$$\frac{g(x) \cdot (b_0 + b_1x + \ldots + b_kx^{k-1})}{a(x)} = \left[\frac{(g(x) \cdot x^k)}{a(x)}\right] \bmod(x^N + 1).$$

(Since $b_k = 1$). That is

$x^k \cdot g(x) = (b_0 + b_1x + \ldots + b_{k-1}x^{k-1} \cdot g(x)) \bmod(x^N + 1)$. Hence $x^k g(x) = b_0 g(x) + b_1[x \cdot g(x)] + \ldots + b_{k-1}[x^{k-1} \cdot g(x)] (\bmod(x^N + $

1)), is a linear combination of $g(x) \bmod(x^N + 1)$, $[x \cdot g(x)] \bmod(x^N + 1)$, …, $[x^{k-1} \cdot g(x)] \bmod(x^N + 1)$ over GF(2).

Now it can be easily proved that $x^i g(x)$ is a linear combination of $g(x) \bmod(x^N + 1)$, $[x \cdot g(x)] \bmod(x^N + 1)$, …, $[x^{k-1} \cdot g(x)] \bmod(x^N + 1)$ for $i > k$.

Hence the space generated by the polynomial $g(x) \bmod(x^N + 1)$, $[x \cdot g(x)] \bmod(x^N + 1)$, …, $[x^{k-1} \cdot g(x)] \bmod(x^N + 1)$ has dimension k; i.e., the rank of the word generated by $\alpha$ is k.

The following two corollaries are obvious.

**COROLLARY 1.1:** *If $\alpha \in GF(2^N)$ is such that its polynomial representation g(x) is relatively prime to $x^N + 1$, then the norm of the word generated by $\alpha$ is N and hence $f(\alpha)$ is invertible.*

*Proof:* Follows from the theorem as $\gcd(g(x), x^N +1) = 1$ has degree 0 and hence rank of $f(\alpha)$ is N.

**COROLLARY 1.2:** *The norms of the vectors corresponding to the polynomials $x +1$ and $x^{N-1} + x^{N-2} + … + x +1$ are respectively $N – 1$ and 1.*

Now we proceed on to define the distance function on $V^N$.

**DEFINITION 1.13:** The distance between two words v, u. in $V^N$ is defined as $d(u, v) = r(u + v)$.

Now we define circulant code of length N.

**DEFINITION 1.14:** *A circulant rank code of length N is defined as a subspace of $V^N$ equipped with the above defined distance function.*

**DEFINITION 1.15:** *A circulant rank code of length N is called cyclic if, whenever $(v_1, v_2, …, v_N)$ is a codeword, then it implies $(v_2, v_3, …, v_N, v_1)$ is also a codeword.*

Now we proceed on to recall the definition of the new class of codes, Almost Maximum Rank Distance codes (AMRD-codes).

**DEFINITION 1.16:** *A linear [n, k] RD code over $GF(2^N)$ is called Almost Maximum Rank Distance (AMRD) code if its minimum distance is greater than or equal to n – k.*

*An AMRD code whose minimum distance is greater than n – k is an MRD code and hence the class of MRD codes is a subclass of the class of AMRD codes.*

We recall the following theorem.

**THEOREM 1.3:** *When n – k is an odd integer,*

1. *The error correcting capability of an [n, k] AMRD code is equal to that of an [n, k] MRD code.*
2. *An [n, k] AMRD code is better than any [n, k] code in Hamming metric for error correction.*

*Proof:* (1) Suppose C be an [n, k] AMRD code such that 'n – k' is an odd integer. The maximum number of errors corrected by C is given by $\frac{(n-k-1)}{2}$. But $\frac{(n-k-1)}{2}$ is equal to the error correcting capability of an [n, k] MRD code (since n – k is odd). Thus when n – k is odd an [n, k] AMRD code is as good as an [n, k] MRD code.

(2) Suppose C be an [n, k] AMRD code such that 'n – k' is odd, then, each codeword of C can correct $L_r(n)$ error vectors where

$$r = \frac{(n-k-1)}{2}$$

and

$$L_r(n) = 1 + \sum_{i=1}^{n} \begin{bmatrix} n \\ i \end{bmatrix} (2^N - 1) \ldots (2^N - 2^{i-1}).$$

Consider the same [n, k] code in Hamming metric. Let it be $C_1$, then the minimum distance of $C_1$ is atmost (n – k +1). The error correcting capability of $C_1$ is

$$\left\lfloor \frac{(n-k+1-1)}{2} \right\rfloor = \frac{(n-k-1)}{2} = r$$

(since n – k is odd).

Hence the number of error vectors corrected by a codeword is given by

$$\sum_{i=0}^{r} \binom{n}{i} (2^N - 1)^i,$$

which is clearly less than $L_r(n)$. Thus the number of error vectors that can be corrected by the [n, k] AMRD code is much greater than that of the same code considered in Hamming metric.

For any given length 'n', a single error correcting AMRD code is one having dimension n – 3 and minimum distance greater than or equal to '3'. We give a characterization for a single error correcting AMRD codes in terms of its parity check matrix. This characterization is based on the condition for the minimum distance proved by Gabidulin in [24, 27].

We just recall the main theorem for more about these properties one can refer [82].

**THEOREM 1.4:** *Let $H = (\alpha_{ij})$ be a $3 \times n$ matrix of rank 3 over $GF(2^N)$, $n \leq N$ which satisfies the following condition. For any two distinct, non empty subsets $P_1$, $P_2$ of $\{1, 2, …, n\}$ there exists $i_1, i_2 \in \{1, 2, 3\}$ such that*

$$\left( \sum_{j \in P_1} \alpha_{i_1 j} \cdot \sum_{k \in P_2} \alpha_{i_2 k} \right) \neq \left( \sum_{j \in P_1} \alpha_{i_2 j} \cdot \sum_{k \in P_2} \alpha_{i_1 k} \right)$$

*then, H as a parity check matrix defines a [n, n – 3] single error correcting AMRD code over $GF(2^N)$.*

*Proof:* Given H is a $3 \times n$ matrix of rank 3 over $GF(2^N)$, so H as a parity check matrix defines a [n, n – 3] RD code C over $GF(2^N)$; where $C = \{x \in V^n \mid xH^T = 0\}$.

It remains to prove that the minimum distance of C is greater than or equal to 3. We will prove that no non-zero codeword of C has rank less than '3'. The proof is by the method of contradiction.

Suppose there exists a non-zero codeword x such that $r(x) \leq 2$; then, x can be written as $x = y \cdot M$ where $y = (y_1, y_2)$; $y_i \in GF(2^N)$ and $M = (m_{ij})$ is a $2 \times n$ matrix of rank 2 over $GF(2)$.

Thus $(y \cdot M)H^T = 0$ implies $y(MH^T) = 0$. Since y is non zero $y(M \cdot H^T) = 0$ implies that the $2 \times 3$ matrix $MH^T$ has rank less than 2 over $GF(2^N)$. Now let $P_1 = \{j$ such that $m_{1j} = 1\}$ and $P_2 = \{j$ such that $m_{2j} = 1\}$. Since $M = (m_{ij})$ is a $2 \times n$ matrix of rank 2, $P_1$ and $P_2$ are disjoint nonempty subsets of $\{1, 2, \ldots, n\}$, and

$$MH^T = \begin{pmatrix} \sum_{j \in P_1} \alpha_{1j} & \sum_{j \in P_1} \alpha_{2j} & \sum_{j \in P_1} \alpha_{3j} \\ \sum_{j \in P_2} \alpha_{1j} & \sum_{j \in P_2} \alpha_{2j} & \sum_{j \in P_2} \alpha_{3j} \end{pmatrix}.$$

But the selection of H is such that there exists $i_1, i_2 \in \{1, 2, 3\}$ such that

$$\left( \sum_{j \in P_1} \alpha_{i_1 j} \cdot \sum_{k \in P_2} \alpha_{i_2 k} \right) \neq \left( \sum_{j \in P_1} \alpha_{i_2 j} \cdot \sum_{k \in P_2} \alpha_{i_1 k} \right).$$

Hence in $MH^T$ there exists a $2 \times 2$ submatrix whose determinant is nonzero; i.e., $r(MH^T) = 2$ over $GF(2^N)$, this contradicts the fact that rank $(MH^T) < 2$. Hence the result.

Analogous to the constant weight codes in Hamming metric, we define the constant rank codes in rank metric. A constant weight code C of length n over a Galois field F is a subset of $F^n$ with the property that all codewords in C have the same Hamming weight [61]. $A(n, d, w)$ denotes the maximum number of vectors in $F^n$, distance atleast d apart from each other and constant Hamming weight 'w'.

Obtaining bounds for A(n, d, w) is one of the problems in the study of constant weight codes. A number of important bounds on A(n, d, w) were obtained in [61].

Here we just define the constant rank codes in rank metric and analyze the function A(n, r, d) which is the analog of the A(n, d, w) and obtain some interesting results.

**DEFINITION 1.17:** *A constant rank code of length n is a subset of a rank space $V^n$ with the property that every codeword has same rank.*

**DEFINITION 1.18:** *A(n, r, d) is defined as the maximum number of vectors in $V^n$, constant rank r and the distance between any two vectors is atleast d.*

(By a (n, r, d) set, we mean a subset of vectors in $V^n$ having constant rank r and distance between any two vectors is at least d).

We analyze the function A(n, r, d).

**THEOREM 1.5:**
1. $A(n, r, 1) = L_r(n)$, the number of vectors of rank r in $V^n$.
2. $A(n, r, d) = 0$ if $r > 0$ or $d > n$ or $d > 2r$.

*Proof.* (1) is obvious from the fact that $L_r(n)$ is the number of vectors of length n, constant rank r and the distance between any two distinct vectors in the rank space $V^n$ is always greater than or equal to one.
(2) Follows immediately from the definition of A(n, r, d).

**THEOREM 1.6:** $A(n, 1, 2) = 2^n - 1$ over any Galois field $GF(2^N)$.

*Proof:* Let $V_1$ denote the set of vectors of rank 1 in $V^n$. We know for each non zero element $\alpha \in GF(2^N)$ there exists $(2^n - 1)$ vectors of rank one having $\alpha$ as a coordinate. Thus the cardinality of $V_1$ is $(2^N - 1)(2^n - 1)$.

Now divide $V_1$ into $(2^n – 1)$ blocks of $(2^N – 1)$ vectors such that each block consists of the same pattern of all non-zero elements of $GF(2^N)$.

Thus from each block atmost one vector can be chosen such that the selected vectors are atleast rank 2 apart from each other. Such a set we call as a (n, 1, 2) set. Also it is always possible to construct such a set. Hence $A(n, 1, 2) = 2^n – 1$.

We give an example of A(n, 1, 2) set for a fixed N and n as follows.

*Example 1.2:* Let $N = 3$, we use the following notation to define $GF(2^3)$. Let 0, 1, 2, 3 be the basic symbols. Then $GF(2^3)=\{0, 1, 2, 3, (12), (13), (23), (123)\}$ (Note that by (ijk), we denote the linear combination of $i + j + k$ over GF(2)).

Suppose $n = 3$, divide $F^3$ into $2^3 – 1$ blocks of $2^3 – 1$ vectors as follows:

| 001 | 010 | 100 | 110 |
|---|---|---|---|
| 002 | 020 | 200 | 220 |
| ⋮ | ⋮ | ⋮ | ⋮ |
| (00) (123) | 0(123)0 | (123)00 | (123)(123)0 |

| 101 | 011 | 111 |
|---|---|---|
| 202 | 022 | 222 |
| ⋮ | ⋮ | ⋮ |
| (123) 0 (123) | 0(123) (123) | (123) (123) (123) |

It is clear from this arrangement that atmost one vector from each block can be selected to form a (3, 1, 2) set. Also, the following set is a (3, 1, 2) set. {001, 020, 300, (12)(12)0, (13)0(13), 0(23)(23), (123)(123)(123)}. Thus $A(3, 1, 2) = 7 = 2^3 – 1$.

Now we prove another interesting theorem.

**THEOREM 1.7:** *$A(n, n, n) = 2^N - 1$ over any $GF(2^N)$.*

*Proof:* Denote by $V_n$ the set of vectors of rank n in the space $V^n$. We know that the cardinality of $V_n$ is $(2^N - 1)(2^N - 2) \ldots (2^N - 2^{n-1})$. By definition in a (n, n, n) set the distance between any two vectors should be n. Thus no two vectors can have a common symbol at a coordinate place i ($1 \leq i \leq n$). This implies that $A(n, n, n) \leq 2^N - 1$.

Now we construct a (n, n, n) set as follows:

Select N vectors from $V_n$ such that,
1. Each basis element of $GF(2^n)$ should occur (can be as a combination) atleast once in each vector.
2. If the $i^{th}$ vector is choosen $(i + 1)^{th}$ vector should be selected such that its rank distance from any linear combination of the previous i vectors is n.

Now the set of all linear combinations of these N vectors over GF(2) will be such that the distance between any two vectors is n.

Hence it is a (n, n, n) set. Also the cardinality of this (n, n, n) set is $2^N - 1$. (We do not count the all zero linear combination). Thus $A(n, n, n) = 2^N - 1$.

We illustrate this by the following example.

*Example 1.3:* Consider $GF(2^N)$ for any $N > 1$. As in example 1.2, we represent $GF(2^N)$ as a linear combination of the symbols 1, 2, …, N over GF(2).

Let n = 2. We construct a (2, 2, 2) set by taking the set of all linear combinations of the N vectors choosen as follows:

We have two cases to be considered separately, when N is an odd integer say 2k + 1 and when N is an even integer say 2k.

Case 1. N is an odd integer say 2k + 1. In this case choose the N vectors as 12, 23, 34, …, (2k +1) (12).

Case 2. N is an even integer say 2k.

In this case choose the N vectors as 1(12), 21, 3(34), 43, …, (2k – 1) ((2k – 1)2k), (2k) (2k – 1).

Consider the set of linear combinations of the N-vectors. It can be verified easily that this set is a (2, 2, 2) set.

Now we obtain the value of A(n, r, d) for a particular triple, when n = 4, r = 2 and d = 4 in the following theorem.

**THEOREM 1.8:** *A(4, 2, 4) = 5 over any Galois field $GF(2^N)$.*

*Proof:* Consider $V^4$, the 4-dimensional space over $GF(2^N)$. We denote the elements of $V^4$ as 4-tuples (a b c d) where a, b, c, d ∈ $GF(2^N)$.

Denote by $V_2$, the set of vectors of rank 2 in $V^4$. The cardinality of $V_2$ is $35 \times (2^N - 1)(2^N - 2)$. (since

$$|V_2| = L_2(4) = \frac{(2^4-1)(2^4-2)(2^N-1)(2^N-2)}{(2^2-1)(2^2-2)}.$$

Thus for each distinct non-zero pair of elements of $GF(2^N)$ there are 35 vectors in $V_2$.

Let a, b ∈ $GF(2^N)$ be such that a ≠ 0, b≠0 and a ≠b. Divide the set of 35 vectors containing a, b and the linear combination (ab) into six blocks as follows:

| I | II | III | IV | V | VI |
|---|----|-----|----|----|-----|
| 00ab<br>0a0b<br>0ab0<br>0aab<br>0aba<br>0baa<br>0ab(ab) | a00b<br>ab00<br>aa0b<br>ab0a<br>ab0b<br>ab0(ab) | a0ab<br>aba(ab)<br>ab(ab)a<br>ab(ab)(ab)<br>a0ba<br>a0bb | aaab<br>abab<br>abba<br>aaab<br>aaba<br>abaa<br>abbb | aab0<br>aba0<br>abb0<br>aab(ab)<br>abbab | ab(ab)b<br>a0b(ab)<br>ab(ab)0<br>a0b0 |

From the arrangement of these six blocks it can be verified easily that atmost 5 vectors can be chosen to form any (4, 2, 4) set.

For example if we choose a vector of pattern 00ab then no vector of any other pattern from block I can be chosen (otherwise distance between the two is < 4).

Now, move to block II. The first pattern in block II cannot be chosen. So select a vector in the second pattern ab00. No other pattern can be selected from block II.

Now move to block III. Here also the first pattern cannot be chosen. So choose a vector of the pattern aba(ab). Similarly from block IV select the pattern abab. In the block V, the first four pattern cannot be selected. Hence select the pattern abb(ab).

Now move to block VI. But no pattern can be selected from block VI since each pattern is at a distance less than four from one of the already selected patterns. Similarly we can exhaust all the possibilities. Hence a (4, 2, 4) set in this space can have atmost five vectors.

Also it is always possible to choose five vectors in different patterns to form a set (4, 2, 4) set. Thus A(4, 2, 4) = 5.

Now we proceed on to give a general bound for A(n, n, d).

**THEOREM 1.9:** *$A(n, n, d) \leq (2^N - 1)(2^N - 2) \ldots (2^N - 2^{n-d})$ over any field $GF(2^N)$.*

*Proof:* Let $V_n$ be the set of vectors of rank n in the space $V^n$. The cardinality of $V_n$ is given by

$$(2^N - 1)(2^N - 2) \ldots (2^N - 2^{n-1}).$$

Now for a (n, n, d) set two vectors should be different atleast by d coordinate places.

Thus the cardinality of any (n, n, d) set is less than or equal to

$$(2^N - 1)(2^N - 2) \ldots (2^N - 2^{n-d}).$$

This proves,

$$A(n, n\ d) \leq (2^N - 1)(2^N - 2) \ldots (2^N - 2^{n-d}).$$

These AMRD codes are useful for error correction in data storage systems.

Chapter Two

# RANK DISTANCE BICODES AND THEIR PROPERTIES

In this chapter we define for the first time the notion of Rank Distance Bicodes; RD-Bicodes and derive some interesting results about them. The properties of bivector spaces and bimatrices can be had from the books [91-93].

As the error correcting capability of a code depends mainly on the distance between codewords, not only choosing an appropriate metric is important but also simultaneous working of a pair of system would be advantageous in this computerized world. This is done by introducing the concepts of rank distance bicodes, maximum rank distance bicodes, circulant rank bicodes, RD-MRD bicodes, RD-circulant bicodes, MRD circulant bicodes, RD-AMRD bicodes, AMRD bicodes and so on. We aim to give certain classes of new bicodes with rank metric.

Let $V^n$ and $V^m$ be n-dimensional and m-dimensional vector spaces over the field $F_{q^N}$; $n \leq N$ and $m \leq N$ ($m \neq n$). That is $V^n = F_{q^N}^n$ and $V^m = F_{q^N}^m$.

We know $V = V^n \cup V^m$ is a (m, n) dimensional vector bispace (or bivector space) over the field $F_{q^N}$. By fixing a bibasis of $V = V^n \cup V^n$ over $F_{q^N}$ we can represent any element $x \cup y \in V = (V^n \cup V^m)$ as a (n, m)-tuple; $(x_1, \ldots, x_n) \cup (y_1, \ldots, y_m)$ where $x_i, y_j \in F_{q^N}$; $1 \leq i \leq n$ and $1 \leq j \leq m$. Again $F_{q^N}$ can be considered as a pseudo false linear bispace of dimension (N, N) over $F_q$ (we say a linear vector bispace $V = V^n \cup V^n$ to be a pseudo false linear bispace if $m = n = N$). Hence the elements $x_i, y_j \in F_{q^N}$ has a representation as N-bituple $(\alpha_{1i}, \ldots, \alpha_{Ni}) \cup (\beta_{1j}, \ldots, \beta_{Nj})$ over $F_q$ with respect to some fixed bibasis. Hence associated with each $x \cup y \in V^n \cup V^m$ ($n \neq m$) there is a bimatrix

$$m(x) \cup m(y) = \begin{bmatrix} a_{11} & \ldots & a_{1n} \\ a_{21} & \ldots & a_{2n} \\ \vdots & & \vdots \\ a_{N1} & \ldots & a_{Nn} \end{bmatrix} \cup \begin{bmatrix} b_{11} & \ldots & b_{1m} \\ b_{21} & \ldots & b_{2m} \\ \vdots & & \vdots \\ b_{N1} & \ldots & b_{Nm} \end{bmatrix}$$

where the $i^{th} \cup j^{th}$ column represents the $i^{th} \cup j^{th}$ coordinate of $x_i \cup y_j$ of $x \cup y$ over $F_q$.

*Remark:* In order to develop the new notion of rank distance bicodes and trying to give the bimatrices and biranks associated with them we are forced to define the notion of pseudo false bivector spaces.

For example; $V = Z_2^5 \cup Z_2^5$ is a false pseudo bivector space over $Z_2$. Likewise $Z_3^7 \cup Z_3^7$ is a pseudo false bivector space over $Z_3$. Also $V = Z_{7^8}^5 \cup Z_{7^8}^5$ is a pseudo false bivector space over $Z_{7^8}$.

However throughout this book we will be using only vector bispaces over $Z_2$ or $Z_{2^N}$ unless otherwise specified.

Now we see to every $x \cup y$ in the bivector space $V^n \cup V^m$ we have an associated bimatrix $m(x) \cup m(y)$.

We now proceed on to define the birank of the bimatrix $m(x) \cup m(y)$ over $F_q$ or $GF(2)$.

**DEFINITION 2.1:** *The birank of an element $x \cup y \in (V^n \cup V^m)$ is defined as the birank of the bimatrix $m(x) \cup m(y)$ over $GF(2)$ or $F_q$. (the birank of the bimatrix $m(x) \cup m(y)$ is the rank of $m(x) \cup$ rank of $m(y)$).*

*We shall denote the birank of $x \cup y$ by $r_1(x) \cup r_2(y) = r(x \cup y)$, we see analogous to the properties of rank we can in case of the birank of a bimatrix prove the following:*

(i) For every $x \cup y \in (V^n \cup V^m)$ ($x \in V^n$ and $y \in V^m$) we have $r(x \cup y) = r_1(x) \cup r_2(y) \geq 0 \cup 0$ (i.e., each $r_1(x) \geq 0$ and each $r_2(y) \geq 0$ for every $x \in V^n$ and $y \in V^m$).

(ii) $r(x \cup y) = r_1(x) \cup r_2(y) = 0 \cup 0$ if and only if $x \cup y = 0 \cup 0$ i.e., $x = 0$ and $y = 0$.

(iii) $r((x_1 + x_2) \cup (y_1 + y_2)) \leq \{r_1(x_1) + r_1(x_2)\} \cup r_2(y_1) + r_2(y_2)$ for every $x_1, x_2 \in V^n$ and $y_1, y_2 \in V^m$. That is we have $r((x_1 + x_2) \cup (y_1 + y_2)) = r_1(x_1 + x_2) \cup r_2(y_1 + y_2) \leq r_1(x_1) + r_1(x_2) \cup r_2(y_1) + r_2(y_2)$; (as we have for every $x_1, x_2 \in V^n$, $r(x_1 + x_2) \leq r_1(x_1) + r_1(x_2)$ and for every $y_1, y_2 \in V^m$; $r_2(y_1 + y_2) \leq r_2(y_1) + (y_2)$).

(iv) $r_1(a_1 x) \cup r_2(a_2 y) = |a_1| r_1(x) \cup |a_2| r_2(y)$ for every $a_1, a_2 \in F_q$ or $GF(2)$ and $x \in V^n$ and $y \in V^m$.

*Thus the bifunction $x \cup y \to r_1(x) \cup r_2(y)$ defines a binorm on $V^n \cup V^m$.*

**DEFINITION 2.2:** *The bimetric induced by the birank binorm is defined as the birank bimetric on $V^n \cup V^m$ and is denoted by*

$d_{R_1} \cup d_{R_2}$. If $x_1 \cup y_1$, $x_2 \cup y_2 \in V^n \cup V^m$ then the birank bidistance between $x_1 \cup y_1$ and $x_2 \cup y_2$ is

$$d_{R_1}(x_1, x_2) \cup d_{R_2}(y_1, y_2) = r_1(x_1 - x_2) \cup r_2(y_1 - y_2)$$

(here $d_{R_1}(x_1, x_2) = r_1(x_1 - x_2)$ for every $x_1$, $x_2$ in $V^n$, the rank distance between $x_1$ and $x_2$ likewise for $y_1, y_2 \in V^m$).

**DEFINITION 2.3:** *A linear bispace $V^n \cup V^m$ over $GF(2^N)$, $N > 1$ of bidimension $n \cup m$ such that $n \leq N$ and $m \leq N$ equipped with the birank bimetric is defined as the birank bispace.*

**DEFINITION 2.4:** *A birank bidistance RD bicode of bilength $n \cup m$ over $GF(2^N)$ is a bisubset of the birank bispace $V^n \cup V^m$ over $GF(2^N)$.*

**DEFINITION 2.5:** *A linear $[n_1, k_1] \cup [n_2, k_2]$ RD bicode is a linear bisubspace of bidimension $k_1 \cup k_2$ in the birank bispace $V^n \cup V^m$. By $C_1[n_1, k_1] \cup C_2[n_2, k_2]$, we denote a linear $[n_1, k_1] \cup [n_2, k_2]$ RD bicode.*

We can equivalently define a RD bicode as follows:

**DEFINITION 2.6:** *Let $V^n$ and $V^m$, $m \neq n$ be rank spaces over $GF(2^N)$, $N > 1$. Suppose $P \subset V^n$ and $Q \subset V^m$ be subsets of the rank spaces over $GF(2^N)$. Then $P \cup Q \subseteq V^n \cup V^m$ is a rank distance bicode of bilength $(n, m)$ over $GF(2^N)$.*

**DEFINITION 2.7:** *Let $C_1[n_1, k_1]$ be $[n_1, k_1]$ RD code (i.e., a linear subspace of dimension $k_1$, in the rank space $V^n$) and $C_2[n_2, k_2]$ be $[n_2, k_2]$ RD code (i.e., a linear subspace of dimension $k_2$ in the rank space $V^m$) ($m \neq n$); then $C_1[n_1, k_1] \cup C_2[n_2, k_2]$ is defined as the linear RD bicode of the linear bisubspace of dimension $(k_1, k_2)$ in the rank bispace $V^n \cup V^m$.*

Now we proceed onto define the notion of the generator bimatrix of a linear $[n_1, k_1] \cup [n_2, k_2]$ RD bicode.

**DEFINITION 2.8:** *A generator bimatrix of a linear $[n_1, k_1] \cup [n_2, k_2]$ RD-bicode $C_1 \cup C_2$ is a $k_1 \times n_1 \cup k_2 \times n_2$ bimatrix over $GF(2^N)$ whose birows form a bibasis for $C_1 \cup C_2$. A generator bimatrix $G = G_1 \cup G_2$ of a linear RD bicode $C_1[n_1, k_1] \cup C_2[n_2, k_2]$ can be brought into the form $G = G_1 \cup G_2 = [I_{k_1}, A_{k_1, n_1-k_1}] \cup [I_{k_2}, A_{k_2 \times n_2 - k_2}]$ where $I_{k_1}$, $I_{k_2}$ is the identity matrix and $A_{k_i \times n_i - k_i}$, $i = 1, 2$ is some matrix over $GF(2^N)$. This form of $G = G_1 \cup G_2$ is called the standard form.*

**DEFINITION 2.9:** *If $G = G_1 \cup G_2$ is a generator bimatrix of $C_1[n_1, k_1] \cup C_2[n_2, k_2]$ then a bimatrix $H = H_1 \cup H_2$ of order $(n_1 - k_1 \times n_1, n_2 - k_2 \times n_2)$ over $GF(2^N)$ such that*
$$GH^T = (G_1 \cup G_2)(H_1 \cup H_2)^T$$
$$= (G_1 \cup G_2)(H_1^T \cup H_2^T)$$
$$= G_1 H_1^T \cup G_2 H_2^T$$
$$= 0 \cup 0$$
*is called a parity check bimatrix of $C_1[n_1, k_1] \cup C_2[n_2, k_2]$.*
*Suppose $C = C_1 \cup C_2$ is a linear $[n_1, k_1] \cup [n_2, k_2]$ RD code with $G = G_1 \cup G_2$ and $H = H_1 \cup H_2$ as its generator and parity check bimatrices respectively, then $C = C_1 \cup C_2$ has two representation,*
*(i) $C = C_1 \cup C_2$ is a row bispace of $G = G_1 \cup G_2$ (i.e., $C_1$ is the row space of $G_1$ and $C_2$ is the row space of $G_2$)*
*(ii) $C = C_1 \cup C_2$ is the solution bispace of $H = H_1 \cup H_2$; i.e., $C_1$ is the solution space of $H_1$ and $C_2$ is the solution space of $H_2$.*

Now we proceed on to define the notion of minimum rank bidistance of a rank distance bicode $C = C_1 \cup C_2$.

**DEFINITION 2.10:** *Let $C = C_1 \cup C_2$ be a rank distance bicode, the minimum-rank bidistance is defined by $d = d_1 \cup d_2$ where,*
$$d_1 = \min\{d_R(x, y) \mid x, y \in C_1, x \neq y\}$$
*and*
$$d_2 = \min\{d_R(x, y) \mid x, y \in C_1, x \neq y\}$$
*i.e.,*
$$d = d_1 \cup d_2$$

$$= \min \{r_1(x) \mid x \in C_1 \text{ and } x \neq 0\} \cup$$
$$\min \{r_2(x) \mid x \in C_2 \text{ and } x \neq 0\}.$$

*If an RD bicode $C = C_1 \cup C_2$ has the minimum rank bidistance $d = d_1 \cup d_2$ then it can correct all bierrors*
$$d = d_1 \cup d_2 \in F_{q^N}^n \cup F_{q^N}^m$$
*with birank*
$$r(e) = (r_1 \cup r_2)(e_1 \cup e_2)$$
$$= r_1(e_1) \cup r_2(e_2) \leq \left\lfloor \frac{d_1 - 1}{2} \right\rfloor \cup \left\lfloor \frac{d_2 - 1}{2} \right\rfloor.$$

*Let $C = C_1 \cup C_2$ denote an $[n_1, k_1] \cup [n_2, k_2]$ RD-bicode over $F_{q^N}$. A generator bimatrix $G = G_1 \cup G_2$ of $C = C_1 \cup C_2$ is a $k_1 \times n_1 \cup k_2 \times n_2$ bimatrix with entries from $F_{q^N}$ whose rows form a bibasis for $C = C_1 \cup C_2$. Then an $(n_1 - k_1) \times n_1 \cup (n_2 - k_2) \times n_2$ bimatrix $H = H_1 \cup H_2$ with entries from $F_{q^N}$ such that*
$$GH^T = (G_1 \cup G_2)(H_1 \cup H_2)^T$$
$$= (G_1 \cup G_2)(H_1^T \cup H_2^T)$$
$$= G_1 H_1^T \cup G_2 H_2^T$$
$$= 0 \cup 0$$

*is called the parity check bimatrix of $C = C_1 \cup C_2$.*

The result analogous to singleton-style bound in case of RD bicode is given in the following:

Result (singleton-style bound) The minimum rank bidistance $d = d_1 \cup d_2$ of any linear $[n_1, k_1] \cup [n_2, k_2]$ RD bicode $C = C_1 \cup C_2 \subseteq F_{q^N}^n \cup F_{q^N}^m$ satisfies the following bounds.
$$d = d_1 \cup d_2 \leq n_1 - k_1 + 1 \cup n_2 - k_2 + 1.$$

Based on this notion we now proceed on to define the new notion of Maximum Rank Distance (MRD) bicodes.

**DEFINITION 2.11:** *An $[n_1, k_1, d_1] \cup [n_2, k_2, d_2]$ RD bicode $C = C_1 \cup C_2$ is called a Maximum Rank Distance (MRD) bicode if the singleton-style bound is reached, that is $d = d_1 \cup d_2 = n_1 - k_1 + 1 \cup n_2 - k_2 + 1$.*

Now we proceed on to briefly give the construction of MRD bicode.

Let $[s] = [s_1] \cup [s_2] = q^{s_1} \cup q^{s_2}$ for any two integers $s_1$ and $s_2$. Let $\{g_1, \ldots, g_n\} \cup \{h_1, h_2, \ldots, h_m\}$ be any set of elements in $F_{q^N}$ that are linearly independent over over $F_q$. A generator bimatrix $G = G_1 \cup G_2$ of an MRD bicode $C = C_1 \cup C_2$ is defined by $G = G_1 \cup G_2$

$$= \begin{bmatrix} g_1 & g_2 & \cdots & g_n \\ g_1^{[1]} & g_2^{[1]} & \cdots & g_n^{[1]} \\ g_1^{[2]} & g_2^{[2]} & \cdots & g_n^{[2]} \\ \vdots & \vdots & & \vdots \\ g_n^{[k_1-1]} & g_n^{[k_1-1]} & \cdots & g_n^{[k_1-1]} \end{bmatrix} \cup \begin{bmatrix} h_1 & h_2 & \cdots & h_m \\ h_1^{[1]} & h_2^{[1]} & \cdots & h_m^{[1]} \\ h_1^{[2]} & h_2^{[2]} & \cdots & h_m^{[2]} \\ \vdots & \vdots & & \vdots \\ h_m^{[k_2-1]} & h_m^{[k_2-1]} & \cdots & h_m^{[k_2-1]} \end{bmatrix}$$

It can be easily proved that the bicode $C = C_1 \cup C_2$ given by the generator bimatrix $G = G_1 \cup G_2$, has the rank bidistance $d = d_1 \cup d_2 = (n_1 - k_1 + 1) \cup (n_2 - k_2 + 1)$. Any bimatrix of the above from is called a Frobenius bimatrix with generating bivector

$$g_C = g_{C_1} \cup h_{C_2}$$
$$= (g_1, \ldots, g_n) \cup (h_1, h_2, \ldots, h_m).$$

One can prove the following theorem:

**THEOREM 2.1:** *Let $C[n, k] = C_1(n_1, k_1) \cup C_2(n_2, k_2)$ be the linear $(n_1, k_1, d_1) \cup (n_2, k_2, d_2)$ MRD bicode with $d_1 = 2t_1 + 1$ and $d_2 = 2t_2 + 1$. Then $C[n, k] = C_1(n_1, k_1) \cup C_2(n_2, k_2)$, bicode corrects all bierrors of birank atmost $t = t_1 \cup t_2$ and detects all bierrors of birank greater than $t = t_1 \cup t_2$.*

Consider the Galois field $GF(2^N)$, $N > 1$. An element $\alpha_1 \cup \beta_2 \in GF(2^N) \cup GF(2^N)$ can be denoted by a biN-tuple $(a_0, \ldots, a_{N-1}) \cup (b_0, b_1, \ldots, b_{N-1})$ as well as by the bipolynomial
$$a_0 + a_1 x + \ldots + a_{N-1} x^{N-1} \cup b_0 + b_1 x + \ldots + b_{N-1} x^{N-1}$$
over $GF(2)$.

We now proceed on to define the new notion of circulant bitranspose.

**DEFINITION 2.12:** *The circulant bitranspose $T_C = T_{C_1}^1 \cup T_{C_2}^2$ of a bivector $\alpha = \alpha_1 \cup \beta_2 = (a_0, \ldots, a_{N-1}) \cup (b_0, b_1, \ldots, b_{N-1}) \in GF(2^N)$ is defined as*
$$\alpha^{T_C} = \alpha_1^{T_{C_1}^1} \cup \beta_2^{T_{C_2}^2} = (a_0, a_1, \ldots, a_{N-1}) \cup (b_0, b_1, \ldots, b_{N-1}).$$
*If $\alpha = \alpha_1 \cup \beta_2 \in GF(2^N) \cup GF(2^N)$ has the bipolynomial representation*
$$a_0 + a_1 x + \ldots + a_{N-1} x^{N-1} \cup b_0 + b_1 x + \ldots + b_{N-1} x^{N-1}$$
*in*
$$\left[ \frac{GF(2)(x)}{\langle x^N + 1 \rangle} \right] \cup \left[ \frac{GF(2)(x)}{\langle x^N + 1 \rangle} \right]$$
*then by; $\alpha_i = \alpha_{1i} \cup \beta_{2i}$ we denote the bivector corresponding to the bipolynomial*

$$[(a_0 + a_1 x + \ldots + a_{N-1} x^{N-1}) \cdot x^i] \bmod (x^n + 1)$$
$$\cup [(b_0 + b_1 x + \ldots + b_{N-1} x^{N-1}) \cdot x^i] \bmod (x^n + 1)$$

*for $i = 0, 1, 2, \ldots, N - 1$. (Note $\alpha_0 = \alpha_1 \cup \beta_2 = \alpha$).*

Now we proceed on to define the biword generated by $\alpha = \alpha_1 \cup \beta_2$.

**DEFINITION 2.13:** *Let $f = f_1 \cup f_2: GF(2^N) \cup GF(2^N) \to [GF(2^N)]^N \cup [GF(2^N)]^N$ be defined as,*
$$f(\alpha) = f_1(\alpha_1) \cup f_2(\beta_2)$$
$$= (\alpha_0^{T_{C_1}^1}, \alpha_1^{T_{C_1}^1}, \ldots, \alpha_{N-1}^{T_{C_1}^1}) \cup (\beta_0^{T_{C_2}^2}, \beta_1^{T_{C_2}^2}, \ldots, \beta_{N-1}^{T_{C_2}^2}).$$

*We call $f(\alpha) = f_1(\alpha_1) \cup f_2(\beta_2)$ as the biword generated by $\alpha = \alpha_1 \cup \beta_2$.*

We analogous to the definition given in MacWilliams and Sloane [61] define circulant bimatrix associated with a bivector in $GF(2^N) \cup GF(2^N)$.

**DEFINITION 2.14:** *A bimatrix of the from*

$$= \begin{bmatrix} a_0 & a_1 & \cdots & a_{N-1} \\ a_{N-1} & a_0 & \cdots & a_{N-2} \\ \vdots & \vdots & & \vdots \\ a_1 & a_2 & \cdots & a_0 \end{bmatrix} \cup \begin{bmatrix} b_0 & b_1 & \cdots & b_{N-1} \\ b_{N-1} & b_0 & \cdots & b_{N-2} \\ \vdots & \vdots & & \vdots \\ b_1 & b_2 & \cdots & b_0 \end{bmatrix}$$

*is called the circulant bimatrix associated with the bivector $(a_0, a_1, ..., a_{N-1}) \cup (b_0, b_1, ..., b_{N-1}) \in GF(2^N) \cup GF(2^N)$. Thus with each $\alpha = \alpha_1 \cup \beta_2 \in GF(2^N) \cup GF(2^N)$, we can associate a circulant bimatrix whose $i^{th}$ bicolumn represents $\alpha_i^{T_{C_1}^1} \cup \beta_i^{T_{C_2}^2}$ ; $i = 0, 1, 2, ..., N - 1$. $f = f_1 \cup f_2$ is nothing but a bimapping of $GF(2^N) \cup GF(2^N)$ on to the pseudo false bialgebra of all $N \times N$ circulant bimatrices over $GF(2)$. Denote the bispace $f(GF(2^N)) = f_1(GF(2^N)) \cup f_2(GF(2^N))$ by $V^N \cup V^N$.*

We define binorm of a biword $v = v_1 \cup v_2 \in V^N \cup V^N$ as follows.

**DEFINITION 2.15:** *The binorm of a biword $v = v_1 \cup v_2 \in V^N \cup V^N$ is defined as the birank of $v = v_1 \cup v_2$ over $GF(2^N)$ (by considering it as a circulant bimatrix over $GF(2)$).*

We denote the binorm of $v = v_1 \cup v_2$ by $r(v) = r_1(v_1) \cup r_2(v_2)$, we prove the following theorem:

**THEOREM 2.2:** *Suppose $\alpha = \alpha_1 \cup \beta_2 \in GF(2^N) \cup GF(2^N)$ has the bipolynomial representation $g_1(x) \cup h_2(x)$ over $GF(2)$ such that gcd $(g_1(x), x^N + 1)$ has degree $N - k_1$ and gcd $(h_2(x), x^N + 1)$* *has degree $N - k_2$ where $1 \leq k_1, k_2 \leq N$; then the binorm of the biword generated by $\alpha = \alpha_1 \cup \beta_2$ is $k_1 \cup k_2$.*

*Proof:* We know the binorm of the biword generated by $\alpha = \alpha_1 \cup \beta_2$ is the birank of the circulant bimatrix

$$= (\alpha_0^{T_{C_1}^1}, \alpha_1^{T_{C_1}^1}, \ldots, \alpha_{N-1}^{T_{C_1}^1}) \cup (\beta_0^{T_{C_2}^2}, \beta_1^{T_{C_2}^2}, \ldots, \beta_{N-1}^{T_{C_2}^2})$$

where $\alpha_i^{T_C} = \alpha_{1i}^{T_{C_1}^1} \cup \beta_{2i}^{T_{C_2}^2}$ represents the bipolynomial $[x^i g_1(x)]$ [mod $x^N + 1$] $\cup$ $[x^i h_2(x)]$ [mod $x^N + 1$] over GF(2). Suppose the bigcd $(g_1(x), x^N+1) \cup (h_2(x), x^N+1)$ has bidegree $N - k_1 \cup N - k_2$, $(0 \leq k_1, k_2 \leq N - 1)$.

To prove that the biword generated by $\alpha = \alpha_1 \cup \beta_2$ has birank $k_1 \cup k_2$. It is enough to prove that the bispace generated by the N-bipolynomials

$$\{g_1(x) \bmod (x^N + 1), (xg_1(x)) [\bmod x^N +1], \ldots,$$
$$[x^{N-1} g_1(x)] \bmod [x^N + 1]\} \cup$$
$$\{h_2(x) \bmod (x^N + 1), (xh_2(x)) [\bmod x^N +1], \ldots,$$
$$[x^{N-1} h_2(x)] \bmod (x^N + 1)\}$$

has bidimension $k_1 \cup k_2$. We will prove that the biset of $k_1 \cup k_2$ bipolynomials

$$\{g_1(x) \bmod (x^N + 1), (xg_1(x)) \bmod (x^N + 1), \ldots,$$
$$(x^{N-1} g_1(x)) \bmod (x^N + 1)\}$$
$$\cup \{h_2(x) \bmod (x^N + 1), (xh_2(x)) (\bmod x^N + 1), \ldots,$$
$$(x^{N-1} h_2(x)) (\bmod x^N + 1)\}$$

forms a bibasis for this bispace.
If possible let
$$a_0(g_1(x)) + a_1(xg_1(x)) + \ldots + a_{k_1-1}(x^{k_1-1} g_1(x)) \cup b_0(h_2(x))$$
$$+ b_1(xh_2(x)) + \ldots + b_{k_2-1}(x^{k_2-1} h_2(x))$$
$$\equiv 0 \cup 0 (\bmod x^N + 1)$$
where $a_i, b_i \in GF(2)$.
This implies $x^{N+1} \cup x^{N+1}$ bidivides
$$(a_0 + a_1 x + \ldots + a_{k_1-1} x^{k_1-1}) g_1(x)$$

$$\cup \ (b_0 + b_1x + \ldots + b_{k_2-2}x^{k_2-1})h_2(x)$$

Now if
$$g_1(x) \cup h_2(x) = p_1(x)\, a_1(x) \cup p_2(x)\, b_2(x)$$

where, $p_1(x) \cup p_2(x)$ is the bigcd $\{(g_1(x), x^N + 1) \cup (h_2(x), x^N + 1)\}$ then $(a_1(x), x^N + 1) \cup (b_2(x), x^N + 1) = 1 \cup 1$. Thus $x^N + 1$ bidivides

$$(a_0 + a_1x + \ldots + a_{k_1-1}x^{k_1-1})g_1(x)$$
$$\cup \ (b_0 + b_1x + \ldots + b_{k_2-2}x^{k_2-1})h_2(x)$$

which inturn implies that the biquotient

$$\left(x^N + 1 \Big/ p_1(x)\right) \cup \left(x^N + 1 \Big/ p_2(x)\right)$$

bidivides

$$(a_0 + a_1x + \ldots + a_{k_1-1}x^{k_1-1})a_1(x)$$
$$\cup \ (b_0 + b_1x + \ldots + b_{k_2-2}x^{k_2-1})b_2(x)$$

That is

$$\left(x^N + 1 \Big/ p_1(x)\right) \cup \left(x^N + 1 \Big/ p_2(x)\right)$$

bidivides
$$(a_0 + a_1x + \ldots + a_{k_1-1}x^{k_1-1}) \cup (b_0 + b_1x + \ldots + b_{k_2-2}x^{k_2-1})$$

which is a contradiction as

$$\left(x^N + 1 \Big/ p_1(x)\right) \cup \left(x^N + 1 \Big/ p_2(x)\right)$$

has bidegree $k_1 \cup k_2$ where as the bipolynomial $(a_0 + a_1x + \ldots + a_{k_1-1}x^{k_1-1}) \cup (b_0 + b_1x + \ldots + b_{k_2-2}x^{k_2-1})$ has bidegree atmost $k_1 -1 \cup k_2 -1$. Hence the bipolynomials

$$\{g_1(x) \bmod (x^N + 1), xg_1(x) \bmod (x^N + 1), \ldots,$$

$$x^{k_1-1} g_1(x) \bmod(x^N + 1)\} \cup$$
$$\{h_2(x) \bmod (x^N + 1), xh_2(x) \bmod (x^N + 1), \ldots,$$
$$x^{k_2-1} h_2(x) \bmod (x^N + 1)\}$$

are bilinearly independent over $GF(2) \cup GF(2)$, i.e.;

$$\{g_1(x) \bmod (x^N + 1), xg_1(x) \bmod (x^N + 1), \ldots,$$
$$x^{k_1-1} g_1(x) \bmod(x^N + 1)\}$$
$$\cup \{h_2(x) \bmod (x^N + 1), xh_2(x) \bmod (x^N + 1), \ldots,$$
$$x^{k_2-1} h_2(x) \bmod (x^N + 1)\}$$

bigenerate the bispace. For this it is enough to prove that $x^i g_1(x) \cup x^i h_2(x)$ is a linear bicombination of these bipolynomials for $k_1 \leq i \leq N-1$ and $k_2 \leq i \leq N-1$. Let $x^N + 1 \cup x^N + 1 = p_1(x)q_1(x) \cup p_2(x)q_2(x)$ where

$$q_1(x) \cup q_2(x) =$$
$$c_0 + c_1 x + \ldots + c_{k_i} x^{k_1} \cup d_0 + d_1 x + \ldots + d_{k_2} x^{k_2}$$

(Note: $c_0 = c_{k_1} = 1$ and $d_0 = d_{k2} = 1$, since $q_1(x) \cup q_2(x)$ bidivides $x^N + 1 \cup x^N + 1$, i.e., $g_1(x)$ divides $x^N + 1$ and $g_2(x)$ divides $x^N + 1$). Also we have

$$g_1(x) = p_1(x) a_1(x)$$

and

$$h_2(x) = p_2(x) b_2(x)$$

i.e.,

$$g_1(x) \cup h_2(x) = p_1(x) a_1(x) \cup p_2(x) b_2(x).$$

Thus

$$x^N + 1 \cup x^N + 1 = \left(g_1(x) \cdot \frac{p_1(x)}{a_1(x)}\right) \cup \left(h_2(x) \cdot \frac{p_2(x)}{b_2(x)}\right)$$

that is

$$g_1(x) \cdot \frac{p_1(x)}{a_1(x)} \equiv 0 \bmod(x^N + 1)$$

and

$$h_2(x) \cdot \frac{p_2(x)}{b_2(x)} \equiv 0 \bmod(x^N + 1)$$

$$x^N + 1 \cup x^N + 1 = \left(g_1(x) \cdot \frac{p_1(x)}{a_1(x)}\right) \cup \left(h_2(x) \cdot \frac{p_2(x)}{b_2(x)}\right).$$

Now

$$\frac{g_1(x)(a_0 + a_1 x + \ldots + a_{k_1-1} x^{k_1-1})}{a_1(x)}$$

$$= \left[\frac{g_1(x) x^{k_1}}{a_1(x)}\right] \mod(x^N + 1) \text{ (since } a_k = 1\text{)}.$$

That is

$$x^{k_1} g_1(x) = (a_0 + a_1 x + \ldots + a_{k_1-1} x^{k_1-1} g_1(x)) \mod(x^N + 1).$$

Hence,

$$x^{k_1} g_1(x) = (a_0 g_1(x) + a_1(x g_1(x)) + \ldots +$$
$$a_{k_1-1}[x^{k_1-1} g_1(x)]) \mod(x^N + 1);$$

a linear combination of

$$g_1(x) \mod(x^N + 1), [x g_1(x)] \mod(x^N + 1), \ldots,$$
$$[x^{k_1-1} g_1(x)] \mod(x^N + 1)$$

over GF(2).

Now it can be easily proved that $x^i g_1(x)$ is a linear combination of

$$g_1(x) \mod(x^N + 1), x g_1(x) \mod(x^N + 1), \ldots,$$
$$x^{k_1-1} g_1(x) \mod(x^N + 1)$$

for $i > k_1$.

Similar argument holds good for $x^i h_2(x)$ where $i > k_2$. Hence the bispace generated by the bipolynomial $\{g_1(x) \mod (x^N + 1), x g_1(x) \mod (x^N + 1), \ldots, x^{k_1-1} g_1(x) \mod(x^N + 1)\} \cup \cup \{h_2(x) \mod (x^N + 1), x h_2(x) \mod (x^N + 1), \ldots, x^{k_2-1} h_2(x) \mod (x^N + 1)\}$ has bidimension $k_1 \cup k_2$; i.e., the birank of the biword generated by $\alpha = \alpha_1 \cup \beta_2$ is $k_1 \cup k_2$.

Now we have the two corollaries to be true.

**COROLLARY 2.1:** *If $\alpha = \alpha_1 \cup \beta_2 \in GF(2^N) \cup GF(2^N)$ is such that its bipolynomial representation $g_1(x) \cup h_2(x)$ is relatively prime to $x^N+1 \cup x^N+1$ (i.e. $g_1(x)$ is relatively prime to $x^N+1$ and* $h_2(x)$ is relatively prime to $x^N + 1$) then the binorm of the biword generated by $\alpha = \alpha_1 \cup \beta_2$ is $(N, N)$ and $f(\alpha) = f_1(\alpha_1) \cup f_2(\beta_2)$ is invertible.

*Proof:* Follows from the theorem $\text{bigcd}(g_1(x) \cup h_2(x), x^N+1 \cup x^N + 1)$
$= \text{bigcd}(g_1(x), x^N + 1) \cup (h_2(x), x^N + 1)$
$= 1 \cup 1$
has bidegree $0 \cup 0$ and hence birank of $f(\alpha) = f_1(\alpha_1) \cup f_2(\beta_2)$ is $(N, N)$.

**COROLLARY 2.2:** *The binorms of the bivectors corresponding to the bipolynomials $x + 1 \cup x + 1$ and $x^{N-1} + x^{N-2} + \ldots + x + 1 \cup x^{N-1} + x^{N-2} + \ldots + x + 1$ are respectively $(N - 1, N - 1)$ and $(1, 1)$.*

Now we proceed on to define the bidistance bifunction on $V^N \cup V^N$.

**DEFINITION 2.16:** *The bidistance between two biwords $u = u_1 \cup u_2$ and $v = v_1 \cup v_2$ in $V^N \cup V^N$ is defined as*
$$d(u, v) = d_1(u_1, v_1) \cup d_2(u_2, v_2)$$
$$= r(u + v)$$
$$= r_1(u_1 + v_1) \cup r_2(u_2 + v_2).$$

Now we proceed on to define the new notion of circulant rank bicodes of bilength $N = N_1 \cup N_2$.

**DEFINITION 2.17:** *Let $C_1$ be a circulant rank code of length $N_1$ which is a subspace of $V^{N_1}$ equipped with the distance function $d_1(u_1, v_1) = r_1(u_1, v_1)$ and $C_2$ be a circulant rank code of length $N_2$ which is the subspace of $V^{N_2}$ equipped with distance function $d_2(u_2, v_2) = r_2(u_2, v_2)$ where $V^{N_1}$ and $V^{N_2}$ are spaces defined over $GF(2^N)$ with $N_1 \neq N_2$. $C = C_1 \cup C_2$ is defined as the circulant birank bicode of bilength $N = N_1 \cup N_2$ defined as a bisubsapce of $V^{N_1} \cup V^{N_2}$ equipped with the bidistance bifunction*
$$d_1(u_1, v_1) \cup d_2(u_2, v_2) = r_1(u_1 + v_1) \cup r_2(u_2 + v_2).$$

**DEFINITION 2.18:** *A circulant birank bicode of bilength $N = N_1 \cup N_2$ is called bicyclic if whenever $(v_1^1...v_{N_1}^1) \cup (u_1^2...u_{N_2}^2)$ is a bicodeword then it implies $(v_2^1 v_3^1 ... v_{N_1}^1 v_1^1) \cup (u_2^2 u_3^2 ... u_{N_2}^2 u_1^2)$ is also a bicodeword.*

Now we proceed on to define semi MRD bicode.

**DEFINITION 2.19:** *Let $C_1[n_1, k_1]$ be a $[n_1, k_1]$ RD code and $C_2[n_2, k_2, d_2]$ be a MRD code. $C_1[n_1, k_1] \cup C_2[n_2, k_2, d_2]$ is defined as the semi MRD bicode if $n_1 \neq n_2$ and $k_1 \neq k_2$ and $C_1[n_1, k_1]$ is only a RD-code and not a MRD code.*

These bicodes will be useful in application where one set of code never attains the singleton-style bound were as the other code attains the singleton style bound.

The special nature of these codes will be very useful in applications were two types of RD codes are used simultaneously. Next we proceed on to define the notion of semi circulant rank bicode of type I and type II.

**DEFINITION 2.20:** *Let $C_1 = C_1[n_1, k_1]$ be a RD-code and $C_2$ be a circulant rank code both are subspaces of rank spaces defined over the same field $GF(2^N)$. Then $C_1 \cup C_2$ is defined as the semi circulant rank bicode of type I.*

These codes find their application where one RD-code which does not attain its single style bound and another need is a circulant rank code. Now we proceed on to define the new concept of semi circulant rank bicode of type II.

**DEFINITION 2.21:** *Let $[n_1, k_1, d_1] = C_1$ be a MRD code and $C_2$ be a circulant rank code both are subspaces or rank spaces defined over the same field $GF(2^N)$ or $F_{q^N}$. $C_1 \cup C_2$ is defined to be semi circulant bicode of type II.*

We see in semi circulant rank bicodes of type I, we use only RD codes which are not MRD codes and in semi circulant rank bicodes of type II, we use only MRD codes which are never RD codes. These rank bicodes will find their applications in special situations. Next we proceed on to define semicyclic circulant rank bicode of type I and type II.

**DEFINITION 2.22:** *Let $C_1[n_1, k_1]$ be a RD-code which is not a MRD code and $C_2$ be a cyclic circulant rank code, both take entries from the same field $GF(2^N)$. $C_1 \cup C_2$ is defined to be a semicyclic circulant rank bicode of type I.*

These also can be used when simultaneous use of two different types of rank codes are needed. Next we proceed on to define the notion of semicyclic circulant rank bicode of type II.

**DEFINITION 2.23:** *Let $[n_1, k_1, d_1] = C_1$ be a MRD-code which is a subspace of $V^N$, $V^N$ defined over $GF(2^N)$. $C_2$ be a cyclic circulant rank code with entries from $GF(2^N)$. $C_1 \cup C_2$ is defined to be the semicyclic circulant rank bicode.*

These bicodes also find their applications when two types of codes are needed simultaneously. Next we proceed on to define semicyclic circulant rank bicode.

**DEFINITION 2.24:** *Let $C_1$ be a circulant rank code and $C_2$ be a cyclic circulant rank code, $C_1 \cup C_2$ is defined as the semicyclic circulant rank bicode, both $C_1$ and $C_2$ take entries from the same field.*

The special semicyclic circulant rank bicodes can be used in communication bichannel having very high error probability for error correction.

Now we proceed on to define yet another class of rank bicodes.

**DEFINITION 2.25:** *Let $C_1[n_1, k_1]$ and $C_2[n_2, k_2]$ be any two distinct Almost Maximum Rank Distance (AMRD) codes with the minimum distances greater than or equal to $n_1 - k_1$ and $n_2 -$* $k_2$ respectively defined over $GF(2^N)$. Then $C_1[n_1, k_1] \cup C_2[n_2, k_2]$ is defined as the Almost Maximum Rank Distance bicode (AMRD-bicode) over $GF(2^N)$.

An AMRD bicode whose minimum bidistance is greater than $n_1 - k_1 \cup n_2 - k_2$ is an MRD bicode and hence the class of MRD bicodes is a subclass of the class of AMRD bicode.

We have an interesting property about these AMRD bicodes.

**THEOREM 2.3:** *When $n_1 - k_1 \cup n_2 - k_2$ is an odd pair of biintegers (i.e., $n_1 - k_1$ and $n_2 - k_2$ are odd integers);*

  (i) *The error correcting capability of the $[n_1, k_1] \cup [n_2, k_2]$ AMRD bicode is equal to that of an $[n_1, k_1] \cup [n_2, k_2]$ MRD bicode.*
  (ii) *An $[n_1, k_1] \cup [n_2, k_2]$ AMRD bicode is better than any $[n_1, k_1] \cup [n_2, k_2]$ bicode in Hamming metric for error correction.*

*Proof:* (i) Suppose $C = C_1 \cup C_2$ is a $[n_1, k_1] \cup [n_2, k_2]$ AMRD bicode such that $n_1 - k_1 \cup n_2 - k_2$ is an odd biinteger (i.e., $n_1 - k_1 \neq n_2 - k_2$ are odd integers). The maximum number of bierrors corrected by $C = C_1 \cup C_2$ is given by

$$\frac{(n_1 - k_1 - 1)}{2} \cup \frac{(n_2 - k_2 - 1)}{2}.$$

But

$$\frac{(n_1 - k_1 - 1)}{2} \cup \frac{(n_2 - k_2 - 1)}{2}$$

is equal to the error correcting capability of an $[n_1, k_1] \cup [n_2, k_2]$ MRD bicode (Since $n_1 - k_1$ and $n_2 - k_2$ are odd). Thus when $n_1 - k_1 \cup n_2 - k_2$ is biodd (i.e., both $n_1 - k_1$ and $n_2 - k_2$ are odd) the $[n_1, k_1] \cup [n_2, k_2]$ AMRD bicode is as good as an $[n_1, k_1] \cup [n_2, k_2]$ MRD bicode.

(ii) Suppose $C = C_1 \cup C_2$ is a $[n_1, k_1] \cup [n_2, k_2]$ AMRD bicode such that $n_1 - k_1$ and $n_2 - k_2$ are odd. Then, each bicodeword of C can correct $L_{r_1}(n_1) \cup L_{r_2}(n_2) = L_r(n)$ error bivectors where

$$r = r_1 \cup r_2 = \frac{(n_1 - k_1 - 1)}{2} \cup \frac{(n_2 - k_2 - 1)}{2}$$

and

$$L_r(n) = L_{r_1}(n_1) \cup L_{r_2}(n_2)$$
$$= 1 + \sum_{i=1}^{n_1} \begin{bmatrix} n_1 \\ i \end{bmatrix}(2^N - 1)\ldots(2^N - 2^{i-1}) \cup$$
$$1 + \sum_{i=1}^{n_2} \begin{bmatrix} n_2 \\ i \end{bmatrix}(2^N - 1)\ldots(2^N - 2^{i-1}).$$

Consider the same $[n_1, k_1] \cup [n_2, k_2]$ bicode in Hamming metric. Let it be $D_1 \cup D_2 = D$, then the minimum bidistance of D is atmost $(n_1 - k_1 + 1) \cup (n_2 - k_2 + 1)$. The error correcting capability of D is

$$\frac{(n_1 - k_1 + 1 - 1)}{2} \cup \frac{(n_2 - k_2 + 1 - 1)}{2} = r_1 \cup r_2$$

($[n_1 - k_1]$ and $[n_2 - k_2]$ are odd).

Hence the number of error bivectors corrected by a codeword is given by

$$\sum_{i=0}^{r_1} \begin{bmatrix} n_1 \\ i \end{bmatrix}(2^N - 1)^i \cup \sum_{i=0}^{r_2} \begin{bmatrix} n_2 \\ i \end{bmatrix}(2^N - 1)^i$$

which is clearly less than $L_{r_1}(n_1) \cup L_{r_2}(n_2)$.

Thus the number of error bivectors that can be corrected by the $[n_1, k_1] \cup [n_2, k_2]$ AMRD bicode is much greater than that of the same bicode considered in Hamming metric.

For a given bilength $n = n_1 \cup n_2$, a single error correcting AMRD bicode is one having bidimension $n_1 - 3 \cup n_2 - 3$ and the minimum distance greater than or equal to $3 \cup 3$. We now proceed on to give a characterization of a single error correcting AMRD bicode in terms of its parity check bimatrices. The characterization is based on the condition for the minimum distance proved by Gabidulin in [24, 27].

**THEOREM 2.4:** Let $H = H_1 \cup H_2 = (\alpha_{ij}^1) \cup (\alpha_{ij}^2)$ be a $3 \times n_1 \cup 3 \times n_2$ bimatrix of birank 3 over $GF(2^N)$; $n_1 \leq N$ and $n_2 \leq N$ which satisfies the following condition. For any two distinct, non empty bisubsets $P_1, P_2$ where $P_1 = P_1^1 \cup P_2^1$ and $P_2 = P_1^2 \cup P_2^2$ of $\{1, 2, ..., n_1\}$ and $\{1, 2, 3, ..., n_2\}$ respectively there exists

$$i_1 = i_1^1 \cup i_2^1, \ i_2 = i_1^2 \cup i_2^2 \in \{1, 2, 3\} \cup \{1, 2, 3\}$$

such that,

$$\left( \sum_{j_1^1 \in P_1^1} \alpha_{i_1^1 j_1^1}^1 \cdot \sum_{k_1^1 \in P_2^1} \alpha_{i_2^1 k_1^1}^1 \right) \cup \left( \sum_{j_1^2 \in P_1^2} \alpha_{i_1^2 j_1^2}^2 \cdot \sum_{k_2^2 \in P_1^2} \alpha_{i_2^2 k_2^2}^2 \right) \neq$$

$$\left( \sum_{j_1^1 \in P_1^1} \alpha_{i_2^1 j_1^1}^1 \cdot \sum_{k_1^1 \in P_2^1} \alpha_{i_1^1 k_1^1}^1 \right) \cup \left( \sum_{j_1^2 \in P_1^2} \alpha_{i_2^2 j_1^2}^2 \cdot \sum_{k_2^2 \in P_1^2} \alpha_{i_1^2 k_2^2}^2 \right).$$

then, $H = H_1 \cup H_2$ as a parity check bimatrix defines a $[n_1, n_1 - 3] \cup [n_2, n_2 - 3]$ single bierror correcting AMRD bicode over $GF(2^N)$.

*Proof:* Given $H = H_1 \cup H_2$ is a $3 \times n_1 \cup 3 \times n_2$ bimatrix of birank $3 \cup 3$ over $GF(2^N)$, so $H = H_1 \cup H_2$ as a parity check bimatrix defines a $[n_1, n_1 - 3] \cup [n_2, n_2 - 3]$ RD bicode $C = C_1 \cup C_2$ over $GF(2^N)$ where,

$$C_1 = \{x \in V^{n_1} \mid xH_1^T = 0\}$$

and

$$C_2 = \{x \in V^{n_2} \mid xH_2^T = 0\}.$$

It remains to prove that the minimum bidistance of $C = C_1 \cup C_2$ is greater than or equal to $3 \cup 3$. We will prove that no non zero bicodeword of $C = C_1 \cup C_2$ has birank less than $3 \cup 3$. The proof is by method of contradiction.

Suppose there exists a non zero bicodeword $x = x_1 \cup x_2$ such that $r_1(x_1) \leq 2$ and $r_2(x_2) \leq 2$, then $x = x_1 \cup x_2$ can be written as $x = x_1 \cup x_2 = (y_1 \cup y_2)(M_1 \cup M_2)$ where $y_1 = (y_1^1, y_2^1)$ and $y_2 = (y_1^2, y_2^2)$; $y_1^1, y_2^1, y_1^2, y_2^2 \in GF(2^N)$ and $M = M_1 \cup M_2 = (m_{ij}^1) \cup (m_{ij}^2)$ is a $2 \times n_1 \cup 2 \times n_2$ bimatrix of birank $2 \cup 2$ over $GF(2)$. Thus

$$(yM)H^T = y_1 M_1 H_1^T \cup y_2 M_2 H_2^T = 0 \cup 0$$

implies that

$$y(MH^T) = y_1(M_1 H_1^T) \cup y_2(M_2 H_2^T) = 0 \cup 0.$$

Since $y = y_1 \cup y_2$ is non zero; $y(MH^T) = 0 \cup 0$; implies $y_1(M_1 H_1^T) = 0$ and $y_2(M_2 H_2^T) = 0$ that is the $2 \times 3$ bimatrix $M_1 H_1^T \cup M_2 H_2^T$ has birank less than 2 over $GF(2^N)$.

Now let

$$P_1 = P_1^1 \cup P_2^1 = \{j_1^1 \text{ such that } m_{1j_1^1}^1 = 1\} \cup \{j_2^2 \text{ such that } m_{1j_2^2}^2 = 1\}$$

and

$$P_2 = P_1^2 \cup P_2^2 = \{j_1^1 \text{ such that } m_{2j_1^1}^1 = 1\} \cup \{j_2^2 \text{ such that } m_{2j_2^2}^2 = 1\}.$$

Since $M = M_1 \cup M_2 = (m_{ij}^1) \cup (m_{ij}^2)$ is a $2 \times 2$ bimatrix of birank $2 \cup 2$. $P_1$ and $P_2$ are disjoint non empty bisubsets of $\{1, 2, \ldots, n_1\} \cup \{1, 2, \ldots, n_2\}$ respectively and

$$M_1 H^T = M_1 H_1^T \cup M_2 H_2^T$$

$$= \begin{pmatrix} \sum_{j_1^1 \in P_1^1} \alpha_{1j_1^1}^1 & \sum_{j_1^1 \in P_1^1} \alpha_{2j_1^1}^1 & \sum_{j_1^1 \in P_1^1} \alpha_{3j_1^1}^1 \\ \sum_{j_1^1 \in P_2^1} \alpha_{1j_1^1}^1 & \sum_{j_1^1 \in P_2^1} \alpha_{2j_1^1}^1 & \sum_{j_1^1 \in P_2^1} \alpha_{3j_1^1}^1 \end{pmatrix} \cup$$

$$\begin{pmatrix} \sum_{j_2^2 \in P_1^2} \alpha_{1j_2^2}^2 & \sum_{j_2^2 \in P_1^2} \alpha_{2j_2^2}^2 & \sum_{j_2^2 \in P_1^2} \alpha_{3j_2^2}^2 \\ \sum_{j_2^2 \in P_2^2} \alpha_{1j_2^2}^2 & \sum_{j_2^2 \in P_2^2} \alpha_{2j_2^2}^2 & \sum_{j_2^2 \in P_2^2} \alpha_{3j_2^2}^2 \end{pmatrix}.$$

But the selection of $H = H_1 \cup H_2$ is such that their exists $i_1^1, i_2^1, i_2^2, i_1^2 \in \{1, 2, 3\}$ such that

$$\sum_{j_1^1 \in P_1^1} \alpha^1_{i_1^1 j_1^1} \cdot \sum_{k_1 \in P_2^1} \alpha^1_{i_2^1 k_1} \cup \sum_{j_2^2 \in P_1^2} \alpha^2_{i_1^2 j_2^2} \cdot \sum_{k_2 \in P_2^2} \alpha^2_{i_2^2 k_2}$$

$$\neq \sum_{j_1^1 \in P_1^1} \alpha^1_{i_2^1 j_1^1} \cdot \sum_{k_1 \in P_2^1} \alpha^1_{i_1^1 k_1} \cup \sum_{j_2^2 \in P_1^2} \alpha^2_{i_2^2 j_2^2} \cdot \sum_{k_2 \in P_2^2} \alpha^2_{i_1^2 k_2}.$$

Hence in $MH^T = M_1 H_1^T \cup M_2 H_2^T$, there exists a $2 \times 2$ subbimatrix whose determinant is non zero;
i.e.,
$$r(MH^T) = r_1(M_1 H_1^T) \cup r_2(M_2 H_2^T)$$
$$= 2 \cup 2$$
over GF(2); this contradicts the fact that birank of
$$MH^T = M_1 H_1^T \cup M_2 H_2^T < 2 \cup 2.$$

Hence the result.

Now using constant rank code, we proceed on to define the notion of constant rank bicodes of bilength $n_1 \cup n_2$.

**DEFINITION 2.26:** *Let $C = C_1 \cup C_2$ be a rank bicode, where $C_1$ is a constant rank code of length $n_1$ (a subset of rank space $V^{n_1}$) and $C_2$ is a constant rank code of length $n_2$ (a subset of rank space $V^{n_2}$) then C is a constant rank bicode of bilength $n_1 \cup n_2$; that is every bicodeword has same birank.*

**DEFINITION 2.27:** *$A(n_1, r_1, d_1) \cup A(n_2, r_2, d_2)$ is defined as the maximum number of bivectors in $V^{n_1} \cup V^{n_2}$, constant birank $r_1 \cup r_2$ and bidistance between any two bivectors is atleast $d_1 \cup d_2$.*

*(By a $(n_1, r_1, d_1) \cup (n_2, r_2, d_2)$ biset we mean a bisubset of bivectors of $V^{n_1} \cup V^{n_2}$ having constant birank $r_1 \cup r_2$ and bidistance between any two bivectors is atleast $d_1 \cup d_2$).*

We analyze the bifunction $A(n_1, r_1, d_1) \cup A(n_2, r_2, d_2)$ by the following theorem.

**THEOREM 2.5:**
(i) $A(n_1, r_1, d_1) \cup A(n_2, r_2, d_2) = L_{r_1}(n_1) \cup L_{r_2}(n_2)$, *the number of bivectors of birank $r_1 \cup r_2$ in $V^{n_1} \cup V^{n_2}$.*
(ii) $A(n_1, r_1, d_1) \cup A(n_2, r_2, d_2) = 0 \cup 0$ *if $r_1 > 0$ and $r_2 > 0$ or $d_1 > n_1$ and $d_2 > n_2$ or $d_1 > 2r_1$ and $d_2 > 2r_2$.*

*Proof:* (i) Follows from the fact that $L_{r_1}(n_1) \cup L_{r_2}(n_2)$ is the number of bivectors of bilength $n_1 \cup n_2$, constant birank $r_1 \cup r_2$ and bidistance between any two distinct bivectors in a rank bispace $V^{n_1} \cup V^{n_2}$, is always greater than or equal to $1 \cup 1$.
(ii) Follows immediately from the definition of $A(n_1, r_1, d_1) \cup A(n_2, r_2, d_2)$.

**THEOREM 2.6:** $A(n_1, 1, 2) \cup A(n_2, 1, 2) = 2^{n_1} - 1 \cup 2^{n_2} - 1$ *over any Galois field $GF(2^N)$.*

*Proof:* Denote by $V_1 \cup V_2$ the set of bivectors of birank $1 \cup 1$ in $V^{n_1} \cup V^{n_2}$. We know for each non zero element $\alpha_1 \cup \alpha_2 \in GF(2^N)$ there exists $(2^{n_1} - 1) \cup (2^{n_2} - 1)$ bivectors of birank $1 \cup 1$ having $\alpha_1 \cup \alpha_2$ as a coordinate. Thus the cardinality of $V_1 \cup V_2$ is $(2^N - 1)(2^{n_1} - 1) \cup (2^N - 1)(2^{n_2} - 1)$. Now bidivide $V_1 \cup V_2$ into $(2^{n_1} - 1) \cup (2^{n_2} - 1)$ blocks of $(2^N - 1) \cup (2^N - 1)$ bivectors such that each block consists of the same pattern of all nonzero bielements of $GF(2^N) \cup GF(2^N)$.

Then from each biblock almost one bivector can be chosen such that the selected bivectors are atleast rank 2 apart from each other. Such a biset we call as a $(n_1, 1, 2) \cup (n_2, 1, 2)$ biset. Also it is always possible to construct such a biset. Thus $A(n_1, 1, 2) \cup A(n_2, 1, 2) = 2^{n_1} - 1 \cup 2^{n_2} - 1$.

**THEOREM 2.7:** $A(n_1, n_1, n_1) \cup A(n_2, n_2, n_2) = 2^N - 1 \cup 2^N - 1$ (i.e. $A(n_i, n_i, n_i) = 2^N - 1$; $i = 1, 2$); over any $GF(2^N)$.

*Proof:* Denote by $V_{n_1} \cup V_{n_2}$; the biset of all bivectors of birank $n_1 \cup n_2$ in the bispace $V^{n_1} \cup V^{n_2}$. We know the bicardinality of $V_{n_1} \cup V_{n_2}$ is $(2^N - 1)(2^N - 2) \ldots (2^N - 2^{n_1-1}) \cup (2^N - 1)(2^N - 2) \ldots (2^N - 2^{n_2-1})$, by the definition of a $(n_1, n_1, n_1) \cup (n_2, n_2, n_2)$ biset the bidistance between any two bivectors should be $n_1 \cup n_2$. Thus no two bivectors can have a common symbol at a coordinate place $i_1 \cup i_2$; $(1 \leq i_1 \leq n_1, 1 \leq i_2 \leq n_2)$. This implies that $A(n_1, n_1, n_1) \cup A(n_2, n_2, n_2) \leq 2^N - 1 \cup 2^N - 1$.

Now we construct a $(n_1, n_1, n_1) \cup (n_2, n_2, n_2)$ biset as follows. Select N bivectors from $V_{n_1} \cup V_{n_2}$ such that

1. Each bibasis bielement of $GF(2^N) \cup GF(2^N)$ should occur (can be as a bicombination) atleast once in each bivector.
2. If the $(i_1^{th}, i_2^{th})$ bivector is chosen $[(i_1 + 1)^{th}, (i_2 + 1)^{th}]$ bivector should be selected such that its birank bidistance from any bilinear combination of the previous $(i_1, i_2)$ bivectors is $n_1 \cup n_2$. Now the set of all bilinear combination of these $N \cup N$ bivectors over $GF(2) \cup GF(2)$, will be such that the bidistance between any two bivectors is $n_1 \cup n_2$. Hence it is a $(n_1, n_1, n_1) \cup (n_2, n_2, n_2)$ biset. Also the bicardinality of this $(n_1, n_1, n_1) \cup (n_2, n_2, n_2)$ biset is $2^N - 1 \cup 2^N - 1$ (we do not count all zero bilinear combination); thus $A(n_1, n_1, n_1) \cup A(n_2, n_2, n_2) = 2^N - 1 \cup 2^N - 1$.

Recall a [n, 1] repetition RD code is code generated by the matrix $G = (1\ 1\ 1\ \ldots\ 1)$ over $F_{q^N}$. Any non zero codeword has rank 1.

**DEFINITION 2.28:** *A $[n_1, 1] \cup [n_2, 1]$ repetition RD bicode is a bicode generated by the bimatrix $G = G_1 \cup G_2 = (11 \ldots 1) \cup (1\ 1\ \ldots\ 1)$ $(G_1 \neq G_2)$ over $F_{2^N}$. Any non zero bicodeword has birank $1 \cup 1$.*

We proceed onto define the notion of covering biradius.

**DEFINITION 2.29:** *Let $C = C_1 \cup C_2$ be a linear $[n_1, k_1] \cup [n_2, k_2]$ RD bicode defined over $F_{2^N}$. The covering biradius of $C = C_1 \cup C_2$ is defined as the smallest pair of integers $(r_1, r_2)$ such that all bivectors in the rank bispace $F_{2^N}^{n_1} \cup F_{2^N}^{n_2}$ are within the rank bidistance $r_1 \cup r_2$ of some bicodeword. The covering biradius of $C = C_1 \cup C_2$ is denoted by $t(C_1) \cup t(C_2)$. In notation, $t(C) = t(C_1) \cup t(C_2)$*

$$= \max_{x_1 \in F_{2^N}^{n_1}} \left\{ \min_{c_1 \in C_1} \{r_1(x_1 + c_1)\} \right\}$$
$$\cup \max_{x_2 \in F_{2^N}^{n_2}} \left\{ \min_{c_2 \in C_2} \{r_2(x_2 + c_2)\} \right\}.$$

**THEOREM 2.8:** *The linear $[n_1, k_1] \cup [n_2, k_2]$ RD bicode $C = C_1 \cup C_2$ satisfies $t(C) = t(C_1) \cup t(C_2) \leq n_1 - k_1 \cup n_2 - k_2$.*

Proof is direct.

**THEOREM 2.9:** *The covering biradius of a $[n_1, 1] \cup [n_2, 1]$ repetition RD-bicode over $F_{2^N}$ is $[n_1 - 1] \cup [n_2 - 1]$.*

Direct from theorem 2.8 as $k_1 = k_2 = 1$. Next we proceed on to define the Cartesian biproduct of two linear RD-bicodes.

The Cartesian biproduct of two linear RD-bicodes $C = C_1[n_1^1, k_1^1] \cup C_2[n_2^1, k_2^1]$ and $D = D_1[n_1^2, k_1^2] \cup D_2[n_2^2, k_2^2]$ over $F_{2^N}$ is given by

$$C \times D = C_1 \times D_1 \cup C_2 \times D_2.$$
$$= \{(a_1^1, b_1^1) \mid a_1^1 \in C_1, b_1^1 \in D_1\} \cup \{(a_1^2, b_1^2) \mid a_1^2 \in C_2, b_1^2 \in D_2\}.$$

$C \times D$ is a $\{(n_1^1 + n_1^2) \cup (n_2^1 + n_2^2), (k_1^1 + k_1^2) \cup (k_2^1 + k_2^2)\}$ linear RD bicode (We assume $= (n_1^1 + n_1^2) \leq N$ and $(n_2^1 + n_2^2) \leq N$).

Now the reader is left with the task of proving the following theorem:

**THEOREM 2.10:** *If $C = C_1 \cup C_2$ and $D = D_1 \cup D_2$ be two linear RD bicodes then $t(C \times D) \leq (t(C_1) + t(D_1)) \cup t(C_2) + t(D_2))$.*

Hint for the proof. If $C = C_1 \cup C_2$ and $D = D_1 \cup D_2$ then $C \times D = \{C_1 \times D_1\} \cup \{C_2 \times D_2\}$ and
$$t(C \times D) = (C_1 \times D_1) \cup (C_2 \times D_2)$$
$$\leq \{t(C_1) + t(D_1)\} \cup \{t(C_2) + t(D_2)\}.$$
.

Now we proceed on to define notion of bidivisble linear RD bicodes.

We just recall the definition of divisible linear RD codes. Let $C(n, k, d)$ be a linear RD code over $F_{q^N}$, $n \leq N$ and $N > 1$ with length n, dimension k and minimum distance d. If there exists m > 1 an integer such that $m/r(c;q)$ for all $0 \neq c \in C$ then the code C is defined to be divisible. ($r(x; q)$ denotes the rank norm of x over the field $F_q$).

**DEFINITION 2.30:** *Let $C = C_1(n_1, k_1, d_1) \cup C_2(n_2, k_2, d_2)$ be a linear RD bicode over $F_{q^N}$, $n_1 \leq N$, $n_2 \leq N$ and $N > 1$. If there exists $(m_1, m_2)$ $(m_1 > 1$ and $m_2 > 1)$ such that*
$$\frac{m_1}{r(c_1;q)} \text{ and } \frac{m_2}{r(c_2;q)}$$
*for all $c_1 \in C_1$ and for all $c_2 \in C_2$ then we say the bicode C is bidivisible.*

**THEOREM 2.11:** *Let $C = [n_1, 1_1, n_1] \cup [n_2, 1_2, n_2]$ $(n_1 \neq n_2)$ be a MRD-bicode for all $n_1 \leq N$ and $n_2 \leq N$; C is a bidivisible bicode.*

*Proof:* Since there cannot exist bicodewords of birank greater than $(n_1, n_2)$ in an $[n_1, 1, n_1] \cup [n_2, 1, n_2]$ MRD-bicode.

**DEFINITION 2.31:** *Let $C_1 = [n_1, k_1]$ be a linear RD-code and $C_2 = (n_2, k_2, d_2)$ linear divisible RD-code defined over $GF(2^N)$. Then the RD-bicode $C = C_1 \cup C_2$ is defined to be a semidivisible RD-bicode.*

**DEFINITION 2.32:** *Let $C_1 = (n_1, k_1, d_1)$ be a MRD code which is not divisible and $C_2 = (n_2, k_2, d_2)$ a divisible RD code defined over $GF(2^N)$; then the bicode $C = C_1 \cup C_2$ is defined to be a semidivisible MRD bicode.*

**DEFINITION 2.33:** *Let $C_1$ be a circulant rank code and $C_2 = (n_2, k_2, d_2)$ a divisible RD code defined over $GF(2^N)$ then $C = C_1 \cup C_2$ is defined to be semidivisible circulant bicode.*

**DEFINITION 2.34:** *Let $C_1$ be a AMRD code and $C_2 = (n_2, k_2, d_2)$ be a divisible RD code defined over $GF(2^N)$ then $C = C_1 \cup C_2$ is defined as a semidivisible AMRD bicode.*

To show the existence of non divisible MRD bicodes, we proceed on to define certain concepts analogous to the ones used in MRD codes.

**DEFINITION 2.35:** *Let $C_1$ be a $(n_2, k_2, d_2)$ MRD code and $C_2 = (n_2, k_2, d_2)$ MRD code over $F_{q^N}$; $n_1 \leq N$ and $n_2 \leq N$; ($n_2 \neq n_2$). $A_{s_1}(n_1, d_1) \cup A_{s_2}(n_2, d_2)$ be the number of bicodewords with rank norms $s_1$ and $s_2$ in the linear $(n_1, k_1, d_1)$ MRD code and $(n_2, k_2, d_2)$ MRD code respectively. Then the bispectrum of the MRD bicode $C_1 \cup C_2$ is described by the formulae;*

$$A_0(n_2, d_2) \cup A_0(n_2, d_2) = 1 \cup 1$$

$$A_{d_1+m_1}(n_1, d_1) \cup A_{d_2+m_2}(n_2, d_2) =$$

$$\begin{bmatrix} n_1 \\ d_1+m_1 \end{bmatrix} \sum_{j_1=0}^{m_1} (-1)^{j_1+m_1} \begin{bmatrix} d_1+m_1 \\ d_1+j_1 \end{bmatrix} q^{\frac{(m_1-j_1)(m_1-j_1-1)}{2}} \left( Q^{j_1+1} - 1 \right)$$

$$\cup \begin{bmatrix} n_2 \\ d_2 + m_2 \end{bmatrix} \sum_{j_2=0}^{m_2} (-1)^{j_2+m_2} \begin{bmatrix} d_2 + m_2 \\ d_2 + j_2 \end{bmatrix} q^{\frac{(m_2-j_2)(m_2-j_2-1)}{2}} (Q^{j_2+1} - 1)$$

$m_1 = 0, 1, \ldots, n_1 - d_1$, $m_2 = 0, 1, \ldots, n_2 - d_2$,
where

$$\begin{bmatrix} n_1 \\ m_1 \end{bmatrix} = \frac{(q^{n_1} - 1)(q^{n_1} - q)\ldots(q^{n_1} - q^{m_1-1})}{(q^{m_1} - 1)(q^{m_1} - q)\ldots(q^{m_1} - q^{m_1-1})}$$

and

$$\begin{bmatrix} n_2 \\ m_2 \end{bmatrix} = \frac{(q^{n_2} - 1)(q^{n_2} - q)\ldots(q^{n_2} - q^{m_2-1})}{(q^{m_2} - 1)(q^{m_2} - q)\ldots(q^{m_2} - q^{m_2-1})}$$

with $Q = q^N$.

Using the bispectrum of a MRD bicode we prove the following theorem:

**THEOREM 2.12:** *All $C_1[n_1, k_1, d_1] \cup C_2[n_2, k_2, d_2]$ MRD bicodes with $d_1 < n_1$ and $d_2 < n_2$ (i.e., with $k_1 \geq 2$ and $k_2 \geq 2$) are non bidivisible.*

*Proof:* This is proved by making use of the bispectrum of the MRD bicodes. Clearly $A_{d_1}(n_1,d_1) \cup A_{d_2}(n_2,d_2) \neq 0 \cup 0$.

If the existence of a bicodeword with birank $d_1 + 1 \cup d_2 + 1$ is established then the proof is complete as bigcd$\{(d_1, d_1 + 1) \cup (d_2, d_2 + 1)\} = 1 \cup 1$.

So the proof is to show that $A_{d_1+1}(n_1,d_1) \cup A_{d_2+1}(n_2,d_2)$ is non zero (i.e., $A_{d_1+1}(n_1,d_1) \neq 0$ and $A_{d_2+1}(n_2,d_2) \neq 0$).
Now
$$A_{d_1+1}(n_1,d_1) \cup A_{d_2+1}(n_2,d_2) =$$
$$\begin{bmatrix} n_1 \\ d_1+1 \end{bmatrix} \left(-\begin{bmatrix} d_1+1 \\ d_1 \end{bmatrix}\right) [(Q-1) + (Q^2 - 1)] \cup$$

$$\begin{bmatrix} n_2 \\ d_2+1 \end{bmatrix}(-\begin{bmatrix} d_2+1 \\ d_2 \end{bmatrix}[(Q-1)+(Q^2-1)] =$$

$$\begin{bmatrix} n_1 \\ d_1+1 \end{bmatrix}(Q-1)+\left(Q+1-\begin{bmatrix} d_1+1 \\ d_1 \end{bmatrix}\right)$$

$$\cup \begin{bmatrix} n_2 \\ d_2+1 \end{bmatrix}(Q-1)+\left(Q+1-\begin{bmatrix} d_2+1 \\ d_2 \end{bmatrix}\right).$$

Suppose that

$$(Q+1)-\begin{bmatrix} d_1+1 \\ d_1 \end{bmatrix} \cup (Q+1)-\begin{bmatrix} d_2+1 \\ d_2 \end{bmatrix} = 0 \cup 0.$$

i.e.,

$$q^N+1 = \frac{q^{d_i+1}-1}{q-1}$$

i.e.,

$$q-1 = \frac{q^{d_i}-1}{q^{N-1}}$$

Clearly

$$\frac{q^{d_i}-1}{q^{N-1}} < 1.$$

For if

$$\frac{q^{d_i}-1}{q^{N-1}} \geq 1$$

then $q^{N-1} < q^{d_i} - 1$, which is not possible as $d_i < n_i \leq N$; i = 1, 2. Thus $q - 1 < 1$ which implies $q < 2$ a contradiction.

Hence $A_{d_1+1}(n_1,d_1) \cup A_{d_2+1}(n_2,d_2)$ is non zero. Thus except $C_1(n_1, 1, n_1) \cup C_2(n_2, 1, n_2)$ MRD bicodes all $C_1[n_1, k_1, d_1] \cup C_2[n_2, k_2, d_2]$ MRD bicodes with $d_1 < n_1$ and $d_2 < n_2$ are non divisible.

Now we finally define the notion of fuzzy rank distance bicodes. Recall von Kaenel [99] introduced the idea of fuzzy codes with Hamming metric. He analysed the distance properties for symmetric error model. Hall and Gur Dial [41] did it for asymmetric and unidirectional error models. Here we define fuzzy RD bicodes.

The study of coding theory resulted from the encounter of noise in communication channels which transmit binary digital data. If a signal 0 or 1 is transmitted electronically it may be distorted into the other signal. A problem occurs when a message in the form of an n-tuple is transmitted, distorted in the channel and received as a new n-tuple representing a different message.

If both $1 \to 0$ and $0 \to 1$ transitions (or errors) appear in a received word with equal probability then the channel is called symmetric channel and the errors are called symmetric error. In an ideal asymmetric channel only one type of error can occur and the error type is know as apriori. Such errors are known as asymmetric. If both $1 \to 0$ and $0 \to 1$ errors can occur in the received words, but in any particular word all error are of one type, then they are called unidirectional errors.

**DEFINITION [41, 99]**: *Let $F_2^n$ denote the n-dimensional vector space of n-tuples over $F_2$. Let $u, v \in F_2^n$ where $u = (u_1, ..., u_n)$ and $v = (v_1, ..., v_n)$. Let p represent the probability that no transition is made and q represent the probability that a transition of the specified type occurs, so that $p + q = 1$. A fuzzy word $f_u$ is the fuzzy subset of $F_2^n$ defined by $f_u = \{(v, f_u(v)) \mid v \in F_2^n\}$ where $f_u(v)$ is the membership function.*

(i) *For the symmetric error model with the Hamming distance, $d = \sum_{i=1}^{n} |u_i - v_i|$, $f_u(v) = p^{n-d}q^d$.*

(ii) *For unidirectional error model,*
$$f_u(v) = \begin{cases} 0 & \text{if } min(k_1, k_2) \neq 0 \\ p^{m-d}q^d & \text{otherwise} \end{cases}$$

*where,*

$$k_1 = \sum_{i=1}^{n} \max(0, u_i - v_i)$$

and

$$k_2 = \sum_{i=1}^{n} \max(0, v_i - u_i)$$

$$d = \begin{cases} k_1 \text{ if } k_2 = 0 \\ k_2 \text{ if } k_1 = 0 \end{cases}$$

$$m = \begin{cases} \sum_{i=1}^{n} u_i \text{ if } k_2 = 0 \\ n - \sum_{i=1}^{n} u_i \text{ if } k_1 = 0 \\ \max\left(\Sigma u_i, n - \Sigma u_i\right) \text{ if } k_1 = k_2 = 0. \end{cases}$$

*For the asymmetric error model*

$$f_u(v) = \begin{cases} 0 \text{ if } \min(k_1, k_2) \neq 0 \\ p^{m-d} q^d \text{ otherwise} \end{cases}$$

*where $d = k_1$ and $m = \Sigma u_i$ for asymmetric $1 \to 0$ error model and $d = k_2$ and $m = n - \Sigma u_i$ for asymmetric $0 \to 1$ error model.*

**DEFINITION 2.36:** *Asymmetric distance $d_a$ between u and v is defined as $d_a(u, v) = \max(k_1, k_2)$, for $u, v \in F_2^n$.*

**DEFINITION 2.37:** *The generalized Hamming distance between fuzzy sets is a metric in the set $f^n = \{f_u: u \in F_2^n\}$ that is $d(f_u, f_v) = \sum_{z \in F_2^n} |f_u(z) - f_v(z)|$ for $f_u, f_v \in f^n$.*

Note that if $u \in F_2^n$ represents a received word and C is a codeword then $f_c(u)$ is the probability that c was transmitted.

**THEOREM 2.13:** *Let $u, v \in F_2^n$ be such that $d_H(u, v) = d$. If $p \neq q$ and $p \neq 0, 1$ then $d(f_u, f_v) = \sum_{i=0}^{d} \binom{d}{i} | p^i q^{d-i} - p^{d-i} q^i |$ for a symmetric error model.*

**THEOREM 2.14:** *Let $u, v \in F_2^n$ be such that $d_a(u, v) = d_a$. If $p \neq q$ and $p \neq 0, 1$ then $d(f_u, f_v) = 2(1 - q^{d_a})$ for an asymmetric error model.*

These two theorems show that the distance between the fuzzy words is dependent only on the Hamming or asymmetric distance (as the case may be) between the base codewords and not on the dimension of the code space. On the other hand it is not so with the unidirectional error model.

Now, we proceed onto recall the definition of fuzzy RD codes and their properties. For more please refer [77, 96].

**DEFINITION 2.38:** *Let $V^n$ denote the n-dimensional vector space of n-tuples over $F_{2^N}$, $n \leq N$ and $N > 1$. Let $u, v \in V^n$ where $u = (u_1, u_2, ..., u_n)$ and $v = (v_1, v_2, ..., v_n)$ with each $u_i, v_i \in F_2^n$. A fuzzy RD word $f_u$ is the fuzzy subset of $V^n$ defined by $f_u = \{(v, f_u(v)) \mid v \in V^n\}$ where $f_u(v)$ is the membership function.*

**DEFINITION 2.39:** *For the symmetric error model, assume p to represent the probability that no transition (i.e., error) is made and q to represent the probability that a rank error occurs so that $p + q = 1$. Then $f_u(v) = p^{n-r} q^r$ where $r = r(u - v; 2) = ||u - v||$.*

**DEFINITION 2.40:** *For the unidirectional and asymmetric error models assume q to represent the probability that $1 \rightarrow 0$ transition or $0 \rightarrow 1$ transition occurs.*

*Then $f_u(v) = \prod_{i=1}^{n} f_{u_i}(v_i)$ where $f_{u_i}(v_i)$ inherits its definition from the above equation for unidirectional and asymmetric* *error models respectively, since each $u_i$ or $v_i$ itself in an N-tuple over $F_2$. That is since $u_i, v_i \in F_{2^N}$ each $u_i$ or $v_i$ itself is an N-tuple from $F_2$.*

*Let $u_i = (u_{i1}, u_{i2}, ..., u_{iN})$ and $v_j = (v_{j1}, v_{j2}, ..., v_{jN})$ where $u_{is}, v_{jt} \in F_2$, $1 \leq s, t \leq N$.*

*Then for unidirectional error model,*

$$f_{u_i}(v_i) = \begin{cases} 0 \text{ if } min(k_{i1}, k_{i2}) \neq 0 \\ p^{m_i - d_i} q^{d_i} \text{ otherwise} \end{cases}$$

*where,*

$$k_{i1} = \sum_{s=1}^{n} max(0, u_{is} - v_{is})$$

*and*

$$k_{i2} = \sum_{s=1}^{n} max(0, v_{is} - u_{is})$$

$$d_i = \begin{cases} k_{i1} \text{ if } k_{i2} = 0 \\ k_{i2} \text{ if } k_{i1} = 0 \end{cases}$$

$$m_i = \begin{cases} \sum_{s=1}^{n} u_{is} \text{ if } k_{i2} = 0 \\ n - \sum_{s=1}^{n} u_{is} \text{ if } k_{i1} = 0 \\ max(\Sigma u_{is}, n - \Sigma u_{is}) \text{ if } k_{i1} = k_{i2} = 0. \end{cases}$$

*For the asymmetric error model*

$$f_{u_i}(v_i) = \begin{cases} 0 \text{ if } min(k_{i1}, k_{i2}) \neq 0 \\ p^{m_i - d_i} q^{d_i} \text{ otherwise} \end{cases}$$

*where $d_i = k_{i1}$ and*

$$m_i = \sum_{s=1}^{N} u_{is}$$

for asymmetric $0 \to 1$ error model $d_i = k_{i2}$ and

$$m_i = n - \sum_{s=1}^{n} u_{is}$$

for asymmetric $0 \to 1$ error model.

**DEFINITION 2.41:** *Let $f^n = \{f_u: u \in V^n\}$. Let $\psi: V_n \to f^n$ be defined as $\psi(u) = f_u$. Clearly $\psi$ is a bijection. Let $C[n, k, d]$ be an RD code which is a subspace of $V^n$ of dimension k and minimum rank distance d. Then $\psi(C) \subseteq f^n$ is a fuzzy RD code. If $c \in C$ then $f_c$ is a fuzzy RD codeword of $\psi(C)$. For any RD code $\psi(C)$ its minimum distance is defined as,*

$$d_{min}(\psi(C)) = \min_{f_a, f_b \in \psi(C)} \{d(f_a, f_b) : f_a \neq f_b\}$$

*where $d(f_a, f_b) = \sum_{z \in V^n} | f_a(z) - f_b(z) |$ is a metric in $f^n$.*

**DEFINITION 2.42:** *If $u \in V^n$ represents a received codeword and $c \in C$ then $f_c(u)$ gives the probability that c was transmitted. Let $\theta(u) = \{f_a \mid a \in C, f_a(u) \geq f_b(u), b \in C\}$. A code for which $|\theta(u)| = 1$ for all $u \in V^n$ is said to be uniquely decodable. In such a case u is decoded as $\psi^{-1}(\theta(u))$.*

Now we proceed on to define the notion of fuzzy RD bicodes.

**DEFINITION 2.43:** *Let $V^{n_1} \cup V^{n_2}$ denote $(n_1, n_2)$ dimensional vector bispace of $(n_1, n_2)$-tuples over $F_{2^N}$; $n_1 \leq N$ and $n_2 \leq N$, $N > 1$. Let $u_1, v_1 \in V^{n_1}, u_2, v_2 \in V^{n_2}$
where*

$$u_1 = (u_1^1, ..., u_{n_1}^1), v_1 = (v_1^1, ..., v_{n_1}^1),$$
$$u_2 = (u_1^2, ..., u_{n_2}^2) \text{ and } v_2 = (v_1^2, ..., v_{n_2}^2)$$

*with $u_i^1, u_i^2, v_i^1, v_i^2 \in F_{2^N}$. A fuzzy RD bicodeword $f_{u_1 \cup u_2} = f_{u_1}^1 \cup f_{u_2}^2$ is a fuzzy bisubset of $V^{n_1} \cup V^{n_2}$ defined by,*

$$f_{u_1 \cup u_2} = f_{u_1}^1 \cup f_{u_2}^2$$

$$= \{(v_1, f^1_{u_1}(v_1)) \mid v_1 \in V^{n_1}\} \cup \{(v_2, f^2_{u_2}(v_2)) \mid v_2 \in V^{n_2}\}$$

where $f^1_{u_1}(v_1) \cup f^2_{u_2}(v_2)$ is the membership bifunction.

**DEFINITION 2.44:** *For the symmetric error bimodel assume $p_1 \cup p_2$ to represent the biprobability that no transition (i.e., error) is made and $q_1 \cup q_2$ to represent the biprobability that a birank error occurs so that $p_1 + q_1 \cup p_2 + q_2 = 1 \cup 1$.*
*Then*

$$f^1_{u_1}(v_1) \cup f^2_{u_2}(v_2) = p_1^{n_1 - r_1} q_1^{r_1} \cup p_2^{n_2 - r_2} q_2^{r_2}$$

*where*

$$r_1 = r_1(u_1 - v_1, 2) = ||u_1 - v_1||$$

*and*

$$r_2 = r_2(u_2 - v_2, 2) = ||u_2 - v_2||.$$

**DEFINITION 2.45:** *For unidirectional and asymmetric error bimodels assume $q_1 \cup q_2$ to represent the probability that $(1 \to 0) \cup (1 \to 0)$ bitransition or $(0 \to 1) \cup (0 \to 1)$ bitransition occurs. Then*

$$f^1_{u_1}(v_1) \cup f^2_{u_2}(v_2) = \prod_{i=1}^{n_1} f^1_{u^1_i}(v^1_i) \cup \prod_{i=1}^{n_2} f^2_{u^2_i}(v^2_i)$$

*where $f^1_{u_1}(v_1) \cup f^2_{u_2}(v_2)$ inherits its definition from the unidirectional and asymmetric bimodels respectively since each $u^1_i \cup u^2_i$ or $v^1_i \cup v^2_i$ itself is an N-bituple over $F_2$.*
*That is since $u^1_i, u^2_i, v^1_i, v^2_i \in F_{2^N}$ each $u^1_i \cup u^2_i$ or $v^1_i \cup v^2_i$ itself is an N-tuple from $F_2$.*

$$u^1_i = (u^1_{i1}, u^1_{i2}, ..., u^1_{iN}), \quad v^1_j = (v^1_{j1}, v^1_{j2}, ..., v^1_{jN}),$$
$$u^2_i = (u^2_{i1}, u^2_{i2}, ..., u^2_{iN}), \quad \text{and} \quad v^2_j = (v^2_{j1}, v^2_{j2}, ..., v^2_{jN})$$

*where $u^1_{is}, u^2_{is}, v^1_{jt}, v^2_{jt} \in F_{2^N}$, $1 \leq s, t \leq N$. Then for unidirectional error bimodel,*

$$f_{u_i^1}^1\left(v_i^1\right) \cup f_{u_i^2}^2\left(v_i^2\right) = \begin{cases} 0 \cup 0 \text{ if } min\left(k_{i1}^1, k_{i2}^1\right) \cup min\left(k_{i1}^2, k_{i2}^2\right) \neq 0 \cup 0 \\ p_1^{m_i^1 - d_i^1} q_1^{d_i^1} \cup p_2^{m_i^2 - d_i^2} q_2^{d_i^2} \text{ otherwise} \end{cases}$$

where

$$k_{i1}^1 = \sum_{s=1}^{N} max\left(0, u_{is}^1 - v_{is}^1\right), \quad k_{i1}^2 = \sum_{s=1}^{N} max\left(0, v_{is}^1 - u_{is}^1\right)$$

$$k_{i2}^1 = \sum_{s=1}^{N} max\left(0, u_{is}^2 - v_{is}^2\right)$$

and

$$k_{i2}^2 = \sum_{s=1}^{N} max\left(0, v_{is}^2 - u_{is}^2\right);$$

where

$$d_i^1 = \begin{cases} k_{i1}^1 \text{ if } k_{i2}^1 = 0 \\ k_{i2}^1 \text{ if } k_{i1}^1 = 0 \end{cases}$$

$$d_i^2 = \begin{cases} k_{i1}^2 \text{ if } k_{i2}^2 = 0 \\ k_{i2}^2 \text{ if } k_{i1}^2 = 0 \end{cases}$$

$$m_i^1 = \begin{cases} \sum_{s=1}^{N} u_{is}^1 \text{ if } k_{i2}^1 = 0 \\ N - \sum_{s=1}^{N} u_{is}^1 \text{ if } k_{i1}^1 = 0 \\ max\left(\Sigma u_{is}^1, N - \Sigma u_{is}^1\right) \text{ if } k_{i1}^1 = k_{i2}^1 = 0 \end{cases}$$

and

$$m_i^2 = \begin{cases} \sum_{s=1}^{N} u_{is}^2 \text{ if } k_{i2}^2 = 0 \\ N - \sum_{s=1}^{N} u_{is}^2 \text{ if } k_{i1}^2 = 0 \\ max\left(\Sigma u_{is}^2, N - \Sigma u_{is}^2\right) \text{ if } k_{i1}^2 = k_{i2}^2 = 0 \end{cases}$$

*for the asymmetric error bimodel*

$$f_{u_i^1}^1(v_i^1) \cup f_{u_i^2}^2(v_i^2) = \begin{cases} 0 \cup 0 \text{ if } \min(k_{i1}^1, k_{i2}^1) \cup \min(k_{i1}^2, k_{i2}^2) \neq 0 \cup 0 \\ p_1^{m_i^1 - d_i^1} q_1^{d_i^1} \cup p_2^{m_i^2 - d_i^2} q_2^{d_i^2} \text{ otherwise} \end{cases}$$

where $d_i^1 = k_{i1}^1$, $d_i^2 = k_{i1}^2$ and

$$m_i^1 = \sum_{s=1}^N u_{is}^1, \quad m_i^2 = \sum_{s=1}^N u_{is}^2$$

*for asymmetric* $(1 \to 0) \cup (1 \to 0)$ *error bimodel and* $d_i^1 = k_{i1}^1$, $d_i^2 = k_{i1}^2$ ,

$$m_i^1 = N - \sum_{s=1}^N u_{is}^1$$

and

$$m_i^2 = N - \sum_{s=1}^N u_{is}^2$$

*for the asymmetric* $(1 \to 0) \cup (0 \to 1)$ *error bimodel.*

The minimum bidistance of a fuzzy RD bicode.

**DEFINITION 2.46:** *Let*
$$f^{n_1} \cup f^{n_2} = \{ f_{u_1}^1 / u_1 \in V^{n_1} \} \cup \{ f_{u_2}^2 / u_2 \in V^{n_2} \}.$$

*Let*
$$\psi = \psi_1 \cup \psi_2 : V^{n_1} \cup V^{n_2} \to f^{n_1} \cup f^{n_2}$$

*defined by* $\psi(u_1) \cup \psi(u_2) = f_{u_1}^1 \cup f_{u_2}^2$. *Clearly $\psi_1$ is a bijection and $\psi_2$ is a bijection. Let $C_1[n_1, k_1, d_1] \cup C_2[n_2, k_2, d_2]$ be an RD bicode which is a subbispace of $V^{n_1} \cup V^{n_2}$ of bidimension $k_1 \cup k_2$ and minimum rank bidistance $d_1 \cup d_2$. Then*
$$\psi_1(C_1) \cup \psi_2(C_2) \subseteq f^{n_1} \cup f^{n_2}$$
*is a fuzzy RD bicode. If $c = c_1 \cup c_2 \in C_1 \cup C_2$ then $f_{c_1} \cup f_{c_2}$ is a fuzzy bicodeword of $\psi_1(C_1) \cup \psi_2(C_2)$.*

*For any fuzzy RD bicode $\psi_1(C_1) \cup \psi_2(C_2)$, its minimum bidistance is defined as,*

$$d_{min}\psi_1(C_1) \cup d_{min}\psi_2(C_2) =$$
$$min\{d(f^1_{a_1}, f^1_{b_1}) / f^1_{a_1}, f^1_{b_1} \in \psi_1(C_1); f^1_{a_1} \neq f^1_{b_1}\} \cup$$
$$min\{d(f^2_{a_2}, f^2_{b_2}) / f^2_{a_2}, f^2_{b_2} \in \psi_2(C_2); f^2_{a_2} \neq f^2_{b_2}\}.$$

*where*

$$d(f^1_{a_1}, f^1_{b_1}) \cup d(f^2_{a_2}, f^2_{b_2}) =$$
$$\sum_{z_1 \in V^{n_1}} |f^1_{a_1}(z_1) - f^1_{b_1}(z_1)| \cup \sum_{z_2 \in V^{n_2}} |f^2_{a_2}(z_2) - f^2_{b_2}(z_2)|$$

*is a bimetric in $f^{n_1} \cup f^{n_2}$.*

**DEFINITION 2.47:** *If $u_1 \cup u_2 \in V^{n_1} \cup V^{n_2}$ represents a received biword and $c_1 \cup c_2 \in C_1 \cup C_2$ then $f^1_{c_1}(u_1) \cup f^2_{c_2}(u_2)$ gives the biprobability that $(c_1 \cup c_2)$ was transmitted.*

$$\theta_1(u_1) \cup \theta_2(u_2) =$$
$$\{f^1_{a_1} | a_1 \in C_1; f^1_{a_1}(u_1) \geq f^1_{b_1}(u_1), b_1 \in C_1\} \cup$$
$$\{f^2_{a_2} | a_2 \in C_2; f^2_{a_2}(u_2) \geq f^2_{b_2}(u_2), b_2 \in C_2\}.$$

*A bicode for which*
$$|\theta_1(u_1) \cup \theta_2(u_2)| = |\theta_1(u_1)| \cup |\theta_2(u_2)|$$
$$= 1 \cup 1$$
*for all $u_1 \in V^{n_1}$ and $u_2 \in V^{n_2}$; is said to be uniquely bidecodable. In such a case $u_1 \cup u_2$ is bicoded as*

$$\psi_1^{-1}(\theta_1(u_1)) \cup \psi_2^{-1}(\theta_2(u_2)).$$

The notions related to m–covering radius of RD–codes can be analogously transformed from the notion of RD–bicodes.

**PROPOSITION 2.1:** *If $C^1_1 \cup C^1_2$ and $C^2_1 \cup C^2_2$ are RD bicodes with $C^1_1 \cup C^1_2 \subseteq C^2_1 \cup C^2_2$ (i.e., $C^1_1 \subseteq C^2_1$ and $C^1_2 \subseteq C^2_2$) then*

$t_{m_1}(C_1^1) \cup t_{m_2}(C_2^1) \geq t_{m_1}(C_1^2) \cup t_{m_2}(C_2^2)$ $[t_{m_1}(C_1^1) \geq t_{m_1}(C_1^2)$ and $t_{m_2}(C_2^1) \geq t_{m_2}(C_2^2)].$

*Proof:* Let $S_1 \subseteq V^{n_1}$ with $|S_1| = m_1$ and $S_2 \subseteq V^{n_2}$ with $|S_2| = m_2$

$$\text{cov}(C_1^2, S_1) \cup \text{cov}(C_2^2, S_2)$$
$$= \min\{\text{cov}(x_1, S_1)/x_1 \in C_1^2\} \cup \min\{\text{cov}(x_2, S_2)/x_2 \in C_2^2\} \leq$$
$$\min\{\text{cov}(x_1, S_1)/x_1 \in C_1^1\} \cup \min\{\text{cov}(x_2, S_2)/x_2 \in C_2^1\}$$
$$= \text{cov}(C_1^1, S_1) \cup \text{cov}(C_2^1, S_2).$$

Thus $t_{m_1}(C_1^2) \cup t_{m_2}(C_2^2) \leq t_{m_1}(C_1^1) \cup t_{m_2}(C_2^1)$.

**PROPOSITION 2.2:** *For any RD bicode $C_1 \cup C_2$ and a pair of positive integers $(m_1, m_2)$,*
$$t_{m_1}(C_1) \cup t_{m_2}(C_2) \leq t_{m_1+1}(C_1) \cup t_{m_2+1}(C_2).$$

*Proof:* Let $S_1 \subseteq V^{n_1}$ with $|S_1| = m_1$ and $S_2 \subseteq V^{n_2}$ with $|S_2| = m_2$, $V^{n_1} \cup V^{n_2}$ is the rank bispace where $S_1 \cup S_2$ is a bisubset of $V^{n_1} \cup V^{n_2}$. Now
$t_{m_1}(C_1) \cup t_{m_2}(C_2)$
$\quad = \quad \max\{\text{cov}(C_1, S_1) \mid S_1 \subseteq V^{n_1}; |S_1| = m_1\} \cup$
$\quad\quad\quad \max\{\text{cov}(C_2, S_2) \mid S_2 \subseteq V^{n_2}; |S_2| = m_2\} \leq$
$\quad\quad\quad \max\{\text{cov}(C_1, S_1) \mid S_1 \subseteq V^{n_1}; |S_1| = m_1 + 1\} \cup$
$\quad\quad\quad \max\{\text{cov}(C_2, S_2) \mid S_2 \subseteq V^{n_2}; |S_2| = m_2 + 1\}$
$\quad = \quad t_{m_1+1}(C_1) \cup t_{m_2+1}(C_2).$

**PROPOSITION 2.3:** *For any biset of positive integers $\{n_1, m_1, k_1, K_1\} \cup \{n_2, m_2, k_2, K_2\}$;*
$$t_{m_1}[n_1, k_1] \cup t_{m_2}[n_2, k_2] \leq t_{m_1+1}[n_1, k_1] \cup t_{m_2+1}[n_2, k_2]$$
*and*
$$t_{m_1}(n_1, K_1) \cup t_{m_2}(n_2, K_2) \leq t_{m_1+1}(n_1, K_1) \cup t_{m_2+1}(n_2, K_2).$$

*Proof:* Given $C_1[n_1, k_1] \cup C_2[n_2, k_2]$ RD bicode, with $C_1 \subseteq V^{n_1}$ and $C_2 \subseteq V^{n_2}$. Now

$$t_{m_1}[n_1,k_1] \cup t_{m_2}[n_2,k_2] =$$
$$\min\{t_{m_1}(C_1)/C_1 \subseteq V^{n_1}; \dim C_1 = k_1\} \cup$$
$$\min\{t_{m_2}(C_2)/C_2 \subseteq V^{n_2}; \dim C_2 = k_2\} \leq$$
$$\min\{t_{m_1+1}(C_1)/C_1 \subseteq V^{n_1}; \dim C_1 = k_1\} \cup$$
$$\min\{t_{m_2+1}(C_2)/C_2 \subseteq V^{n_2}; \dim C_2 = k_2\}$$
$$= t_{m_1+1}[n_1,k_1] \cup t_{m_2+1}[n_2,k_2].$$

Similarly we have

$$t_{m_1}(n_1,K_1) \cup t_{m_2}(n_2,K_2) \leq t_{m_1+1}(n_1,K_1) \cup t_{m_2+1}(n_2,K_2).$$
That is
$$t_{m_1}(n_1,K_1) \cup t_{m_2}(n_2,K_2) =$$
$$\min\{t_{m_1}(C_1)/C_1 \subseteq V^{n_1}; |C_1| = K_1\} \cup$$
$$\min\{t_{m_2}(C_2)/C_2 \subseteq V^{n_2}; |C_2| = K_2\}$$
$$\min\{t_{m_1+1}(C_1)/C_1 \subseteq V^{n_1}; |C_1| = K_1\} \cup$$
$$\min\{t_{m_2+1}(C_2)/C_2 \subseteq V^{n_2}; |C_2| = K_2\}$$
$$\leq t_{m_1+1}(n_1,K_1) \cup t_{m_2+1}(n_2,K_2).$$

**PROPOSITION 2.4:** *For any biset of positive integers $\{n_1, m_1, k_1, K_1\} \cup \{n_2, m_2, k_2, K_2\}$;*
$$t_{m_1}[n_1,k_1] \cup t_{m_2}[n_2,k_2] \geq t_{m_1}[n_1,k_1+1] \cup t_{m_2}[n_2,k_2+1]$$
*and*
$$t_{m_1}(n_1,K_1) \cup t_{m_2}(n_2,K_2) \geq t_{m_1}(n_1,K_1+1) \cup t_{m_2}(n_2,K_2+1).$$

*Proof:* Given $C = C_1 \cup C_2$ is a RD bicode hence a bisubspace of $V^{n_1} \cup V^{n_2}$.
Consider
$$t_{m_1}[n_1,k_1+1] \cup t_{m_2}[n_2,k_2+1] =$$

$$\min\{t_{m_1}(C_1)/C_1 \subseteq V^{n_1}; \dim C_1 = k_1 + 1\} \cup$$
$$\min\{t_{m_2}(C_2)/C_2 \subseteq V^{n_2}; \dim C_2 = k_2 + 1\} \le$$
$$\min\{t_{m_1}(C_1)/C_1 \subseteq V^{n_1}; \dim C_1 = k_1\} \cup$$
$$\min\{t_{m_2}(C_2)/C_2 \subseteq V^{n_2}; \dim C_2 = k_2\}.$$

(since for each $C_1 \cup C_2 \subseteq C_{12} \cup C_{22}$;
$t_{m_1}(C_{12}) \cup t_{m_2}(C_{22}) \le t_{m_1}(C_1) \cup t_{m_2}(C_2)$)
$= t_{m_1}(n_1, k_1) \cup t_{m_2}(n_2, k_2)$.
Similarly,
$t_{m_1}[n_1, K_1 + 1] \cup t_{m_2}[n_2, K_2 + 1] \le t_{m_1}(n_1, K_1) \cup t_{m_2}(n_2, K_2)$.

Using these results and the fact $k_{1m_1}[n_1, t_1]$ denotes the smallest dimension of a linear RD code of length $n_1$ and $m_1$-covering radius $t_1$ and $k_{1m_1}[n_1, t_1]$ denotes the least cardinality of the RD codes of length $n_1$ and $m_1$-covering radius $t_1$.

The following results can be easily proved.

*Result 1:* For any biset of positive integers $\{n_1, m_1, t_1\} \cup \{n_2, m_2, t_2\}$ and
$$k_{m_1}(n_1, t_1) \cup k_{m_2}(n_2, t_2) \le k_{m_1+1}(n_1, t_1) \cup k_{m_2+1}(n_2, t_2)$$
and
$$K_{m_1}(n_1, t_1) \cup K_{m_2}(n_2, t_2) \le K_{m_1+1}(n_1, t_1) \cup K_{m_2+1}(n_2, t_2).$$

*Result 2:* For any biset of positive integers $\{n_1, m_1, t_1\} \cup \{n_2, m_2, t_2\}$ we have
$$k_{m_1}[n_1, t_1] \cup k_{m_2}[n_2, t_2] \ge k_{m_1}[n_1, t_1 + 1] \cup k_{m_2}[n_2, t_2 + 1]$$
and
$$K_{m_1}(n_1, t_1) \cup K_{m_2}(n_2, t_2) \ge K_{m_1}(n_1, t_1 + 1) \cup K_{m_2}(n_2, t_2 + 1).$$
We say a bifunction $f_1 \cup f_2$ is a non-decreasing function in some bivariable say $(x_1 \cup x_2)$ if both $f_1$ and $f_2$ happen to be a non-decreasing function in the same variable $x_1$ and $x_2$ respectively.

With this understanding we can say the $(m_1, m_2)$- covering biradius of a fixed RD bicode $C_1 \cup C_2$,

$$t_{m_1}[n_1, k_1] \cup t_{m_2}[n_2, k_2], \quad t_{m_1}(n_1, K_1) \cup t_{m_2}(n_2, K_2),$$

$$k_{m_1}[n_1, t_1] \cup k_{m_2}[n_2, t_2] \text{ and } K_{m_1}(n_1, t_1) \cup K_{m_2}(n_2, t_2),$$

are non decreasing bifunctions of $(m_1, m_2)$. The relationship between the multicovering biradii of two RD bicodes and bicodes that are built using them are described.

Let $C^i = C_1^i \cup C_2^i$, $i = 1, 2$ be a $[n_1^1, k_1^1, d_1^1] \cup [n_1^2, k_1^2, d_1^2]$, $[n_2^1, k_2^1, d_2^1] \cup [n_2^2, k_2^2, d_2^2]$ RD bicodes over $F_{2^N}$ with $n_1^1, n_1^2, n_2^1, n_2^2, n_1^1 + n_2^1, n_1^2 + n_2^2 \leq N$.

**PROPOSITION 2.5:** *Let $C^1 = C_1^1 \cup C_2^1$ and $C^2 = C_1^2 \cup C_2^2$ be RD bicodes described above.*

$$C = C^1 \times C^2 = C_1^1 \times C_1^2 \cup C_2^1 \times C_2^2$$

$$= \{(x_1 \mid y_1)/x_1 \in C_1^1, y_1 \in C_1^2\} \cup \{(x_2 \mid y_2)/x_2 \in C_2^1, y_2 \in C_2^2\}.$$

*Then $C^1 \times C^2$ is a $[n_1^1 + n_2^1 \cup n_1^2 + n_2^2, k_1^1 + k_2^1 \cup k_1^2 + k_2^2, \min\{d_1^1, d_2^1\} \cup \min\{d_1^2, d_2^2\}]$ rank distance bicode over $F_{2^N}$ and*

$$t_{m_1}(C_1^1 \times C_2^1) \cup t_{m_2}(C_1^2 \times C_2^2) \leq$$
$$t_{m_1}(C_1^2) + t_{m_1}(C_2^1) \cup t_{m_2}(C_1^2) + t_{m_2}(C_2^2).$$

*Proof:* Let

$$S_1 \subseteq V^{n_1^1 + n_2^1} \text{ and } S_2 \subseteq V^{n_1^2 + n_2^2}$$

and

$$S_1 = \{s_1^1, \ldots, s_{m_1}^1\} \text{ and } S_2 = \{s_1^2, \ldots, s_{m_2}^2\}$$

with

$$s_i^1 = (x_{1i}/y_{1i}) \text{ and } s_i^2 = (x_{2i}/y_{2i})$$

$x_{1i} \in V^{n_1^1}, y_{1i} \in V^{n_2^1}, x_{2i} \in V^{n_1^2}$ and $y_{2i} \in V^{n_2^2}$.

Let

$$S_1^1 = \{x_{11} \ldots x_{1m_1}\}, S_1^2 = \{y_{11} \ldots y_{1m_1}\},$$
$$S_2^1 = \{x_{21} \ldots x_{1m_2}\} \text{ and } S_2^2 = \{y_{21} \ldots y_{2m_2}\}.$$

Now $t_{m_1}(C_1^1)$ being the $m_1$-covering radius of $(C_1^1)$ there exists $c_1^1 \in C_1^1$ such that $S_1^1 \subseteq B_{t_{m_1}(C_1^1)}^1(C_1^1)$. This implies $r_1(x_{1i} + c_1^1) \leq t_{m_1}(C_1)$ for all $x_{1i} \in S_1^1$. The same argument is true for $(C_1^2)$. Now consider $(C_2^1)$, this code has $m_2$ covering radius $t_{m_2}(C_2^1)$ such that there exists $c_2^1 \in C_2^1$ such that $S_2^1 \subseteq B_{t_{m_2}(C_2^1)}^1(C_2^1)$. This implies $r_2(x_{2i} + c_2^1) \leq t_{m_2}(C_2)$ for all $x_{2i} \in S_2^1$.

Now
$$C = (C_1^1 | C_2^1) \cup (C_1^2 | C_2^2) = C^1 \cup C^2.$$

Here
$$r_1(s_{1i} + C^1) = r_1((x_{1i} | y_{1i}) + (C_1^1 | C_2^1))$$
$$= r_1(x_{1i} + C_1^1 | y_{1i} + C_2^1) \leq r_1(x_{1i} + C_1^1) + r_1(y_{1i} + C_2^1)$$
$$\leq t_{m_1}(C_1^1) + t_{m_1}(C_1^2).$$

Similarly we have,
$$r_2(s_{2i} + C^2) = r_2((x_{2i} | y_{2i}) + (C_1^2 | C_2^2))$$
$$= r_2(x_{2i} + C_1^2 | y_{2i} + C_2^2) \leq r_2(x_{2i} + C_1^2) + r_2(y_{2i} + C_2^2)$$
$$\leq t_{m_2}(C_2^1) + t_{m_2}(C_2^2).$$

Thus
$$t_m(C) = t_{m_1}(C_1^1 \times C_1^2) \cup t_{m_2}(C_2^1 \times C_2^2) \leq$$
$$t_{m_1}(C_1^1) + t_{m_1}(C_1^2) \cup\ t_{m_2}(C_2^1) + t_{m_2}(C_2^2).$$

For any positive integer r, the r-fold repetition RD code $C_1$ is the code $C = \{(c | c | \ldots | c) | c \in C_1\}$ where the code word c is concatenated r-times. This is a $[rn_1, k_1, d_1]$ rank distance code. Note that here $n_1 \leq N$ is choosen so that $rn_1 \leq N$.

We proceed on to define (r, r)-fold repetition of RD bicode.

**DEFINITION 2.48:** *For any (r, r) (r any positive integer), the (r, r) repetition of RD bicode $C_1 \cup D_1$ is the bicode $C = \{(c | c | \ldots | c) | c \in C_1\} \cup D = \{(d | d | \ldots | d) | d \in D_1\}$ where the bicode word $c \cup d$ is concatenated r-times this is a $[rn_1, k_1, d_1] \cup [rn_2, k_2, d_2]$ rank distance bicode word with $n_i \leq N$ and $rn_i \leq N$; $i = 1$,*

2. *Thus any bicode word in $C \cup D$ would be of the form $(c \mid c \mid ... \mid c) \cup (d \mid d \mid ... \mid d)$ where $c \in C_1$ and $d \in D_1$.*

We can also define $(r_1, r_2)$ fold repetition bicode $(r_1 \neq r_2)$.

**DEFINITION 2.49:** *Let $C_1 \cup D_1$ be a $[n_1, k_1, d_1] \cup [n_2, k_2, d_2]$ RD-code. Let $C = \{(c \mid c \mid ... \mid c) \mid c \in C_1\}$ be a $r_1$-fold repetition RD code $C_1$ and $D = \{(d \mid d \mid ... \mid d) \mid d \in D_1\}$ be a $r_2$-fold repetition RD code $C_2$ $(r_1 \neq r_2)$. Then $C \cup D$ is defined as the $(r_1, r_2)$-fold repetition bicode.*

We prove the following interesting result.

**PROPOSITION 2.6:** *For an $(r, r)$ fold repetition RD-bicode*
$$C \cup D, \ t_{m_1}(C) \cup t_{m_2}(D) = t_{m_1}(C_1) \cup t_{m_2}(D_1).$$

*Proof:* Let $S_1 = \{v_1 ... v_{m_1}\} \subseteq V^{n_1}$ be such that $cov(C_1, S_1) = t_{m_1}(C_1)$. Let $S_2 = \{u_1 ... u_{m_2}\} \subseteq V^{n_2}$ be such that $cov(D_1, S_2) = t_{m_2}(D_1)$.

Let $v_i^1 = (v_i \mid v_i \mid ... \mid v_i)$. Let $S_1^1 = \{v_1^1 \mid v_2^1 \mid ... \mid v_{m_1}^1\}$ be a set of $m_1$-vectors of length $rn_1$ each. A r-fold repetition of any RD code word retains the same rank weight.

Hence $(C, S_1^1) = t_{m_1}(C_1)$.

Since $t_{m_1}(C) \geq cov(C, S_1^1)$, it follows that
$$t_{m_1}(C) \geq t_{m_1}(C_1) \qquad \text{------} \qquad I$$

Conversely let $S_1 = \{v_1 ... v_{m_1}\}$ be a set of m-vectors of length $rn_1$ with $v_i = (v_i^1 \mid v_i^1 \mid ... \mid v_i^1)$; $v_i^1 \in V^{n_1}$. Then there exists $c \in C_1$ such that $d_{R_1}(c, v_i^1) \leq t_{m_1}(C_1)$ for every $i$ $(1 \leq i \leq m_1)$. This implies $d_{R_1}((c \mid c \mid ... \mid c), v_i) \leq t_{m_1}(C_1)$ for every $i$ $(1 \leq i \leq m_1)$.

Thus
$$t_{m_1}(C) \leq t_{m_1}(C_1) \qquad \text{------} \qquad II$$
From I and II,
$$t_{m_1}(C) = t_{m_1}(C_1).$$

On similar lines we can prove,
$$t_{m_2}(D) = t_{m_2}(D_1)$$

where $cov(D_1, S_2) = t_{m_2}(D_1)$.

Hence
$$t_{m_1}(C) \cup t_{m_2}(D) = t_{m_1}(C_1) \cup t_{m_2}(D_1).$$

Multi-covering bibounds for RD-bicodes is discussed and a few interesting properties in this direction are given. The $(m_1, m_2)$ covering biradius $t_{m_1}(C_1) \cup t_{m_2}(C_2)$ of a RD-bicode $C = C_1 \cup C_2$ is a non-decreasing bifunction of $m_1 \cup m_2$ (proved earlier). Thus a lower bi-bound for $t_{m_1}(C_1) \cup t_{m_2}(C_2)$ implies a bibound for $t_{m_1+1}(C_1) \cup t_{m_2+1}(C_2)$. First bibound exhibits $m_1 \cup m_2 \geq 2 \cup 2$ then situation of $(m_1, m_2)$-covering biradii is quite different for ordinary covering radii.

**PROPOSITION 2.7:** *If $m_1 \cup m_2 \geq 2 \cup 2$ then the $(m_1, m_2)$-covering biradii of a RD bicode $C = C_1 \cup C_2$ of bilength $(n_1, n_2)$ is atleast $\left\lceil \dfrac{n_1}{2} \right\rceil \cup \left\lceil \dfrac{n_2}{2} \right\rceil$.*

*Proof:* Let $C = C_1 \cup C_2$ be a RD bicode of bilength $(n_1, n_2)$ over $GF(2^N)$. Let $m_1 \cup m_2 \geq 2 \cup 2$, let $t_1, t_2$ be the 2-covering biradii of the RD code $C = C_1 \cup C_2$. Let $x = x_1 \cup x_2 \in V^{n_1} \cup V^{n_2}$.

Choose $y = y_1 \cup y_2 \in V^{n_1} \cup V^{n_2}$ such that all the $(n_1, n_2)$ coordinates of $x - y = (x_1 - y_1) \cup (x_2 - y_2)$ are linearly independent, that is
$$d_R(x, y) = d_R(x_1 \cup x_2, y_1 \cup y_2)$$
$$= d_{R_1}(x_1, y_1) \cup d_{R_2}(x_2, y_2)$$
$$(R = R_1 \cup R_2 \text{ and } d_{R_1} = d_{R_1} \cup d_{R_2})$$
$$= n_1 \cup n_2.$$
Then for any $c = c_1 \cup c_2 \in C = C_1 \cup C_2$,
$$d_R(x, c) + d_R(c, y)$$

$$= d_{R_1}(x_1,c_1) + d_{R_1}(c_1,y_1) \cup d_{R_2}(x_2,c_2) + d_{R_2}(c_2,y_2)$$
$$\geq d_{R_1}(x_1,y_1) \cup d_{R_2}(x_2,y_2)$$
$$= n_1 \cup n_2.$$

This implies that one of
$$d_{R_1}(x_1,c_1) \cup d_{R_2}(x_2,c_2)$$
and
$$d_{R_1}(c_1,y_1) \cup d_{R_2}(c_2,y_2)$$
is at least $\frac{n_1}{2} \cup \frac{n_2}{2}$ (that is one of $d_{R_1}(x_1,c_1)$ and $d_{R_1}(c_1,y_1)$ is atleast $\frac{n_1}{2}$ and one of $d_{R_2}(x_2,c_2)$ and $d_{R_2}(c_2,y_2)$ is at least $\frac{n_2}{2}$) and hence

$$t = t_1 \cup t_2 \geq \left\lceil \frac{n_1}{2} \right\rceil \cup \left\lceil \frac{n_2}{2} \right\rceil.$$

Since t is non decreasing bifunction of $m_1 \cup m_2$ it follows that
$$t_m(C) = t_{m_1}(C_1) \cup t_{m_2}(C_2) \geq \left\lceil \frac{n_1}{2} \right\rceil \cup \left\lceil \frac{n_2}{2} \right\rceil$$
for $m_1 \cup m_2 \geq 2 \cup 2$.

Bibounds of the multi-covering biradius of $V^{n_1} \cup V^{n_2}$ can be used to obtain bibounds on the multi covering biradii of arbitary bicodes. Thus a relationship between $(m_1, m_2)$-covering biradii of an RD bicode and that of its ambient bispace $V^{n_1} \cup V^{n_2}$ is established

**THEOREM 2.15:** *Let $C = C_1 \cup C_2$ be any RD-code of bilength $n_1 \cup n_2$ over $F_{2^N} \cup F_{2^N}$. Then for any pair of positive integers $(m_1, m_2)$;*
$$t_{m_1}^1(C_1) \cup t_{m_2}^2(C_2) \leq t_1^1(C_1) + t_{m_1}^1(V^{n_1}) \cup t_1^2(C_2) + t_{m_2}^2(V^{n_2}).$$

*Proof:* Let $S = S_1 \cup S_2 \subseteq V^{n_1} \cup V^{n_2}$ (i.e., $S_1 \subseteq V^{n_1}$ and $S_2 \subseteq V^{n_2}$) with $|S| = |S_1| \cup |S_2| = m_1 \cup m_2$. Then there exists $u = u_1 \cup u_2 \in V^{n_1} \cup V^{n_2}$ such that

$$\text{cov}(u, S) = \text{cov}(u_1, S_1) \cup \text{cov}(u_2, S_2) \leq t^1_{m_1}(V^{n_1}) \cup t^1_{m_1}(V^{n_2}).$$

Also there is a $c = c_1 \cup c_2 \in C_1 \cup C_2$ such that

$$d_R(c, u) = d_{R_1}(c_1, u_1) \cup d_{R_2}(c_2, u_2) \leq t^1_1(C_1) \cup t^2_1(C_2).$$

Now

$$\text{cov}(c, S) = \text{cov}(c_1, S_1) \cup \text{cov}(c_2, S_2)$$

$$= \max\{d_{R_1}(c_1, y_1)/y_1 \in S_1\} \cup \max\{d_{R_2}(c_2, y_2)/y_2 \in S_2\}$$

$$\leq \max\{d_{R_1}(c_1, u_1) + d_{R_1}(u_1, y_1)/y_1 \in S_1\} \cup$$
$$\max\{d_{R_2}(c_2, u_2) + d_{R_2}(u_2, y_2)/y_2 \in S_2\}$$

$$= d_{R_1}(c_1, u_1) + \text{cov}(u_1, S_1) \cup d_{R_2}(c_2, u_2) + \text{cov}(u_2, S_2)$$

$$\leq t^1_1(C_1) + t^1_{m_1}(V^{n_1}) \cup t^2_1(C_2) + t^2_{m_2}(V^{n_2}).$$

Thus for every $S = S_1 \cup S_2 \subseteq V^{n_1} \cup V^{n_2}$ with $|S| = m = |S_1| \cup |S_2| = m_1 \cup m_2$ one can find a $c = c_1 \cup c_2 \in C_1 \cup C_2$ such that,

$$\text{cov}(c, S) = \text{cov}(c_1, S_1) \cup \text{cov}(c_2, S_2)$$
$$\leq t^1_1(C_1) + t^1_{m_1}(V^{n_1}) \cup t^2_1(C_2) + t^2_{m_2}(V^{n_2}).$$

Since

$$\text{cov}(c, S) = \text{cov}(c_1, S_1) \cup \text{cov}(c_2, S_2)$$
$$= \min\{\text{cov}(a_1, S_1)/a_1 \in C_1\} \cup \min\{\text{cov}(a_2, S_2)/a_2 \in C_2\}$$
$$\leq \{t^1_1(C_1) + t^1_{m_1}(V^{n_1})\} \cup \{t^2_1(C_2) + t^2_{m_2}(V^{n_2})\};$$

for all $S = S_1 \cup S_2 \subseteq V^{n_1} \cup V^{n_2}$ with $|S| = |S_1| \cup |S_2| = m_1 \cup m_2$, it follows that

$$t^1_{m_1}(C_1) \cup t^2_{m_2}(C_2) = \max\{\text{cov}(C_1, S_1)/S_1 \subseteq V^{n_1}; |S_1| = m_1\} \cup$$

$$\max\left\{\text{cov}(C_2, S_2)/S_2 \subseteq V^{n_2}; |S_2| = m_2\right\}$$
$$\leq \left\{t_1^1(C_1) + t_{m_1}^1(V^{n_1})\right\} \cup \left\{t_1^2(C_2) + t_{m_2}^2(V^{n_2})\right\}.$$

**PROPOSITION 2.8:** *For any pair of integers $(n_1, n_2)$; $n_1 \cup n_2 \geq 2 \cup 2$, $t_2^1(V^{n_1}) \cup t_2^2(V^{n_2}) \leq n_1 - 1 \cup n_2 - 1$; where $V^{n_1} = F_{2^N}^{n_1}$, $V^{n_2} = F_{2^N}^{n_2}$; $n_1 \leq N$ and $n_2 \leq N$.*

*Proof:* Let
$$x_1 = (x_1^1, x_2^1, ..., x_{n_1}^1), \ y_1 = (y_1^1, y_2^1, ..., y_{n_1}^1) \in V^{n_1}$$
and
$$x_2 = (x_1^2, x_2^2, ..., x_{n_2}^2), \ y_2 = (y_1^2, y_2^2, ..., y_{n_2}^2) \in V^{n_2}.$$
Let
$$u = u_1 \cup u_2 \in V^{n_1} \cup V^{n_2}$$
where
$$u_1 = (x_1^1 u_2^1 u_3^1 ... u_{n_1-1}^1 y_{n_1}^1)$$
and
$$u_2 = (x_1^2 u_2^2 u_3^2 ... u_{n_2-1}^2 y_{n_2}^2).$$

Thus $u = u_1 \cup u_2$ bicovers $x_1 \cup x_2$ and $y_1 \cup y_2 \in V^{n_1} \cup V^{n_2}$ with in a biradius $n_1 - 1 \cup n_2 - 1$ as
$$d_{R_1}(u_1, x_1) \leq n_1 - 1 \text{ and } d_{R_2}(u_2, x_2) \leq n_2 - 1$$
and
$$d_{R_1}(u_1, y_1) \leq n_1 - 1 \text{ and } d_{R_2}(u_2, y_2) \leq n_2 - 1.$$

Thus for any pair of bivectors $x_1 \cup x_2$ and $y_1 \cup y_2$ in $V^{n_1} \cup V^{n_2}$ there always exists a bivector namely $u = u_1 \cup u_2$ which bicovers $x_1 \cup x_2$ and $y_1 \cup y_2$ within a biradius $n_1 - 1 \cup n_2 - 1$. Hence
$$t_2^1(V^{n_1}) \cup t_2^2(V^{n_2}) \leq (n_1 - 1 \cup n_2 - 1).$$

Now we proceed on to describe the notion of generalized sphere bicovering bibounds for RD bicodes. A natural question is for a given $t^1 \cup t^2$, $m_1 \cup m_2$ and $n_1 \cup n_2$ what is the smallest RD bicode whose $m_1 \cup m_2$ bicovering biradius is atmost $t^1 \cup t^2$.

As it turns out even for $m_1 \cup m_2 \geq 2 \cup 2$, it is necessary that $t^1 \cup t^2$ be atleast $\frac{n_1}{2} \cup \frac{n_2}{2}$. Infact the minimal $t^1 \cup t^2$ for which such a bicode exists is the $(m_1, m_2)$ bicovering biradius of $C_1 \cup C_2 = F_{2^N}^{n_1} \cup F_{2^N}^{n_2}$.

Various external values associated with this notion are $t_{m_1}^1(V^{n_1}) \cup t_{m_2}^2(V^{n_2})$ the smallest $(m_1, m_2)$-covering biradius among bilength $n_1 \cup n_2$ RD bicodes $t_{m_1}^1(n_1, K^1) \cup t_{m_2}^2(n_2, K^2)$, the smallest $(m_1, m_2)$ covering biradius among all $(n_1, K^1) \cup (n_2, K^2)$ RD bicodes. $K_{m_1}^1(n_1, t^1) \cup K_{m_2}^2(n_2, t^2)$ is the smallest bicardinality of bilength $n_1 \cup n_2$ RD bicode with $m_1 \cup m_2$ covering biradius $t^1 \cup t^2$ and so on. It is the latter quality that is studied in the book for deriving new lower bibounds.

From the earlier results $K_{m_1}^1(n_1, t^1) \cup K_{m_2}^2(n_2, t^2)$ is undefined if
$$t^1 \cup t^2 < \frac{n_1}{2} \cup \frac{n_2}{2}.$$
When this is the case, it is accepted to say
$$K_{m_1}^1(n_1, t^1) \cup K_{m_2}^2(n_2, t^2) = \infty \cup \infty.$$
There are other circumstances when $K_{m_1}^1(n_1, t^1) \cup K_{m_2}^2(n_2, t^2)$ is undefined.

For instance
$$K_{2^{Nn_1}}^1(n_1, n_1 - 1) \cup K_{2^{Nn_2}}^2(n_2, n_2 - 1) = \infty \cup \infty.$$
Also
$$K_{m_1}^1(n_1, t^1) \cup K_{m_2}^2(n_2, t^2) = \infty \cup \infty,$$
$$m_1 > V(n_1, t^1)$$
and
$$m_2 > V(n_2, t^2),$$
since in this case no biball of biradius $t^1 \cup t^2$ covers any biset of $m_1 \cup m_2$ distinct bivectors. More generally one has the fundamental issue of whether
$$K_{m_1}^1(n_1, t^1) \cup K_{m_2}^2(n_2, t^2)$$

is bifinite for a given $n_1$, $m_1$, $t^1$ and $n_2$, $m_2$, $t^2$.

This is the case if and only if

$$t^1_{m_1}(V^{n_1}) \leq t^1 \text{ and } t^2_{m_2}(V^{n_2}) \leq t^2,$$

since $t^1_{m_1}(V^{n_1}) \cup t^2_{m_2}(V^{n_2})$ lower bibounds the $(m_1, m_2)$ covering biradii of all other bicodes of bidimension $n_1 \cup n_2$. When $t_1 \cup t_2 = n_1 \cup n_2$ every bicode word bicovers every bivector, so a bicode of size $1 \cup 1$ will $(m_1, m_2)$ bicover $V^{n_1} \cup V^{n_2}$ for every $m_1 \cup m_2$.

Thus $K^1_{m_1}(n_1, n_1) \cup K^2_{m_2}(n_2, n_2) = 1 \cup 1$ for every $m_1 \cup m_2$. If $t^1 \cup t^2$ is $n_1 - 1 \cup n_2 - 1$ what happens to $K^1_{m_1}(n_1, t^1) \cup K^2_{m_2}(n_2, t^2)$?

When $m_1 = m_2 = 1$,

$$K^1_1(n_1, n_1 - 1) \cup K^2_1(n_2, n_2 - 1) \leq 1 + L_{n_1}(n_1) \cup 1 + L_{n_2}(n_2).$$

For $\overline{0} \cup \overline{0} = (0 \ldots 0) \cup (0 \ldots 0)$ will cover all bivectors of birank binorm less than or equal $n_1 - 1 \cup n_2 - 1$ within biradius $n_1 - 1 \cup n_2 - 1$. That is $\overline{0} \cup \overline{0} = (0, 0, \ldots, 0) \cup (0, 0, \ldots, 0)$ will bicover all binorm $n_1 - 1 \cup n_2 - 1$ bivectors within the biradius $n_1 - 1 \cup n_2 - 1$.

Hence remaining bivectors are rank $n_1 \cup n_2$ bivectors. Thus $\overline{0} \cup \overline{0} = (0, 0, \ldots, 0) \cup (0, 0, \ldots, 0)$ and these birank-$(n_1 \cup n_2)$ bivectors can bicover the ambient bispace within the biradius $n_1 - 1 \cup n_2 - 1$.

Therefore

$$K^1_1(n_1, n_1 - 1) \cup K^2_1(n_2, n_2 - 1) \leq 1 + L_{n_1}(n_1) \cup 1 + L_{n_2}(n_2).$$

**PROPOSITION 2.9:** *For any RD bicode of bilength $n_1 \cup n_2$ over $F_{2^N} \cup F_{2^N}$*

$$K^1_{m_1}(n_1, n_1 - 1) \cup K^2_{m_2}(n_2, n_2 - 1) \leq m_1 L_{n_1}(n_1) + 1 \cup m_2 L_{n_2}(n_2) + 1$$

*provided $m_1 \cup m_2$ is such that*

$$m_1 L_{n_1}(n_1) + 1 \cup m_2 L_{n_2}(n_2) + 1 \leq |V^{n_1}| + |V^{n_2}|.$$

*Proof:* Consider a RD-bicode $C = C_1 \cup C_2$ such that $|C| = |C_1| \cup |C_2| = m_1 L_{n_1}(n_1) + 1 \cup m_2 L_{n_2}(n_2) + 1$. Each bivector in $V^{n_1} \cup V^{n_2}$ has $L_{n_1}(n_1) \cup L_{n_2}(n_2)$ rank complements, that is from each bivector $v_1 \cup v_2 \in V^{n_1} \cup V^{n_2}$; there are $L_{n_1}(n_1) \cup L_{n_2}(n_2)$ bivectors at rank bidistance $n_1 \cup n_2$. This means for any set $S_1 \cup S_2 \subseteq V^{n_1} \cup V^{n_2}$ of $(m_1, m_2)$ bivectors there always exists a $c_1 \cup c_2 \in C_1 \cup C_2$ which bicovers $S_1 \cup S_2$ birank distance $n_1 - 1 \cup n_2 - 1$.

Thus,
$$\text{cov}(c_1, S_1) \cup \text{cov}(c_2, S_2) \leq n_1 - 1 \cup n_2 - 1$$

which implies $\text{cov}(C_1, S_1) \cup \text{cov}(C_2, S_2) \leq n_1 - 1 \cup n_2 - 1$.
Hence

$$K^1_{m_1}(n_1, n_1 - 1) \cup K^2_{m_2}(n_2, n_2 - 1) \leq m_1 L_{n_1}(n_1) + 1 \cup m_2 L_{n_2}(n_2) + 1.$$

By bounding the number of $(m_1, m_2)$ bisets that can be covered by a given bicode word, one obtains a straight forward generalization of the classical sphere bibound.

**THEOREM 2.16:** *(Generalized Sphere Bound for RD bicodes) For any $(n_1, K^1) \cup (n_2, K^2)$ RD bicode $C = C_1 \cup C_2$,*

$$K^1 \binom{V(n_1, t_{m_1}(C_1))}{m_1} \cup K^2 \binom{V(n_2, t_{m_2}(C_2))}{m_2} \geq \binom{2^{N_{n_1}}}{m_1} \cup \binom{2^{N_{n_2}}}{m_2}.$$

*Hence for any $n_1$, $t_1$ and $m_1$, $n_2$, $t_2$ and $m_2$*

$$K^1_{m_1}(n_1, t_1) \cup K^2_{m_2}(n_2, t_2) \geq \frac{\left[\begin{array}{c} 2^{N_{n_1}} \\ m_1 \end{array}\right]}{\left[\begin{array}{c} V(n_1, t_1) \\ m_1 \end{array}\right]} \cup \frac{\left[\begin{array}{c} 2^{N_{n_2}} \\ m_2 \end{array}\right]}{\left[\begin{array}{c} V(n_2, t_2) \\ m_2 \end{array}\right]}$$

*where*

$$V(n_1,t_1) \cup V(n_1,t_1) = \sum_{i_1=0}^{t_1} L^1_{i_1}(n_1) \cup \sum_{i_2=0}^{t_2} L^2_{i_2}(n_2),$$

*number of bivectors in a sphere of biradius $t^1 \cup t^2$ and $L^1_{i_1}(n_1) \cup L^2_{i_2}(n_2)$ is the number of bivectors in $V^{n_1} \cup V^{n_2}$ whose rank binorm is $i_1 \cup i_2$.*

*Proof:* Each set of $(m_1, m_2)$ bivectors in $V^{n_1} \cup V^{n_2}$ must occur in a sphere of biradius $t_{m_1}(C_1) \cup t_{m_2}(C_2)$ around at least one code biword.

Total number of such bisets is $|V^{n_1}| \cup |V^{n_2}|$, choose $m_1 \cup m_2$, where

$$|V^{n_1}| \cup |V^{n_2}| = 2^{N_{n_1}} \cup 2^{N_{n_2}}.$$

The number of bisets of $(m_1, m_2)$ bivectors in a neighborhood of biradius $t_{m_1}(C_1) \cup t_{m_2}(C_2)$ is

$$V(n_1, t_{m_1}(C_1)) \cup V(n_2, t_{m_2}(C_2)).$$

Choose $m_1 \cup m_2$. There are K code biwords.

Hence

$$K^1 \binom{V(n_1, t_{m_1}(C_1))}{m_1} \cup K^2 \binom{V(n_2, t_{m_2}(C_2))}{m_2} \geq \binom{2^{N_{n_1}}}{m_1} \cup \binom{2^{N_{n_2}}}{m_2}.$$

Thus for any $n_1 \cup n_2$, $t^1 \cup t^2$ and $m_1 \cup m_2$,

$$K^1_{m_1}(n_1,t_1) \cup K^2_{m_2}(n_2,t_2) \geq \frac{\left[\begin{array}{c} 2^{N_{n_1}} \\ m_1 \end{array}\right]}{\left[\begin{array}{c} V(n_1,t_1) \\ m_1 \end{array}\right]} \cup \frac{\left[\begin{array}{c} 2^{N_{n_2}} \\ m_2 \end{array}\right]}{\left[\begin{array}{c} V(n_2,t_2) \\ m_2 \end{array}\right]}.$$

**COROLLARY 2.3:** *If*

$$\begin{bmatrix} 2^{N_{n_1}} \\ m_1 \end{bmatrix} \cup \begin{bmatrix} 2^{N_{n_2}} \\ m_2 \end{bmatrix} > 2^{N_{n_1}} \begin{pmatrix} V(n_1,t_1) \\ m_1 \end{pmatrix} \cup 2^{N_{n_2}} \begin{pmatrix} V(n_2,t_2) \\ m_2 \end{pmatrix},$$

*then* $K^1_{m_1}(n_1,t_1) \cup K^2_{m_2}(n_2,t_2) = \infty \cup \infty$.

Chapter Three

# RANK DISTANCE m-CODES

In this chapter we introduce the new notion of rank distance m-codes and describe some of their properties.

**DEFINITION 3.1:** *Let $C_1 = C_1[n_1, k_1]$, $C_2 = C_2[n_2, k_2]$, ..., $C_m = C_m[n_m, k_m]$, be m distinct RD codes such that $C_i = C_i[n_i, k_i] \neq C_j = C_j[n_j, k_j]$ if $i \neq j$ and $C_i = C_i[n_i, k_i] \not\subset C_j = C_j[n_j, k_j]$ or $C_j = C_j[n_j, k_j] \subseteq C_i = C_i[n_i, k_i]$ for $1 \leq i, j \leq m$ if $i \neq j$; be subspaces of the rank spaces $V^{n_1}, V^{n_2}, ..., V^{n_m}$ over the field $GF(2^N)$ or $F_{q^N}$ where $n_1, n_2, ..., n_m \leq N$, i.e., each $n_i \leq N$ for $i = 1, 2, ..., m$. $C = C_1[n_1, k_1] \cup C_2[n_2, k_2] \cup ... \cup C_m[n_m, k_m]$ is defined as the Rank Distance m-code ($m \geq 3$). If $m = 3$ we call the Rank Distance 3-code as the Rank Distance tricode. We say $V^{n_1} \cup V^{n_2} \cup ... \cup V^{n_m}$ to be $(n_1, n_2, ..., n_m)$ dimensional vector m-space over the field $F_{q^N}$. So we can say $C_1 \cup C_2 \cup ... \cup C_m$ is a m-subspace of the* *vector m-space $V^{n_1} \cup V^{n_2} \cup ... \cup V^{n_m}$. We represent any element of $V^{n_1} \cup V^{n_2} \cup ... \cup V^{n_m}$ by*

$$x_1 \cup x_2 \cup ... \cup x_n =$$
$$(x_1^1, x_2^1, ..., x_{n_1}^1) \cup (x_1^2, x_2^2, ..., x_{n_2}^2) \cup ... \cup (x_1^m, x_2^m, ..., x_{n_m}^m),$$

*where $x_{i_j}^i \in F_{q^n}$; $1 \le j \le m$ and $1 \le i_j \le n_j$; $j = 1, 2, 3, ..., m$. Also $F_{q^N}$ can be considered as a pseudo false m-space of dimension $\underbrace{(N, ..., N)}_{m-times}$ over $F_q$.*

Thus elements $x_i, y_i \in F_{q^N}$ has N - m-tuple representation as

$$(\alpha_{1j}^1, \alpha_{2j}^1, ..., \alpha_{Nj_1}^1) \cup (\alpha_{1j}^2, \alpha_{2j}^2, ..., \alpha_{Nj_2}^2) \cup ...$$
$$\cup (\alpha_{1j}^m, \alpha_{2i}^m, ..., \alpha_{Nj_m}^m)$$

over $F_q$ with respect to some m-basis. Hence associated with each $x^1 \cup x^2 \cup ... \cup x^m \in V^{n_1} \cup V^{n_2} \cup ... \cup V^{n_m}$ ($n_i \ne n_j$ if $i \ne j$, $1 \le i, j \le m$) there is a m-matrix.

$$m_1(x^1) \cup ... \cup m_m(x^m) =$$

$$\begin{bmatrix} a_{11}^1 & \cdots & a_{1n_1}^1 \\ a_{21}^1 & \cdots & a_{2n_1}^1 \\ \vdots & & \vdots \\ a_{N1}^1 & \cdots & a_{Nn_1}^1 \end{bmatrix} \cup \begin{bmatrix} a_{11}^2 & \cdots & a_{1n_2}^2 \\ a_{21}^2 & \cdots & a_{2n_2}^2 \\ \vdots & & \vdots \\ a_{N2}^2 & \cdots & a_{Nn_2}^2 \end{bmatrix} \cup ... \cup \begin{bmatrix} a_{11}^m & \cdots & a_{1n_m}^m \\ a_{21}^m & \cdots & a_{2n_m}^m \\ \vdots & & \vdots \\ a_{Nm}^m & \cdots & a_{Nn_m}^m \end{bmatrix}$$

where $i_1^{th} \cup i_2^{th} \cup ... \cup i_m^{th}$ m-column represents the $i_1^{th} \cup i_2^{th} \cup ... \cup i_m^{th}$ m–coordinate of $x_{i1}^1 \cup x_{i2}^2 \cup ... \cup x_{im}^m$ of $x^1 \cup x^2 \cup ... \cup x^m$ over $F_q$.

It is important and interesting to note in order to develop the new notion rank distance m-codes ($m \ge 3$) and while trying to give m- matrices and m- ranks associated with them we are forced to define the notion of pseudo false m- vector spaces. For example, $Z_2^5 \cup Z_2^5 \cup ... \cup Z_2^5$ is a false pseudo m-vector space over $Z_2$. $Z_5^9 \cup Z_5^9 \cup Z_5^9 \cup Z_5^9 \cup Z_5^9 \cup Z_5^9$ is a false pseudo 6-vector space over $Z_5$.

However we will in this book use only m-vector spaces over $Z_2$ or $Z_{2^N}$ and denote it by GF(2) or GF($2^N$) or $F_{2^N}$, unless otherwise specified. Now we see every $x^1 \cup \ldots \cup x^m$ in the m-vector space $V^{n_1} \cup V^{n_2} \cup \ldots \cup V^{n_m}$ have an associated m-matrix $m_1(x^1) \cup \ldots \cup m_m(x^m)$.

We proceed to define m- rank of the m- matrix over $F_q$ or GF(2).

**DEFINITION 3.2:** *The m-rank of an element $x^1 \cup x^2 \cup \ldots \cup x^m \in V^{n_1} \cup V^{n_2} \cup \ldots \cup V^{n_m}$ is defined as the m-rank of the m-matrix $m(x^1) \cup m(x^2) \cup \ldots \cup m(x^m)$ over GF(2) or $F_q$ [i.e., the m-rank of $m(x^1) \cup m(x^2) \cup \ldots \cup m(x^m)$ is the rank of $m(x^1) \cup$ rank of $m(x^2) \cup \ldots \cup$ rank of $m(x^m)$]. We shall denote the m-rank of $x^1 \cup x^2 \cup \ldots \cup x^m$ by $r_1(x^1) \cup r_2(x^2) \cup \ldots \cup r_m(x^m)$, we can in case of m-rank of a m-matrix prove the following;*

(i) *For every $x^1 \cup x^2 \cup \ldots \cup x^m \in V^{n_1} \cup V^{n_2} \cup \ldots \cup V^{n_m}$ ($x_i \in V^{n_i}$, $1 \leq i \leq m$) we have, $r(x^1 \cup \ldots \cup x^m) = r_1(x^1) \cup \ldots \cup r_m(x^m) \geq 0 \cup 0 \cup \ldots \cup 0$ (i.e., each $r_i(x^i) \geq 0$ for every $x^i \in V^{n_i}$; $i = 1, 2, 3, \ldots, m$).*

(ii) *$r(x^1 \cup x^2 \cup \ldots \cup x^m) = r_1(x^1) \cup r_2(x^2) \cup \ldots \cup r_m(x^m) = 0 \cup 0 \cup \ldots \cup 0$ if and only if $x^i = 0$ for $i = 1, 2, 3, \ldots, m$.*

$$r[(x_1^1 + x_2^1) \cup (x_1^2 + x_2^2) \cup \ldots \cup (x_1^m + x_2^m)]$$
$$\leq r_1(x_1^1) + r_1(x_2^1) \cup r_2(x_1^2) + r_2(x_2^2) \cup \ldots \cup r_m(x_1^m) + r_m(x_2^m)$$

*for every $x_1^i, x_2^i \in V^{n_i}$; $i = 1, 2, 3, \ldots, m$. That is we have,*

(iii) $\qquad r[(x_1^1 + x_2^1) \cup (x_1^2 + x_2^2) \cup \ldots \cup (x_1^m + x_2^m)]$
$$= r_1(x_1^1 + x_2^1) \cup r_2(x_1^2 + x_2^2) \cup \ldots \cup r_m(x_1^m + x_2^m)$$
$$\leq r_1(x_1^1) + r_1(x_2^1) \cup r_2(x_1^2) + r_2(x_2^2) \cup \ldots \cup r_m(x_1^m) + r_m(x_2^m)$$
*(as we have for every*
$$x_1^i, x_2^i \in V^{n_i}; \; r_i(x_1^i + x_2^i) \leq r_i(x_1^i) + r_i(x_2^i)$$
*for $i = 1, 2, 3, \ldots, m$).*

(iv) $r_1(a_1x_1) \cup r_2(a_2x_2) \cup \ldots \cup r_m(a_mx_m) = |a_1|r_1(x_1) \cup |a_2|r_2(x_2) \cup \ldots \cup |a_m|r_m(x_m)$ for every $a_1, a_2, \ldots, a_m \in F_q$ or $GF(2)$ and for every $x_i \in V^{n_i}$; $i = 1, 2, 3, \ldots, m$.

Thus the m-function $x_1 \cup x_2 \cup \ldots \cup x_m \to r_1(x_1) \cup r_2(x_2) \cup \ldots \cup r_m(x_m)$ defines a m-norm on $V^{n_1} \cup V^{n_2} \cup \ldots \cup V^{n_m}$.

**DEFINITION 3.3:** *The m-metric induced by the m- rank m- norm is defined as the m- rank m-metric on $V^{n_1} \cup V^{n_2} \cup \ldots \cup V^{n_m}$ and is denoted by $d_{R_1} \cup d_{R_2} \cup \ldots \cup d_{R_m}$. If $x_1^1 \cup x_2^1 \cup \ldots \cup x_m^1$, $y_1^1 \cup y_2^1 \cup \ldots \cup y_m^1 \in V^{n_1} \cup V^{n_2} \cup \ldots \cup V^{n_m}$ then the m-rank m-distance between $x_1^1 \cup x_2^1 \cup \ldots \cup x_m^1$ and $y_1^1 \cup y_2^1 \cup \ldots \cup y_m^1$ is*
$$d_{R_1}(x_1^1, y_1^1) \cup d_{R_2}(x_2^1, y_2^1) \cup \ldots \cup d_{R_m}(x_m^1, y_m^1) =$$
$$r_1(x_1^1 - y_1^1) \cup r_2(x_2^1 - y_2^1) \cup \ldots \cup r_m(x_m^1 - y_m^1)$$
*for every $x_i^1, y_i^1 \in V^{n_i}$, $i = 1, 2, 3, \ldots, m$ (Here $d_{R_i}(x_i^1, y_i^1) = r_i(x_i^1 - y_i^1)$ for every $x_i^1, y_i^1 \in V^{n_i}$; for $i = 1, 2, 3, \ldots, m$).*

**DEFINITION 3.4:** *A linear m-space $V^{n_1} \cup V^{n_2} \cup \ldots \cup V^{n_m}$ over $GF(2^N)$, $N > 1$ of m-dimension $n_1 \cup n_2 \cup \ldots \cup n_m$ such that $n_i \leq N$ for $i = 1, 2, 3, \ldots, m$ equipped with the m-rank m-metric is defined as the m-rank m-space.*

**DEFINITION 3.5:** *A m-rank m-distance RD m-code of m-length $n_1 \cup n_2 \cup \ldots \cup n_m$ over $GF(2^N)$ is a m-subset of the m-rank m-space $V^{n_1} \cup V^{n_2} \cup \ldots \cup V^{n_m}$ over $GF(2^N)$.*

**DEFINITION 3.6:** *A linear $[n_1, k_1] \cup [n_2, k_2] \cup \ldots \cup [n_m, k_m]$ RD-m-code is a linear m-subspace of m-dimension $k_1 \cup k_2 \cup \ldots \cup k_m$ in the m-rank m-space $V^{n_1} \cup V^{n_2} \cup \ldots \cup V^{n_m}$. By $C_1[n_1, k_1] \cup C_2[n_2, k_2] \cup \ldots \cup C_m[n_m, k_m]$ we denote a linear $[n_1, k_1] \cup [n_2, k_2] \cup \ldots \cup [n_m, k_m]$ RD m-code.*

We can equivalently define a RD m-code as follows:

**DEFINITION 3.7:** *Let $V^{n_1}, V^{n_2}, \ldots, V^{n_m}$ be rank spaces $n_i \neq n_j$ if $i \neq j$ over $GF(2^N)$, $N > 1$. Suppose $P_i \subseteq V^{n_i}$, $i = 1, 2, 3, \ldots, m$ be subset of the rank spaces over $GF(2^N)$. Then $P_1 \cup P_2 \cup \ldots \cup P_m \subseteq V^{n_1} \cup V^{n_2} \cup \ldots \cup V^{n_m}$ is a rank distance m-code of m-length $(n_1, n_2, \ldots, n_m)$ over $GF(2^N)$.*

**DEFINITION 3.8:** *A generator m-matrix of a linear $[n_1, k_1] \cup [n_2, k_2] \cup \ldots \cup [n_m, k_m]$ RD-m-code $C_1 \cup C_2 \cup \ldots \cup C_m$ is a $k_1 \times n_1 \cup k_2 \times n_2 \cup \ldots \cup k_m \times n_m$, m-matrix over $GF(2^N)$ whose m-rows form a m-basis for $C_1 \cup C_2 \cup \ldots \cup C_m$.*

*A generator m-matrix $G = G_1 \cup G_2 \cup \ldots \cup G_m$ of a linear RD m-code $C_1[n_1, k_1] \cup C_2[n_2, k_2] \cup \ldots \cup C_m[n_m, k_m]$ can be brought into the form $G = G_1 \cup G_2 \cup \ldots \cup G_m = [I_{k_1}, A_{k_1 \times n_1 - k_1}] \cup [I_{k_2}, A_{k_2 \times n_2 - k_2}] \cup \ldots \cup [I_{k_m}, A_{k_m \times n_m - k_m}]$ where $I_{k_1}, I_{k_2}, \ldots, I_{k_m}$ is the identity matrix and $A_{k_i \times n_i - k_i}$ is some matrix over $GF(2^N)$; $i = 1, 2, 3, \ldots, m$; this form of $G = G_1 \cup G_2 \cup \ldots \cup G_m$ is called the standard form.*

**DEFINITION 3.9:** *If $G = G_1 \cup G_2 \cup \ldots \cup G_m$ is a generator m-matrix of $C_1[n_1, k_1] \cup C_2[n_2, k_2] \cup \ldots \cup C_m[n_m, k_m]$ then a m-matrix $H = H_1 \cup H_2 \cup \ldots \cup H_m$ of order $(n_1 - k_1 \times n_1, n_2 - k_2 \times n_2, \ldots, n_m - k_m \times n_m)$ over $GF(2^N)$ such that,*

$$GH^T = (G_1 \cup G_2 \cup \ldots \cup G_m)(H_1 \cup H_2 \cup \ldots \cup H_m)^T$$
$$= (G_1 \cup G_2 \cup \ldots \cup G_m)(H_1^T \cup H_2^T \cup \ldots \cup H_m^T)$$
$$= G_1 H_1^T \cup G_2 H_2^T \cup \ldots \cup G_m H_m^T$$
$$= 0 \cup 0 \cup \ldots \cup 0$$

*is called a parity check m-matrix of $C_1[n_1, k_1] \cup C_2[n_2, k_2] \cup \ldots \cup C_m[n_m, k_m]$. Suppose $C = C_1 \cup C_2 \cup \ldots \cup C_m$ is a linear $[n_1, k_1] \cup [n_2, k_2] \cup \ldots \cup [n_m, k_m]$ RD m-code with $G = G_1 \cup G_2 \cup \ldots \cup G_m$ and $H = H_1 \cup H_2 \cup \ldots \cup H_m$ as its generator and parity check m-matrix respectively, then $C = C_1 \cup C_2 \cup \ldots \cup C_m$ has two representations.*

(a) $C = C_1 \cup C_2 \cup \ldots \cup C_m$ is a row m-space of $G = G_1 \cup G_2 \cup \ldots \cup G_m$ (i.e., $C_i$ is the row space of $G_i$ for $i = 1, 2, 3, \ldots, m$).

(b) $C = C_1 \cup C_2 \cup \ldots \cup C_m$ is the solution m-space of $H = H_1 \cup H_2 \cup \ldots \cup H_m$ (i.e.; $C_i$ is the solution space of $H_i$ for $i = 1, 2, 3, \ldots, m$).

Now we proceed on to define the notion of minimum rank m-distance of the rank distance m-code $C = C_1 \cup C_2 \cup \ldots \cup C_m$.

**DEFINITION 3.10:** *Let $C = C_1 \cup C_2 \cup \ldots \cup C_m$ be a rank distance m-code, the minimum rank m-distance $d = d_1 \cup d_2 \cup \ldots \cup d_m$ is defined by*

$$d_i = \min\left\{ d_{R_i}(x_i, y_i) \,\middle|\, \begin{matrix} x_i, y_i \in C_i \\ x_i \neq y_i \end{matrix} \right\}$$

$i = 1, 2, 3, \ldots, m$. *That is*

$$d = d_1 \cup d_2 \cup \ldots \cup d_m$$
$$= \min\{ r_1(x_1)/x_1 \in C_1 \text{ and } x_1 \neq 0 \} \cup$$
$$\min\{ r_2(x_2)/x_2 \in C_2 \text{ and } x_2 \neq 0 \} \cup \ldots$$
$$\cup \min\{ r_m(x_m)/x_m \in C_m \text{ and } x_m \neq 0 \}.$$

*If an RD m-code $C = C_1 \cup C_2 \cup \ldots \cup C_m$ has the minimum rank m-distance $d = d_1 \cup d_2 \cup \ldots \cup d_m$ then it can correct all m-errors $e = e_1 \cup e_2 \cup \ldots \cup e_m \in F_{q^N}^{n_1} \cup F_{q^N}^{n_2} \cup \ldots \cup F_{q^N}^{n_m}$ with m-rank*

$$r(e) = (r_1 \cup r_2 \cup \ldots \cup r_m)(e_1 \cup e_2 \cup \ldots \cup e_m)$$
$$= r_1(e_1) \cup r_2(e_2) \cup \ldots \cup r_m(e_m)$$
$$\leq \left\lfloor \frac{d_1 - 1}{2} \right\rfloor \cup \left\lfloor \frac{d_2 - 1}{2} \right\rfloor \cup \ldots \cup \left\lfloor \frac{d_n - 1}{2} \right\rfloor.$$

*Let $C = C_1 \cup C_2 \cup \ldots \cup C_m$ denote an $[n_1, k_1] \cup [n_2, k_2] \cup \ldots \cup [n_m, k_m]$ RD m-code over $F_{q^N}$.*

*A generator m-matrix $G = G_1 \cup G_2 \cup \ldots \cup G_m$ of $C = C_1 \cup C_2 \cup \ldots \cup C_m$ is a $k_1 \times n_1 \cup k_2 \times n_2 \cup \ldots \cup k_m \times n_m$ m-matrix* *with entries from $F_{q^N}$ whose rows form a m-basis for $C = C_1 \cup C_2 \cup \ldots \cup C_m$.*

*Then an $(n_1 - k_1) \times n_1 \cup (n_2 - k_2) \times n_2 \cup \ldots \cup (n_m - k_m) \times n_m$ m-matrix $H = H_1 \cup H_2 \cup \ldots \cup H_m$ with entries from $F_{q^N}$ such that,*

$$GH^T = (G_1 \cup G_2 \cup \ldots \cup G_m)(H_1 \cup H_2 \cup \ldots \cup H_m)^T$$
$$= (G_1 \cup G_2 \cup \ldots \cup G_m)(H_1^T \cup H_2^T \cup \ldots \cup H_m^T)$$
$$= G_1 H_1^T \cup G_2 H_2^T \cup \ldots \cup G_m H_m^T$$
$$= 0 \cup 0 \cup \ldots \cup 0$$

*is called the parity check m-matrix of the RD m-code $C = C_1 \cup C_2 \cup \ldots \cup C_m$.*

The result analogous to Singleton–Style bound in case of RD m-code is given in the following.

*Result:* (Singleton–Style bound). The minimum rank m-distance $d = d_1 \cup d_2 \cup \ldots \cup d_m$ of any linear $[n_1, k_1] \cup [n_2, k_2] \cup \ldots \cup [n_m, k_m]$ RD m-code $C = C_1 \cup C_2 \cup \ldots \cup C_m \subseteq F_{q^N}^{n_1} \cup F_{q^N}^{n_2} \cup \ldots \cup F_{q^N}^{n_m}$ satisfies the following bounds $d = d_1 \cup d_2 \cup \ldots \cup d_m \leq n_1 - k_1 + 1 \cup n_2 - k_2 + 1 \cup \ldots \cup n_m - k_m + 1$.

Based on this notion we now proceed on to define the new notion of Maximum Rank Distance m-codes.

**DEFINITION 3.11:** *An $[n_1, k_1, d_1] \cup [n_2, k_2, d_2] \cup \ldots \cup [n_m, k_m, d_m]$ RD-m-code $C = C_1 \cup C_2 \cup \ldots \cup C_m$ is called a Maximum Rank Distance (MRD) m-code if the Singleton Style bound is reached that is $d = d_1 \cup d_2 \cup \ldots \cup d_m = n_1 - k_1 + 1 \cup n_2 - k_2 + 1 \cup \ldots \cup n_m - k_m + 1$.*

Now we proceed on to briefly give the construction of MRD m-code.

Let $[s] = [s_1] \cup [s_2] \cup \ldots \cup [s_m] = q^{s_1} \cup q^{s_2} \cup \ldots \cup q^{s_m}$ for any m integers $s_1, s_2, \ldots, s_m$. Let $\{g_1, \ldots, g_{n_1}\} \cup \{h_1, \ldots, h_{n_2}\}$

$\cup ... \cup \{p_1, ..., p_{n_m}\}$ be any m-set of elements in $F_{q^N}$ that are linearly independent over $F_q$. A generator m-matrix $G = G_1 \cup G_2 \cup ... \cup G_m$ of an MRD m-code $C = C_1 \cup C_2 \cup ... \cup C_m$ is defined by $G = G_1 \cup G_2 \cup ... \cup G_m$

$$= \begin{bmatrix} g_1 & g_2 & \cdots & g_{n_1} \\ g_1^{[1]} & g_2^{[1]} & \cdots & g_{n_1}^{[1]} \\ g_1^{[2]} & g_2^{[2]} & \cdots & g_{n_1}^{[2]} \\ \vdots & \vdots & & \vdots \\ g_1^{[k_1-1]} & g_2^{[k_1-1]} & \cdots & g_{n_1}^{[k_1-1]} \end{bmatrix} \cup \begin{bmatrix} h_1 & h_2 & \cdots & h_{n_2} \\ h_1^{[1]} & h_2^{[1]} & \cdots & h_{n_2}^{[1]} \\ h_1^{[2]} & h_2^{[2]} & \cdots & h_{n_2}^{[2]} \\ \vdots & \vdots & & \vdots \\ h_1^{[k_2-1]} & h_2^{[k_2-1]} & \cdots & h_{n_2}^{[k_2-1]} \end{bmatrix}$$

$$\cup ... \cup \begin{bmatrix} p_1 & p_2 & \cdots & p_{n_m} \\ p_1^{[1]} & p_2^{[1]} & \cdots & p_{n_m}^{[1]} \\ p_1^{[2]} & p_2^{[2]} & \cdots & p_{n_m}^{[2]} \\ \vdots & \vdots & & \vdots \\ p_1^{[k_m-1]} & p_2^{[k_m-1]} & \cdots & p_{n_m}^{[k_m-1]} \end{bmatrix}.$$

It can be easily proved that the m-code $C = C_1 \cup C_2 \cup ... \cup C_m$ given by the generator m-matrix $G = G_1 \cup G_2 \cup ... \cup G_m$ has the rank m-distance $d = d_1 \cup d_2 \cup ... \cup d_m$. Any m-matrix of the above form is called a Frobenius m-matrix with generating m-vector

$$g_c = g_{c_1} \cup g_{c_2} \cup ... \cup g_{c_m}$$
$$= (g_1, g_2, ..., g_{n_1}) \cup (h_1, h_2, ..., h_{n_2}) \cup ... \cup (p_1, p_2, ..., p_{m_n}).$$

Interested reader can prove the following theorem:

**THEOREM 3.1:** *Let $C[n, k] = C_1[n_1, k_1] \cup C_2[n_2, k_2] \cup ... \cup C_m[n_m, k_m]$ be the linear $[n_1, k_1, d_1] \cup [n_2, k_2, d_2] \cup ... \cup [n_m, k_m, d_m]$ MRD m-code with $d_i = 2t_i + 1$ for $i = 1, 2, 3, ..., m$. Then $C[n, k] = C_1[n_1, k_1] \cup C_2[n_2, k_2] \cup ... \cup C_m[n_m, k_m]$ m-code* *corrects all m-errors of m-rank atmost $t = t_1 \cup t_2 \cup ... \cup t_m$ and detects all m-errors of m-rank greater than $t = t_1 \cup t_2 \cup ... \cup t_m$.*
*Consider the Galois field $GF(2^N)$ ; $N > 1$. An element*
$$\alpha = \alpha_1 \cup \alpha_2 \cup ... \cup \alpha_m \in \underbrace{GF(2^N) \cup GF(2^N) \cup ... \cup GF(2^N)}_{m-times}$$
*can be denoted by m-N-tuple*

$$(a_0^1, a_1^1, ..., a_{N-1}^1) \cup (a_0^2, a_1^2, ..., a_{N-1}^2) \cup ... \cup (a_0^m, a_1^m, ..., a_{N-1}^m)$$

*as well as by the m-polynomial*
$$a_0^1 + a_1^1 x + ... + a_{N-1}^1 x^{N-1} \cup a_0^2 + a_1^2 x + ... + a_{N-1}^2 x^{N-1}$$
$$\cup ... \cup a_0^m + a_1^m x + ... + a_{N-1}^m x^{N-1}$$
*over GF(2)*

We now proceed on to define the new notion of circulant m-transpose.

**DEFINITION 3.12:** *The circulant m-transpose*
$$T_C = T_{C_1}^1 \cup T_{C_2}^2 \cup ... \cup T_{C_m}^m$$
*of a m-vector*
$$\alpha = \alpha_1 \cup \alpha_2 \cup ... \cup \alpha_m =$$
$$(a_0^1, a_1^1, ..., a_{N-1}^1) \cup (a_0^2, a_1^2, ..., a_{N-1}^2) \cup ... \cup (a_0^m, a_1^m, ..., a_{N-1}^m)$$

$\in GF(2^N)$ *is defined as,*
$$\alpha^{T_C} = \alpha_1^{T_{C_1}^1} \cup \alpha_2^{T_{C_2}^2} \cup ... \cup \alpha_m^{T_{C_m}^m} =$$
$$(a_0^1, a_1^1, ..., a_{n_1}^1) \cup (a_0^2, a_1^2, ..., a_{n_2}^2) \cup ... \cup (a_0^m, a_1^m, ..., a_{n_m}^m)$$
*If $\alpha = \alpha_1 \cup \alpha_2 \cup ... \cup \alpha_m \in \underbrace{GF(2^N) \cup GF(2^N) \cup ... \cup GF(2^N)}_{m-times}$*
*has the m-polynomial representation*
$$a_0^1 + a_1^1 x + ... + a_{N-1}^1 x^{N-1} \cup a_0^2 + a_1^2 x + ... + a_{N-1}^2 x^{N-1}$$
$$\cup ... \cup a_0^m + a_1^m x + ... + a_{N-1}^m x^{N-1}$$
*in*
$$\frac{GF(2)[x]}{(x^N + 1)} \cup \frac{GF(2)[x]}{(x^N + 1)} \cup ... \cup \frac{GF(2)[x]}{(x^N + 1)}$$

*then by $\alpha = \alpha_1 \cup \alpha_2 \cup ... \cup \alpha_m$, we denote the m-vector corresponding to the m-polynomial*

$$[a_0^1 + a_1^1 x + ... + a_{N-1}^1 x^{N-1} \cdot x^i] \, mod(\, x^N + 1\,) \cup$$
$$[a_0^2 + a_1^2 x + ... + a_{N-1}^2 x^{N-1} \cdot x^i] \, mod(\, x^N + 1\,) \cup ... \cup$$
$$[a_0^m + a_1^m x + ... + a_{N-1}^m x^{N-1} \cdot x^i] \, mod(\, x^N + 1\,)$$

*for $i = 0, 1, 2, 3, ..., N-1$ ($\alpha = \alpha_1 \cup \alpha_2 \cup ... \cup \alpha_m = \alpha_o$).*

We now proceed onto define the m-word generated by $\alpha = \alpha_1 \cup \alpha_2 \cup ... \cup \alpha_m$.

**DEFINITION 3.13:** *Let $f = f_1 \cup f_2 \cup ... \cup f_m$:*

$$\underbrace{GF(2^N) \cup GF(2^N) \cup ... \cup GF(2^N)}_{m-times} \longrightarrow$$
$$\underbrace{[GF(2^N)]^N \cup [GF(2^N)]^N \cup ... \cup [GF(2^N)]^N}_{m-times}$$

*be defined as*

$$f(\alpha) = f_1(\alpha_1) \cup f_2(\alpha_2) \cup ... \cup f_m(\alpha_m)$$
$$= (\alpha_{10}^{T_{C_1}^1}, \alpha_{11}^{T_{C_1}^1}, ..., \alpha_{1N-1}^{T_{C_1}^1}) \cup (\alpha_{20}^{T_{C_2}^2}, \alpha_{21}^{T_{C_2}^2}, ..., \alpha_{2N-1}^{T_{C_2}^2}) \cup ...$$
$$\cup (\alpha_{m0}^{T_{C_m}^m}, \alpha_{m1}^{T_{C_m}^m}, ..., \alpha_{mN-1}^{T_{C_m}^m}).$$

*We call $f(\alpha) = f_1(\alpha_1) \cup f_2(\alpha_2) \cup ... \cup f_m(\alpha_m)$ as the m-code word generated by $\alpha = \alpha_1 \cup \alpha_2 \cup ... \cup \alpha_m$.*

We analogous to the definition of Macwilliams and Solane define circulant m-matrix associated with a m-vector in $GF(2^N) \cup GF(2^N) \cup ... \cup GF(2^N)$.

**DEFINITION 3.14:** *A m-matrix of the form*

$$\begin{bmatrix} a_0^1 & a_1^1 & ... & a_{N-1}^1 \\ a_{N-1}^1 & a_0^1 & ... & a_{N-2}^1 \\ \vdots & \vdots & & \vdots \\ a_1^1 & a_2^1 & ... & a_0^1 \end{bmatrix} \cup \begin{bmatrix} a_0^2 & a_1^2 & ... & a_{N-1}^2 \\ a_{N-1}^2 & a_0^2 & ... & a_{N-2}^2 \\ \vdots & \vdots & & \vdots \\ a_1^2 & a_2^2 & ... & a_0^2 \end{bmatrix}$$

$$\cup \ldots \cup \begin{bmatrix} a_0^m & a_1^m & \ldots & a_{N-1}^m \\ a_{N-1}^m & a_0^m & \ldots & a_{N-2}^m \\ \vdots & \vdots & & \vdots \\ a_1^m & a_2^m & \ldots & a_0^m \end{bmatrix}$$

*is called the circulant m-matrix associated with the m-vector*
$(a_0^1, a_1^1, \ldots, a_{N-1}^1) \cup (a_0^2, a_1^2, \ldots, a_{N-1}^2) \cup \ldots \cup (a_0^m, a_1^m, \ldots, a_{N-1}^m) \in$
*$GF(2^N) \cup GF(2^N) \cup \ldots \cup GF(2^N)$.*

Thus with each $\alpha = \alpha_1 \cup \alpha_2 \cup \ldots \cup \alpha_m \in GF(2^N) \cup GF(2^N) \cup \ldots \cup GF(2^N)$ we can associate a circulant m-matrix whose $i^{th}$ m-columns represents $\alpha_{1i}^{T_{c_1}^1} \cup \alpha_{2i}^{T_{c_2}^2} \cup \ldots \cup \alpha_{mi}^{T_{c_m}^m}$ ; i = 0, 1, 2, …, N–1. $f = f_1 \cup f_2 \cup \ldots \cup f_m$ is nothing but a m-mapping of $GF(2^N) \cup GF(2^N) \cup \ldots \cup GF(2^N)$ onto the pseudo false m-algebra of all N × N circulant m-matrices over GF(2). Denote the m-space of $f(GF(2^N)) = f_1(GF(2^N)) \cup f_2(GF(2^N)) \cup \ldots \cup f_m(GF(2^N))$ by

$$\underbrace{V^N \cup V^N \cup \ldots \cup V^N}_{m-times}.$$

We define the m-norm of a m-word
$$v = v_1 \cup v_2 \cup \ldots \cup v_m \in \underbrace{V^N \cup V^N \cup \ldots \cup V^N}_{m-times}$$

as follows :

**DEFINITION 3.15:** *The m-norm of a m-word $v = v_1 \cup v_2 \cup \ldots \cup v_m \in V^N \cup V^N \cup \ldots \cup V^N$ is defined as the m-rank of $v = v_1 \cup v_2 \cup \ldots \cup v_m$ over GF(2) [By considering it as a circulant m-matrix over GF(2)].*

We denote the m-norm of $v = v_1 \cup v_2 \cup \ldots \cup v_m$ by $r(v) = r_1(v_1) \cup r_2(v_2) \cup \ldots \cup r_m(v_m)$, we prove the following theorem:

**THEOREM 3.2:** *Suppose*
$$\alpha = \alpha_1 \cup \alpha_2 \cup \ldots \cup \alpha_m \in \underbrace{GF(2^N) \cup GF(2^N) \cup \ldots \cup GF(2^N)}_{m-times}$$

*has the m-polynomial representation $g(x) = g_1(x) \cup g_2(x) \cup \ldots \cup g_m(x)$ over GF(2) such that $\gcd(g_i(x), x^N + 1)$ has degree $N - k_i$ for $i = 1, 2, 3, \ldots, m$; $1 \leq k_1, k_2, \ldots, k_m \leq N$. Then the m-norm of the m-word generated by $\alpha = \alpha_1 \cup \alpha_2 \cup \ldots \cup \alpha_m$ is $k_1 \cup k_2 \cup \ldots \cup k_m$.*

*Proof:* We know the m-norm of the m-word generated by $\alpha = \alpha_1 \cup \alpha_2 \cup \ldots \cup \alpha_m$ is the m-rank of the circulant m-matrix

$$(\alpha_{10}^{T_{C_1}^1}, \alpha_{11}^{T_{C_1}^1}, \ldots, \alpha_{1N-1}^{T_{C_1}^1}) \cup (\alpha_{20}^{T_{C_2}^2}, \alpha_{21}^{T_{C_2}^2}, \ldots, \alpha_{2N-1}^{T_{C_2}^2}) \cup \ldots \cup$$
$$(\alpha_{m0}^{T_{C_m}^m}, \alpha_{m1}^{T_{C_m}^m}, \ldots, \alpha_{mN-1}^{T_{C_m}^m})$$

where $\alpha_i^{T_C} = \alpha_{1i}^{T_{C_1}^1} \cup \alpha_{2i}^{T_{C_2}^2} \cup \ldots \cup \alpha_{mi}^{T_{C_m}^m}$ represents a m-polynomial

$$[x^i g_1(x)] \bmod (x^N + 1) \cup [x^i g_2(x)] \bmod (x^N + 1)$$
$$\cup \ldots \cup [x^i g_m(x)] \bmod (x^N + 1)$$

over GF(2).

Suppose the m-GCD $\{(g_1(x), x^N + 1) \cup (g_2(x), x^N + 1) \cup \ldots \cup (g_m(x), x^N + 1)\}$ has m-degree $N - k_1 \cup N - k_2 \cup \ldots \cup N - k_m$. To prove that the m-word generated by $\alpha = \alpha_1 \cup \alpha_2 \cup \ldots \cup \alpha_m$ has m-rank $k_1 \cup k_2 \cup \ldots \cup k_m$. It is enough to prove that the m-space generated by the N-polynomials $\{g_1(x) \bmod(x^N + 1), x \cdot g_1(x) \bmod(x^N + 1), \ldots, x^{N-1} \cdot g_1(x) \bmod(x^N + 1)\} \cup \{g_2(x) \bmod(x^N + 1), x \cdot g_2(x) \bmod(x^N + 1), \ldots, x^{N-1} \cdot g_2(x) \bmod(x^N + 1)\} \cup \ldots \cup \{g_m(x) \bmod(x^N + 1), x \cdot g_m(x) \bmod(x^N + 1), \ldots, x^{N-1} \cdot g_m(x) \bmod(x^N + 1)\}$ has m-dimension $k_1 \cup k_2 \cup \ldots \cup k_m$. We will prove that the m-set of $k_1 \cup k_2 \cup \ldots \cup k_m$, m-polynomials $\{g_1(x) \bmod x^N + 1, x \cdot g_1(x) \bmod x^N + 1, \ldots, x^{N-1} \cdot g_1(x) \bmod x^N + 1\} \cup \{g_2(x) \bmod x^N + 1, x \cdot g_2(x) \bmod x^N + 1, \ldots, x^{N-1} \cdot g_2(x) \bmod x^N + 1\} \cup \ldots \cup \{g_m(x) \bmod x^N + 1, x \cdot g_m(x) \bmod x^N + 1, \ldots, x^{N-1} \cdot g_m(x) \bmod x^N + 1\}$ forms a m-basis for the m-space. If possible let,

$$a_0^1(g_1(x)) + a_1^1 x(g_1(x)) + \ldots + a_{k_1-1}^1(x^{k_1-1} g_1(x)) \cup$$
$$a_0^2(g_2(x)) + a_1^2 x(g_2(x)) + \ldots + a_{k_2-1}^2(x^{k_2-1} g_2(x)) \cup \ldots \cup$$
$$a_0^m(g_m(x)) + a_1^m x(g_m(x)) + \ldots + a_{k_m-1}^m(x^{k_m-1} g_m(x))$$

$$= 0 \cup 0 \cup \ldots \cup 0 \ (\text{mod } x^N + 1);$$

where $a_{j_i}^i \in GF(2)$, $1 \leq j_i \leq k_i - 1$ and $i = 1, 2, \ldots, m$.

This implies

$$\underbrace{x^{N+1} \cup x^{N+1} \cup \ldots \cup x^{N+1}}_{m-\text{times}}$$

m-divides

$$(a_0^1 + a_1^1 x + \ldots + a_{k_1-1}^1 x^{k_1-1})g_1(x) \cup$$
$$(a_0^2 + a_1^2 x + \ldots + a_{k_2-1}^2 x^{k_2-1})g_2(x) \cup \ldots \cup$$
$$(a_0^m + a_1^m x + \ldots + a_{k_m-1}^m x^{k_m-1})g_m(x).$$

Now if $g_1(x) \cup g_2(x) \cup \ldots \cup g_m(x) = p_1(x)a_1(x) \cup p_2(x)a_2(x) \cup \ldots \cup p_m(x)a_m(x)$ where $p_i(x)$ is the $\gcd(g_i(x), x^N + 1)$; $i = 1, 2, \ldots, m$, then $(a_i(x), x^N + 1) = 1$. Thus $x^N + 1$ m-divides

$$(a_0^1 + a_1^1 x + \ldots + a_{k_1-1}^1 x^{k_1-1})g_1(x) \cup \ldots \cup$$
$$(a_0^m + a_1^m x + \ldots + a_{k_m-1}^m x^{k_m-1})g_m(x)$$

implies the m-quotient

$$\frac{(x^N + 1)}{p_1(x)} \cup \ldots \cup \frac{(x^N + 1)}{p_m(x)}$$

m-divides

$$(a_0^1 + a_1^1 x + \ldots + a_{k_1-1}^1 x^{k_1-1})a_1(x) \cup$$
$$(a_0^2 + a_1^2 x + \ldots + a_{k_2-1}^2 x^{k_2-1})a_2(x) \cup \ldots \cup$$
$$(a_0^m + a_1^m x + \ldots + a_{k_m-1}^1 x^{k_m-1})a_m(x).$$

That is

$$\left[\frac{(x^N + 1)}{p_1(x)}\right] \cup \left[\frac{(x^N + 1)}{p_2(x)}\right] \cup \ldots \cup \left[\frac{(x^N + 1)}{p_m(x)}\right]$$

m-divides

$$(a_0^1 + a_1^1 x + \ldots + a_{k_1-1}^1 x^{k_1-1}) \cup$$
$$(a_0^2 + a_1^2 x + \ldots + a_{k_2-1}^2 x^{k_2-1}) \cup \ldots \cup$$
$$(a_0^m + a_1^m x + \ldots + a_{k_m-1}^1 x^{k_m-1})$$

which is a contradiction, as

$$\left[\frac{(x^N+1)}{p_1(x)}\right] \cup \left[\frac{(x^N+1)}{p_2(x)}\right] \cup ... \cup \left[\frac{(x^N+1)}{p_m(x)}\right]$$

has m-degree $(k_1, k_2, ..., k_m)$ where as the m-polynomial

$$(a_0^1 + a_1^1 x + ... + a_{k_1-1}^1 x^{k_1-1}) \cup$$
$$(a_0^2 + a_1^2 x + ... + a_{k_2-1}^2 x^{k_2-1}) \cup ... \cup$$
$$(a_0^m + a_1^m x + ... + a_{k_m-1}^1 x^{k_m-1})$$

has m-degree atmost $((k_1 - 1), (k_2 - 1), ..., (k_m - 1))$. Hence the m-polynomials $\{g_1(x) \bmod(x^N + 1), x \cdot g_1(x) \bmod(x^N + 1), ..., x^{k_1-1} \cdot g_1(x) \bmod(x^N + 1)\} \cup \{g_2(x) \bmod(x^N + 1), x \cdot g_2(x) \bmod(x^N + 1), ..., x^{k_2-1} \cdot g_2(x) \bmod(x^N + 1)\} \cup ... \cup \{g_m(x) \bmod(x^N + 1), x \cdot g_m(x) \bmod(x^N + 1), ..., x^{k_m-1} \cdot g_m(x) \bmod(x^N + 1)\}$ are m-linearly independent over GF(2).

We will prove, $\{g_1(x) \bmod(x^N + 1), x \cdot g_1(x) \bmod(x^N + 1), ..., x^{k_1-1} \cdot g_1(x) \bmod(x^N + 1)\} \cup \{g_2(x) \bmod(x^N + 1), x \cdot g_2(x) \bmod(x^N + 1), ..., x^{k_2-1} \cdot g_2(x) \bmod(x^N + 1)\} \cup ... \cup \{g_m(x) \bmod(x^N + 1), x \cdot g_m(x) \bmod(x^N + 1), ..., x^{k_m-1} \cdot g_m(x) \bmod(x^N + 1)\}$ generate the m-space.

For this it is enough to prove that $x^i g_1(x) \cup x^i g_2(x) \cup ... \cup x^i g_m(x)$ is a linear combination of these m-polynomials for $k_j \leq i \leq N - 1; j = 1, 2, 3, ..., m$.

$$\underbrace{x^N + 1 \cup x^N + 1 \cup ... \cup x^N + 1}_{m-\text{times}}$$

$$= p_1(x)b_1(x) \cup p_2(x)b_2(x) \cup ... \cup p_m(x)b_m(x)$$

where $b_i(x) = b_0^i + b_1^i x + ... + b_{k_i}^i x^{k_i}$; $i = 1, 2, 3, ..., m$. (Note that $b_0^i = b_{k_i}^i = 1$ since $b_i(x)$ divides $x^N + 1$, $i = 1, 2, 3, ..., m$). Also we have $g_i(x) = p_i(x)a_i(x)$ for $i = 1, 2, 3, ..., m$. Thus

$$x^N + 1 = \left( \frac{g_i(x) b_i(x)}{a_i(x)} \right)$$

for i = 1, 2, 3, …, m. That is

$$\frac{g_i(x)(b_0^i + b_1^i x + \ldots + b_{k_i}^i x^{k_i})}{a_i(x)} = 0 \bmod(x^N + 1)$$

(i = 1, 2, 3, …, m) that is

$$\frac{g_i(x)(b_0^i + b_1^i x + \ldots + b_{k_i-1}^i x^{k_i-1})}{a_i(x)} = \left[ \frac{g_i(x) \cdot x^{k_i}}{a_i(x)} \right] \bmod(x^N + 1)$$

true for i = 1, 2, 3, ..., m. Hence
$$x^{k_i} g_i(x) = (b_0^i g_i(x) + b_1^i x g_i(x) + \ldots + b_{k_i-1}^i x^{k_i-1} g_i(x)) \bmod(x^N + 1)$$
a linear combination of $\{g_i(x) \bmod(x^N + 1), [x \cdot g_i(x)] \bmod(x^N + 1), \ldots, x^{k_i-1} \cdot g_m(x) \bmod(x^N + 1)\}$ over GF(2).

This is true for each i, for i = 1, 2, 3, …, m. Now it can be easily proved that $x^i g_1(x) \cup x^i g_2(x) \cup \ldots \cup x^i g_m(x)$ is a m-linear combination of $\{g_1(x) \bmod(x^N + 1), x \cdot g_1(x) \bmod(x^N + 1), \ldots, x^{k_1-1} \cdot g_1(x) \bmod(x^N + 1)\} \cup \{g_2(x) \bmod(x^N + 1), x \cdot g_2(x) \bmod(x^N + 1), \ldots, x^{k_2-1} \cdot g_2(x) \bmod(x^N + 1)\} \cup \ldots \cup \{g_m(x) \bmod(x^N + 1), x \cdot g_m(x) \bmod(x^N + 1), \ldots, x^{k_m-1} \cdot g_m(x) \bmod(x^N + 1)\}$ for $i > k_j$; j = 1, 2, …, m.

Hence the m-space generated by the m-polynomials $\{g_1(x) \bmod(x^N + 1), x \cdot g_1(x) \bmod(x^N + 1), \ldots, x^{k_1-1} \cdot g_1(x) \bmod(x^N + 1)\} \cup \{g_2(x) \bmod(x^N + 1), x \cdot g_2(x) \bmod(x^N + 1), \ldots, x^{k_2-1} \cdot g_2(x) \bmod(x^N + 1)\} \cup \ldots \cup \{g_m(x) \bmod(x^N + 1), x \cdot g_m(x) \bmod(x^N + 1), \ldots, x^{k_m-1} \cdot g_m(x) \bmod(x^N + 1)\}$ has m-dimension $(k_1, k_2, \ldots, k_m)$. That is the m-rank of a m-word generated by $\alpha = \alpha_1 \cup \alpha_2 \cup \ldots \cup \alpha_m$ is $(k_1, k_2, \ldots, k_m)$.

**COROLLARY 3.1:** *If* $\alpha = \alpha_1 \cup \alpha_2 \cup \ldots \cup \alpha_m \in \underbrace{GF(2^N) \cup \ldots \cup GF(2^N)}_{m-times}$ *then the m-norm of the m-word generated by* $\alpha = \alpha_1 \cup \alpha_2 \cup \ldots \cup \alpha_m$ *is (N, N, …, N) and hence* $f(\alpha) = f_1(\alpha_1) \cup f_2(\alpha_2) \cup \ldots \cup f_m(\alpha_m)$ *is m-invertible. (We say* $f(\alpha)$ is m-invertible if each $f_i(\alpha)$ is invertible for $i = 1, 2, 3, \ldots, m$).

*Proof:* The corollary follows immediately from the theorem since $\gcd(g_i(x), x^N+1) = 1$ has degree 0 for $i = 1, 2, 3, \ldots, m$; hence the m-rank of $f(\alpha) = f_1(\alpha_1) \cup f_2(\alpha_2) \cup \ldots \cup f_m(\alpha_m)$ is $(N, N, \ldots, N)$.

**DEFINITION 3.16:** *The m-distance between two m-words $u, v \in V^N \cup \ldots \cup V^N$ is defined as, $d(u, v) = d_1(u_1, v_1) \cup \ldots \cup d_m(u_m, v_m) = r_1(u_1 + v_1) \cup \ldots \cup r_m(u_m + v_m)$ where $u = u_1 \cup u_2 \cup \ldots \cup u_m$ and $v = v_1 \cup v_2 \cup \ldots \cup v_m$.*

**DEFINITION 3.17:** *Let $C = C_1 \cup C_2 \cup \ldots \cup C_m$ be a circulant rank m-code of m-length $N_1 \cup N_2 \cup \ldots \cup N_m$ which is a m-subspace of $V^{N_1} \cup V^{N_2} \cup \ldots \cup V^{N_m}$ equipped with the m-distance m-function $d_1(u_1, v_1) \cup d_2(u_2, v_2) \cup \ldots \cup d_m(u_m, v_m) = r_1(u_1 + v_1) \cup r_2(u_2 + v_2) \cup \ldots \cup r_m(u_m + v_m)$ where $V^{N_1}, V^{N_2}, \ldots, V^{N_m}$ are rank spaces defined over $GF(2^N)$ with $N_i \neq N_j$ if $i \neq j$. $C = C_1 \cup C_2 \cup \ldots \cup C_m$ is defined as the circulant m-code of m-length $N_1 \cup N_2 \cup \ldots \cup N_m$ defined as a m-subspace of $V^{N_1} \cup V^{N_2} \cup \ldots \cup V^{N_m}$ equipped with the m-distance m-function.*

**DEFINITION 3.18:** *A circulant m-rank m-code of m-length $N_1 \cup N_2 \cup \ldots \cup N_m$ is called m-cyclic if whenever $(v_1^1, \ldots, v_{N_1}^1) \cup (v_1^2, \ldots, v_{N_2}^2) \cup \ldots \cup (v_1^m, \ldots, v_{N_m}^m)$ is a m-code word then it implies $(v_2^1, v_3^1, \ldots, v_{N_1}^1, v_1^1) \cup (v_2^2, v_3^2, \ldots, v_{N_2}^2, v_1^2) \cup \ldots \cup (v_2^m, v_3^m, \ldots, v_{N_m}^m, v_1^m)$ is also a m-code word.*

Now we proceed on to define quasi MRD-m-codes.

**DEFINITION 3.19:** *Let $C = C_1 \cup C_2 \cup \ldots \cup C_m$ be a RD rank m-code where each $C_i \neq C_j$ if $i \neq j$. If some of the $C_i$'s are MRD codes and others are RD codes then we call $C = C_1 \cup C_2 \cup \ldots \cup C_m$ to be a quasi MRD m-code.*

*Note:* If $r_1$ are MRD codes and $r_2$ are RD codes $r_1 + r_2 = m$, ($r_1 \geq 1$, $r_2 \geq 1$). Then we call C to be a quasi ($r_1$, $r_2$) MRD m-code. We can say $C_1, \ldots, C_{r_1}$ are $C_i[n_i, k_i]$ RD-codes; i = 1, 2, 3, …, $r_1$ and $C_j[n_j, k_j, d_j]$ are MRD codes for j = 1, 2, 3, …, $r_2$ with $r_1 + r_2 = m$.

Thus
$$C = C_1[n_1, k_1] \cup \ldots \cup C_{r_1}[n_{r_1}, k_{r_1}]$$
$$\cup C_1[n_1, k_1, d_1] \cup \ldots \cup C_{r_2}[n_{r_2}, k_{r_2}, d_{r_2}]$$

is a quasi ($r_1$, $r_2$) MRD m-code. Any m-code word in C would be of the form $C = C_1 \cup C_2 \cup \ldots \cup C_m$. The m-codes can be used in multi channel simultaneously when one needs both MRD codes and RD codes. This will be useful in applications in such type of channels.

We proceed on to define the notion of quasi circulant m-codes of type I.

**DEFINITION 3.20:** *Let $C_1, C_2, \ldots, C_m$ be m distinct codes some circulant rank codes and others linear RD-codes defined over $GF(2^N)$. $C = C_1 \cup C_2 \cup \ldots \cup C_m$ is defined as the quasi circulant m-code of type I. If in the quasi circulant m-code of type I some of the $C_i$'s are MRD codes i.e., $C_1, C_2, \ldots, C_m$ is a collection of RD codes, MRD codes and circulant rank codes then we define $C = C_1 \cup C_2 \cup \ldots \cup C_m$ to be a quasi circulant m-code of type II.*

We can define also mixed quasi circulant rank m-codes.

**DEFINITION 3.21:** *Let $C_1, C_2, \ldots, C_m$ be a collection of m-codes, all of them distinct $C_i \neq C_j$ if $i \neq j$ and $C_i \not\subset C_j$ or $C_j \not\subset C_i$ if $i \neq j$. If this collection of codes $C_1, C_2, \ldots, C_m$ are such that some of them are RD-codes, some MRD codes, some cyclic circulant rank codes and some only circulant codes then we define $C = C_1 \cup C_2 \cup \ldots \cup C_m$ to be a mixed quasi circulant rank m-code.*

**DEFINITION 3.22:** *If $C = C_1 \cup C_2 \cup \ldots \cup C_m$ be a collection of distinct circulant codes some $C_i$'s are circulant codes and some* *of them are cyclic circulant codes then C is defined to be a mixed circulant m-code.*

These codes will find applications in multi channels which have very high error probability and error correction. These multi channels (m-channels) are such that some of the channels have to work only with circulant codes not cyclic circulant codes and some only with cyclic circulant codes these mixed circulant m-codes will be appropriate.

Now we proceed on to define the notion of almost maximum rank distance m-codes.

**DEFINITION 3.23:** *Let $C_1[n_1, k_1] \cup C_2[n_2, k_2] \cup \ldots \cup C_m[n_m, k_m]$ be a collection of m distinct almost maximum rank distance codes with minimum distances greater than equal to $n_1 - k_1 \cup n_2 - k_2 \cup \ldots \cup n_m - k_m$ defined over $GF(2^N)$. C is defined as the Almost Maximum Distance Rank-m-code or (AMRD-m-code) over $GF(2^N)$. An AMRD m-code is called a AMRD – tricode, if m = 3. An AMRD m-code whose minimum distance is greater than $n_1 - k_1 \cup n_2 - k_2 \cup \ldots \cup n_m - k_m$ is an MRD m-code hence the class of MRD m-codes is a subclass of the class of AMRD m-codes.*

We have an interesting property about the AMRD m-codes.

**THEOREM 3.3:** *When $(n_1 - k_1) \cup (n_2 - k_2) \cup \ldots \cup (n_m - k_m)$ AMRD m-code $C = C_1[n_1, k_1] \cup C_2[n_2, k_2] \cup \ldots \cup C_m[n_m, k_m]$ is such that each $(n_i - k_i)$ is odd for i = 1, 2, 3, …, m; then*
  1. *The error correcting capability of the $[n_1, k_1] \cup [n_2, k_2] \cup \ldots \cup [n_m, k_m]$ AMRD m-code is equal to that of an $[n_1, k_1] \cup [n_2, k_2] \cup \ldots \cup [n_m, k_m]$ MRD m-code.*
  2. *An $[n_1, k_1] \cup [n_2, k_2] \cup \ldots \cup [n_m, k_m]$ AMRD m-code is better than any $[n_1, k_1] \cup [n_2, k_2] \cup \ldots \cup [n_m, k_m]$ m-code in Hamming metric for error correction.*

*Proof: (1)* Suppose $C = C_1 \cup C_2 \cup \ldots \cup C_m$ is a $[n_1, k_1] \cup [n_2, k_2] \cup \ldots \cup [n_m, k_m]$ AMRD m-code such that $(n_1 - k_1) \cup (n_2 -$ $k_2) \cup \ldots \cup (n_m - k_m)$ is an odd m-integer (i.e., $n_i - k_i \neq n_j - k_j$ if $i \neq j$ are odd integer $1 \leq i, j \leq m$). The maximum number of m-errors corrected by $C = C_1 \cup C_2 \cup \ldots \cup C_m$ is given by

$$\frac{(n_1 - k_1 - 1)}{2} \cup \frac{(n_2 - k_2 - 1)}{2} \cup \ldots \cup \frac{(n_m - k_m - 1)}{2}.$$

But

$$\frac{(n_1 - k_1 - 1)}{2} \cup \frac{(n_2 - k_2 - 1)}{2} \cup \ldots \cup \frac{(n_m - k_m - 1)}{2}$$

is equal to the error correcting capability of an $[n_1, k_1] \cup [n_2, k_2] \cup \ldots \cup [n_m, k_m]$; MRD m-code (since $(n_1 - k_1)$, $(n_2 - k_2)$, …, $(n_m - k_m)$ are odd). That is, $(n_1 - k_1) \cup (n_2 - k_2) \cup \ldots \cup (n_m - k_m)$ is said to be m-odd if each $n_i - k_i$ is odd for $i = 1, 2, 3, \ldots, m$. Thus a $[n_1, k_1] \cup [n_2, k_2] \cup \ldots \cup [n_m, k_m]$ AMRD m-code is as good as an $[n_1, k_1] \cup [n_2, k_2] \cup \ldots \cup [n_m, k_m]$ MRD m-code.

*Proof: (2)* Suppose $C = C_1 \cup C_2 \cup \ldots \cup C_m$ is a $[n_1, k_1] \cup [n_2, k_2] \cup \ldots \cup [n_m, k_m]$ AMRD m-code such that $(n_1 - k_1) \cup (n_2 - k_2) \cup \ldots \cup (n_m - k_m)$ are odd; then each m-code word of C can correct $L_{r_1}(n_1) \cup L_{r_2}(n_2) \cup \ldots \cup L_{r_m}(n_m) = L_r(n)$ error m-vectors where

$$r = r_1 \cup r_2 \cup \ldots \cup r_m$$
$$= \frac{(n_1 - k_1 - 1)}{2} \cup \frac{(n_2 - k_2 - 1)}{2} \cup \ldots \cup \frac{(n_m - k_m - 1)}{2}$$

and

$$L_r(n) = L_{r_1}(n_1) \cup L_{r_2}(n_2) \cup \ldots \cup L_{r_m}(n_m)$$
$$= 1 + \sum_{i=1}^{n_1} \begin{bmatrix} n_1 \\ i \end{bmatrix} (2^N - 1) \ldots (2^N - 2^{i-1}) \cup$$
$$1 + \sum_{i=1}^{n_2} \begin{bmatrix} n_2 \\ i \end{bmatrix} (2^N - 1) \ldots (2^N - 2^{i-1}) \cup$$
$$\ldots \cup 1 + \sum_{i=1}^{n_m} \begin{bmatrix} n_m \\ i \end{bmatrix} (2^N - 1) \ldots (2^N - 2^{i-1}).$$

Consider the same $[n_1, k_1] \cup [n_2, k_2] \cup \ldots \cup [n_m, k_m]$ m-code in Hamming metric. Let it be denoted by $D = D_1 \cup D_2 \cup \ldots \cup D_m$ then the minimum m-distance of D is atmost $(n_1 - k_1 + 1) \cup (n_2 - k_2 + 1) \cup \ldots \cup (n_m - k_m + 1)$. The error correcting capability of D is

$$\left\lfloor \frac{n_1 - k_1 + 1 - 1}{2} \right\rfloor \cup \left\lfloor \frac{n_2 - k_2 + 1 - 1}{2} \right\rfloor \cup \ldots \cup \left\lfloor \frac{n_m - k_m + 1 - 1}{2} \right\rfloor$$

$= r_1 \cup r_2 \cup \ldots \cup r_m$ (since $(n_1 - k_1) \cup (n_2 - k_2) \cup \ldots \cup (n_m - k_m)$ are odd). Hence the number of error m-vectors corrected by the m-code word is given by

$$\sum_{i=0}^{r_1} \begin{bmatrix} n_1 \\ i \end{bmatrix}(2^N - 1)^i \cup \sum_{i=0}^{r_2} \begin{bmatrix} n_2 \\ i \end{bmatrix}(2^N - 1)^i \cup \ldots \cup \sum_{i=0}^{r_m} \begin{bmatrix} n_m \\ i \end{bmatrix}(2^N - 1)^i$$

which is clearly less than $L_{r_1}(n_1) \cup L_{r_2}(n_2) \cup \ldots \cup L_{r_m}(n_m)$. Thus the number of error m-vectors that can be corrected by the $[n_1, k_1] \cup [n_2, k_2] \cup \ldots \cup [n_m, k_m]$ AMRD m-code is much greater than that of the same m-code considered in Hamming metric.

For a given m-length $n = n_1 \cup n_2 \cup \ldots \cup n_m$ a single error correcting AMRD m-code is one having m-dimension $(n_1 - 3) \cup (n_2 - 3) \cup \ldots \cup (n_m - 3)$ and the minimum m-distance greater than or equal to $3 \cup 3 \cup \ldots \cup 3$.

We now proceed on to give a characterization of a single error correcting AMRD m-codes in terms of its parity check m-matrices.

The characterization is based on the condition for the minimum distance proved by Gabidulin in [24, 27].

**THEOREM 3.4:** *Let* $H = H_1 \cup H_2 \cup \ldots \cup H_m = (\alpha_{ij}^1) \cup (\alpha_{ij}^2) \cup \ldots \cup (\alpha_{ij}^m)$ *be a* $(3 \times n_1) \cup (3 \times n_2) \cup \ldots \cup (3 \times n_m)$ *m-matrix of m-rank 3 over* $GF(2^N)$; $n_1 \leq N$ *and* $n_2 \leq N$ *which satisfies the following condition.*

*For any two distinct, non empty m-subsets* $P_1, P_2, \ldots, P_m$ *where* $P_1 = P_1^1 \cup P_2^1 \cup \ldots \cup P_m^1$ *and* $P_2 = P_1^2 \cup P_2^2 \cup \ldots \cup P_m^2$ *of* $\{1, 2, 3, \ldots, n_1\}$ *and* $\{1, 2, 3, \ldots, n_2\}$ *respectively; there exists*

$i_1 = i_1^1 \cup i_2^1$, $i_2 = i_1^2 \cup i_2^2$, ..., $i_m = i_1^m \cup i_2^m \in \{1, 2, 3\} \cup \{1, 2, 3\} \cup$ ... $\cup \{1, 2, 3\}$ such that

$$\left(\sum_{j_1^1 \in P_1^1} \alpha_{i_1^1 j_1^1}^1 \cdot \sum_{k_1^1 \in P_2^1} \alpha_{i_2^1 k_1^1}^1\right) \cup \left(\sum_{j_1^2 \in P_1^2} \alpha_{i_1^2 j_1^2}^2 \cdot \sum_{k_1^2 \in P_2^2} \alpha_{i_2^2 k_1^2}^2\right)$$

$$\cup ... \cup \left(\sum_{j_1^m \in P_1^m} \alpha_{i_1^m j_1^m}^m \cdot \sum_{k_1^m \in P_2^m} \alpha_{i_2^m k_1^m}^m\right)$$

$$\neq \left(\sum_{j_1^1 \in P_1^1} \alpha_{i_2^1 j_1^1}^1 \cdot \sum_{k_1^1 \in P_2^1} \alpha_{i_1^1 k_1^1}^1\right) \cup \left(\sum_{j_1^2 \in P_1^2} \alpha_{i_2^2 j_1^2}^2 \cdot \sum_{k_1^2 \in P_2^2} \alpha_{i_1^2 k_1^2}^2\right)$$

$$\cup ... \cup \left(\sum_{j_1^m \in P_1^m} \alpha_{i_2^m j_1^m}^m \cdot \sum_{k_1^m \in P_2^m} \alpha_{i_1^m k_1^m}^m\right).$$

*Then $H = H_1 \cup H_2 \cup ... \cup H_m$ as a parity check m-matrix defines a $(n_1, n_1 - 3) \cup (n_2, n_2 - 3) \cup ... \cup (n_m, n_m - 3)$ single m-error correcting AMRD m-code over GF(2).*

*Proof:* Given $H = H_1 \cup H_2 \cup ... \cup H_m$ is a $(3 \times n_1) \cup (3 \times n_2) \cup ... \cup (3 \times n_m)$ m-matrix of m-rank $3 \cup 3 \cup ... \cup 3$ over $GF(2^N)$, so that $H = H_1 \cup H_2 \cup ... \cup H_m$ as a parity check m-matrix defines a $(n_1, n_1 - 3) \cup (n_2, n_2 - 3) \cup ... \cup (n_m, n_m - 3)$ RD m-code, where

$$C_1 = \{x \in V^{n_1} / xH_1^T = 0\},$$
$$C_2 = \{x \in V^{n_2} / xH_2^T = 0\}, ...,$$
$$C_m = \{x \in V^{n_m} / xH_m^T = 0\}.$$

It remains to prove that the minimum m-distance of $C = C_1 \cup C_2 \cup ... \cup C_m$ is greater than or equal to $3 \cup 3 \cup ... \cup 3$. We will prove that no non zero m-word of $C = C_1 \cup C_2 \cup ... \cup C_m$ has m-rank less than $3 \cup 3 \cup ... \cup 3$.

The proof is by the method of contradiction.

Suppose there exists a non zero m-code word $x = x_1 \cup x_2 \cup ... \cup x_m$ such that $r_1(x_1) \leq 2$, $r_2(x_2) \leq 2$, ..., $r_m(x_m) \leq 2$, then $x = x_1$ $\cup\, x_2 \cup \ldots \cup x_m$ can be written as $x = x_1 \cup x_2 \cup \ldots \cup x_m = (y_1 \cup y_2 \cup \ldots \cup y_m)(M_1 \cup M_2 \cup \ldots \cup M_m)$ where

$$y_1 = (y_1^1, y_2^1),\ y_2 = (y_1^2, y_2^2),\ \ldots,\ y_m = (y_1^m, y_2^m);$$

$y_1^i, y_2^i \in GF(2^N)$; $1 \le i \le m$ and $M = M_1 \cup M_2 \cup \ldots \cup M_m = (m_{ij}^1) \cup (m_{ij}^2) \cup \ldots \cup (m_{ij}^m)$ is a $(2 \times n_1) \cup (2 \times n_2) \cup \ldots \cup (2 \times n_m)$ m-matrix of m-rank $2 \cup 2 \cup \ldots \cup 2$ over GF(2). Thus

$$(yM)H^T = y_1 M_1 H_1^T \cup y_2 M_2 H_2^T \cup \ldots \cup y_m M_m H_m^T$$
$$= 0 \cup 0 \cup \ldots \cup 0$$

implies that

$$y(MH^T) = y_1(M_1 H_1^T) \cup y_2(M_2 H_2^T) \cup \ldots \cup y_m(M_m H_m^T)$$
$$= 0 \cup 0 \cup \ldots \cup 0.$$

Since $y = y_1 \cup y_2 \cup \ldots \cup y_m$ is non zero $y(MH^T) = (0 \cup 0 \cup \ldots \cup 0)$ implies $y_i(M_i H_i^T) = 0$ for $i = 1, 2, 3, \ldots, m$; that is the $2 \times 3$ m-matrix $M_1 H_1^T \cup M_2 H_2^T \cup \ldots \cup M_m H_m^T$ has m-rank less than 2 over $GF(2^N)$. Now let

$$P_1 = P_1^1 \cup P_2^1 \cup \ldots \cup P_m^1$$
$$= \{j_1^1 \text{ such that } m_{1j_1^1}^1 = 1\} \cup \{j_2^2 \text{ such that } m_{1j_2^2}^2 = 1\} \cup \ldots$$
$$\cup \{j_m^m \text{ such that } m_{1j_m^m}^m = 1\}$$

and

$$P_2 = P_1^2 \cup P_2^2 \cup \ldots \cup P_m^2$$
$$= \{j_1^1 \text{ such that } m_{2j_1^1}^1 = 1\} \cup \{j_2^2 \text{ such that } m_{2j_2^2}^2 = 1\} \cup \ldots$$
$$\cup \{j_m^m \text{ such that } m_{2j_m^m}^m = 1\}.$$

Since $M = M_1 \cup M_2 \cup \ldots \cup M_m = (m_{ij}^1) \cup (m_{ij}^2) \cup \ldots \cup (m_{ij}^m)$ is a $2 \times n_1 \cup 2 \times n_2 \cup \ldots \cup 2 \times n_m$ m-matrix of m-rank $2 \cup 2 \cup \ldots \cup 2$ and $P_1$ and $P_2$ are disjoint non empty m-subsets of $\{1, 2, \ldots, n_1\} \cup \{1, 2, \ldots, n_2\} \cup \ldots \cup \{1, 2, \ldots, n_m\}$ respectively and

$$MH^T = M_1 H_1^T \cup M_2 H_2^T \cup \ldots \cup M_m H_m^T$$

$$= \begin{pmatrix} \sum_{j_1^1 \in P_1^1} \alpha^1_{1j_1^1} & \sum_{j_1^1 \in P_1^1} \alpha^1_{2j_1^1} & \sum_{j_1^1 \in P_1^1} \alpha^1_{3j_1^1} \\ \\ \sum_{j_1^1 \in P_2^1} \alpha^1_{1j_1^1} & \sum_{j_1^1 \in P_2^1} \alpha^1_{2j_1^1} & \sum_{j_1^1 \in P_2^1} \alpha^1_{3j_1^1} \end{pmatrix} \cup$$

$$\begin{pmatrix} \sum_{j_2^2 \in P_1^2} \alpha^2_{1j_2^2} & \sum_{j_2^2 \in P_1^2} \alpha^2_{2j_2^2} & \sum_{j_2^2 \in P_1^2} \alpha^2_{3j_2^2} \\ \\ \sum_{j_2^2 \in P_2^2} \alpha^2_{1j_2^2} & \sum_{j_2^2 \in P_2^2} \alpha^2_{2j_2^2} & \sum_{j_2^2 \in P_2^2} \alpha^2_{3j_2^2} \end{pmatrix} \cup \ldots \cup$$

$$\begin{pmatrix} \sum_{j_1^1 \in P_1^1} \alpha^m_{1j_1^1} & \sum_{j_1^1 \in P_1^1} \alpha^m_{2j_1^1} & \sum_{j_1^1 \in P_1^1} \alpha^m_{3j_1^1} \\ \\ \sum_{j_2^2 \in P_2^2} \alpha^m_{1j_2^2} & \sum_{j_2^2 \in P_2^2} \alpha^m_{2j_2^2} & \sum_{j_2^2 \in P_2^2} \alpha^m_{3j_2^2} \end{pmatrix}.$$

But the selection of $H = H_1 \cup H_2 \cup \ldots \cup H_m$ is such that there exists $i_1^p, i_2^p \in \{1, 2, 3\}$; $p = 1, 2, \ldots, m$ such that

$$\sum_{j_1^1 \in P_1} \alpha_{i_1^1 j_1^1} \cdot \sum_{k_1 \in P_2} \alpha_{i_2^1 k_1} \cup \sum_{j_2^2 \in P_1} \alpha_{i_1^2 j_1^1} \cdot \sum_{k_2 \in P_2} \alpha_{i_2^2 k_2} \cup \ldots \cup$$
$$\sum_{j_2^2 \in P_1} \alpha_{i_1^m j_1^1} \cdot \sum_{k_m \in P_2} \alpha_{i_2^m k_m}$$
$$\neq \sum_{j_1^1 \in P_1} \alpha_{i_2^1 j_1^1} \cdot \sum_{k_1 \in P_2} \alpha_{i_1^1 k_1} \cup \sum_{j_2^2 \in P_1} \alpha_{i_2^2 j_1^1} \cdot \sum_{k_2 \in P_2} \alpha_{i_1^2 k_2} \cup \ldots \cup$$
$$\sum_{j_1^1 \in P_1} \alpha_{i_2^m j_1^1} \cdot \sum_{k_m \in P_2} \alpha_{i_1^m k_m}.$$

Hence in $MH^T$ there exists a $2 \times 2$ m-submatrices whose determinant is non zero; i.e.,

$$r(MH^T) = r_1(M_1 H_1^T) \cup r_2(M_2 H_2^T) \cup \ldots \cup r_m(M_m H_m^T)$$

over $GF(2^N)$. But this is a contradiction to fact that,

$$\mathrm{rank}(MH^T) = \mathrm{rank}(M_1 H_1^T) \cup \mathrm{rank}(M_2 H_2^T) \cup ... \cup \mathrm{rank}(M_m H_m^T)$$
$$< 2 \cup 2 \cup ... \cup 2.$$

Hence the proof.

Now using constant rank code we proceed on to define the notion of constant rank m-codes of m-length $n_1 \cup n_2 \cup ... \cup n_m$.

**DEFINITION 3.24:** *Let $C_1 \cup C_2 \cup ... \cup C_n$ be a RD-m-code where $C_i$ is a constant rank code of length $n_i$, $i = 1, 2, ..., m$ (Each $C_i$ is a subset of the rank space $V^{n_i}$; $i = 1, 2, ..., m$) then $C = C_1 \cup C_2 \cup ... \cup C_m$ is a constant m-rank code of m-length $n_1 \cup n_2 \cup ... \cup n_m$; that is every m-code word has same m-rank.*

**DEFINITION 3.25:** *$A(n_1, r_1, d_1) \cup A(n_2, r_2, d_2) \cup ... \cup A(n_m, r_m, d_m)$ is defined as the maximum number of m-vectors in $V^{n_1} \cup V^{n_2} \cup ... \cup V^{n_m}$ of constant m-rank, $r_1 \cup r_2 \cup ... \cup r_m$ and m-distance between any two m-vectors is at least $d_1 \cup d_2 \cup ... \cup d_m$ [By $(n_1, r_1, d_1) \cup (n_2, r_2, d_2) \cup ... \cup (n_m, r_m, d_m)$ m-set we mean a m-subset of m-vectors of $V^{n_1} \cup V^{n_2} \cup ... \cup V^{n_m}$ having constant m-rank $r_1 \cup r_2 \cup ... \cup r_m$ and m-distance between any two m-vectors is atleast $d_1 \cup d_2 \cup ... \cup d_m$].*

We analyze the m-function $A(n_1, r_1, d_1) \cup A(n_2, r_2, d_2) \cup ... \cup A(n_m, r_m, d_m)$ by the following theorem:

**THEOREM 3.5:**
1. $A(n_1, r_1, 1) \cup A(n_2, r_2, 1) \cup ... \cup A(n_m, r_m, 1) = L_{r_1}(n_1) \cup L_{r_2}(n_2) \cup ... \cup L_{r_m}(n_m)$, the number of m-vectors of m-rank $r_1 \cup r_2 \cup ... \cup r_m$ in $V^{n_1} \cup V^{n_2} \cup ... \cup V^{n_m}$.
2. $A(n_1, r_1, d_1) \cup A(n_2, r_2, d_2) \cup ... \cup A(n_m, r_m, d_m) = 0 \cup 0 \cup ... \cup 0$ if $r_i > 0$ or $d_i > n_i$ and $d_i > 2r_i$ ($i = 1, 2, 3, ..., m$),

*Proof:*
(1) Follows from the fact that $L_{r_1}(n_1) \cup L_{r_2}(n_2) \cup ... \cup L_{r_m}(n_m)$ is the number of m-vectors of m-length $n_1 \cup n_2 \cup ... \cup n_m$, constant m-rank $r_1 \cup r_2 \cup ... \cup r_m$ and m-distance between any two distinct m-vectors in a m-rank space $V^{n_1} \cup V^{n_2} \cup ... \cup V^{n_m}$ is always greater than or equal to $1 \cup 1 \cup ... \cup 1$.

(2) Follows immediately from the definition of $A(n_1, r_1, d_1) \cup A(n_2, r_2, d_2) \cup ... \cup A(n_m, r_m, d_m)$.

**THEOREM 3.6:** $A(n_1, 1, 2) \cup A(n_2, 1, 2) \cup ... \cup A(n_m, 1, 2) = 2^{n_1} - 1 \cup 2^{n_2} - 1 \cup .... \cup 2^{n_m} - 1$ over any Galois field $GF(2^N)$.

*Proof:* Denote by $V_1 \cup V_2 \cup ... \cup V_m$ the set of m-vectors of m-rank $1 \cup 1 \cup ... \cup 1$ in $V^{n_1} \cup V^{n_2} \cup ... \cup V^{n_m}$; we know each non zero element $\alpha_1 \cup \alpha_2 \cup ... \cup \alpha_m \in GF(2^N)$, there exists $(2^{n_1} - 1) \cup (2^{n_2} - 1) \cup ... \cup (2^{n_m} - 1)$ m-vectors of m-rank one having $\alpha_1 \cup \alpha_2 \cup ... \cup \alpha_m$ as a coordinate. Thus the m-cardinality of $V_1 \cup V_2 \cup ... \cup V_m$ is $(2^N - 1)(2^{n_1} - 1) \cup (2^N - 1)(2^{n_2} - 1) \cup ... \cup (2^N - 1)(2^{n_m} - 1)$. Now m-divide $V_1 \cup V_2 \cup ... \cup V_m$ into $(2^{n_1} - 1) \cup (2^{n_2} - 1) \cup ... \cup (2^{n_m} - 1)$ blocks of $(2^N - 1) \cup (2^N - 1) \cup ... \cup (2^N - 1)$ m-vectors such that each block consists of the same pattern of all nonzero m-elements of $GF(2^N) \cup GF(2^N) \cup ... \cup GF(2^N)$.

Then from each m-block element almost one m-vector can be choosen such that the selected m-vectors are atleast rank 2 apart from each other. Such a m-set we call as $(n_1, 1, 2) \cup (n_2, 1, 2) \cup ... \cup (n_m, 1, 2)$ m-set. Also it is always possible to construct such a m-set. Thus $A(n_1, 1, 2) \cup A(n_2, 1, 2) \cup ... \cup A(n_m, 1, 2) = (2^{n_1} - 1) \cup (2^{n_2} - 1) \cup ... \cup (2^{n_m} - 1)$.

**THEOREM 3.7:** $A(n_1, n_1, n_1) \cup A(n_2, n_2, n_2) \cup ... \cup A(n_m, n_m, n_m) = (2^N - 1) \cup (2^N - 1) \cup ... \cup (2^N - 1)$ (i.e., $A(n_i, n_i, n_i) = 2^N - 1$; $i = 1, 2, 3, ..., m$) over $GF(2^N)$.

*Proof:* Denote by $V_{n_1} \cup V_{n_2} \cup ... \cup V_{n_m}$ the m-set of all m-vectors of m-rank $n_1 \cup n_2 \cup ... \cup n_m$ in the m-space $V^{n_1} \cup V^{n_2} \cup ... \cup V^{n_m}$. We know the m-cardinality of $V_{n_1} \cup V_{n_2} \cup ... \cup V_{n_m}$ is

$$(2^N - 1)(2^N - 2) \ldots (2^N - 2^{n_1-1}) \cup (2^N - 1)(2^N - 2) \ldots$$
$$(2^N - 2^{n_2-1}) \cup \ldots \cup (2^N - 1)(2^N - 2) \ldots (2^N - 2^{n_m-1})$$

and by the definition in a $(n_1, n_1, n_1) \cup (n_2, n_2, n_2) \cup \ldots \cup (n_m, n_m, n_m)$ m-set, the m-distance between any two m-vector should be $n_1 \cup n_2 \cup \ldots \cup n_m$. Thus no two m-vectors can have a common symbol at a co-ordinate place $i_1 \cup i_2 \cup \ldots \cup i_m$; ($1 \leq i_1 \leq n_1$, $1 \leq i_2 \leq n_2$, …, $1 \leq i_m \leq n_m$). This implies that $A(n_1, n_1, n_1) \cup A(n_2, n_2, n_2) \cup \ldots \cup A(n_m, n_m, n_m) \leq (2^N - 1) \cup (2^N - 1) \cup \ldots \cup (2^N - 1)$.

Now we construct a $(n_1, n_1, n_1) \cup (n_2, n_2, n_2) \cup \ldots \cup (n_m, n_m, n_m)$ m-set as follows:

Select N m-vectors from $V_{n_1} \cup V_{n_2} \cup \ldots \cup V_{n_m}$ such that

i. Each m-basis m-elements of $GF(2^{n_1}) \cup GF(2^{n_2}) \cup \ldots \cup GF(2^{n_m})$ should occur (can be as a m-combination) atleast once in each m-vector.

ii. If the $(i_1^{th}, i_2^{th}, \ldots, i_m^{th})$ m-vector is choosen $((i_1 + 1)^{th}, (i_2 + 1)^{th}, \ldots, (i_m + 1)^{th})$ m-vector should be selected such that its m-rank m-distance from any m-linear combination of the previous $(i_1, i_2, \ldots, i_m)$ m-vectors is $n_1 \cup n_2 \cup \ldots \cup n_m$. Now the set of all m-linear combinations of these $N \cup N \cup \ldots \cup N$, m-vectors over $GF(2) \cup GF(2) \cup \ldots \cup GF(2)$ will be such that the m-distance between any two m-vectors is $n_1 \cup n_2 \cup \ldots \cup n_m$. Hence it is $(n_1, n_1, n_1) \cup (n_2, n_2, n_2) \cup \ldots \cup (n_m, n_m, n_m)$ m-set. Also the m-cardinally of this $(n_1, n_1, n_1) \cup (n_2, n_2, n_2) \cup \ldots \cup (n_m, n_m, n_m)$ m-sets is $(2^N - 1) \cup (2^N - 1) \cup \ldots \cup (2^N - 1)$ (we do not count all zero m-linear combinations). Thus $A(n_1, n_1, n_1) \cup A(n_2, n_2, n_2) \cup \ldots \cup A(n_m, n_m, n_m) = (2^N - 1) \cup (2^N - 1) \cup \ldots \cup (2^N - 1)$.

Recall a [n, 1] repetition RD code is a code generated by the matrix $G = (1, 1, \ldots, 1)$ over $F_{2^N}$. Any non zero code word has rank 1.

**DEFINITION 3.26:** *A $[n_1, 1] \cup [n_2, 1] \cup ... \cup [n_m, 1]$ repetition RD m-code is a m-code generated by the m-matrix $G = G_1 \cup G_2 \cup ... \cup G_m = (1\ 1\ ...\ 1) \cup (1\ 1\ ...\ 1) \cup ... \cup (1\ 1\ ...\ 1)$ $(G_i \neq G_j$ if $i \neq j$; $1 \leq i, j \leq m)$ over $F_{2^N}$. Any non zero m-code word has m-rank $1 \cup 1 \cup ... \cup 1$.*

**DEFINITION 3.27:** *Let $C = C_1 \cup C_2 \cup ... \cup C_m$ be a linear $[n_1, k_1] \cup [n_2, k_2] \cup ... \cup [n_m, k_m]$ RD m-code defined over $F_{2^N}$. The covering m-radius of $C = C_1 \cup C_2 \cup ... \cup C_m$ is defined as the smallest m-tuple of integers $(r_1, r_2, ..., r_m)$ such that all m-vectors in the rank m-space $F_{2^N}^{n_1} \cup F_{2^N}^{n_2} \cup ... \cup F_{2^N}^{n_m}$ are with in the rank m-distance $r_1 \cup r_2 \cup ... \cup r_m$ of some m-code word.*

The covering m-radius of $C = C_1 \cup C_2 \cup ... \cup C_m$ is denoted by
$$t(C_1) \cup t(C_2) \cup ... \cup t(C_m) = t(C)$$

$$= \max_{x_1 \in F_{2^N}^{n_1}} \left\{ \begin{matrix} \min(r_1(x_1 + C_1)) \\ c_1 \in C_1 \end{matrix} \right\} \cup \max_{x_2 \in F_{2^N}^{n_2}} \left\{ \begin{matrix} \min(r_2(x_2 + C_2)) \\ c_2 \in C_2 \end{matrix} \right\}$$
$$\cup ... \cup \max_{x_m \in F_{2^N}^{n_m}} \left\{ \begin{matrix} \min(r_m(x_m + C_m)) \\ c_m \in C_m \end{matrix} \right\}.$$

**THEOREM 3.8:** *The linear $[n_1, k_1] \cup [n_2, k_2] \cup ... \cup [n_m, k_m]$ RD-m-code $C = C_1 \cup C_2 \cup ... \cup C_m$ satisfies $t(C) = t(C_1) \cup t(C_2) \cup ... \cup t(C_m) \leq (n_1 - k_1) \cup (n_2 - k_2) \cup ... \cup (n_m - k_m)$.*

*Proof:* Let $C = C_1 \cup C_2 \cup ... \cup C_n$ be a $(n_1, k_1) \cup ... \cup (n_m, k_m)$ RD-m-code. Consider the m-generator m-matrix
$$G = G_1 \cup G_2 \cup ... \cup G_m =$$
$$(I_{k_1}, A_{k_1, n_1 - k_1}) \cup (I_{k_2}, A_{k_2, n_2 - k_2}) \cup ... \cup (I_{k_m}, A_{k_m, n_m - k_m}).$$
Suppose
$$x = x_1 \cup x_2 \cup ... \cup x_m$$
$$= (x_1^1, x_2^1, ..., x_{k_1}^1, x_{k_1+1}^1, ..., x_{n_1}^1) \cup ...$$
$$\cup (x_1^m, x_2^m, ..., x_{k_m}^m, x_{k_m+1}^m, ..., x_{n_m}^m)$$

be any m vector in $V^{n_1} \cup \ldots \cup V^{n_m}$.
Let
$$C = (C_1 \cup C_2 \cup \ldots \cup C_n)(G_1 \cup G_2 \cup \ldots \cup G_m)$$
$$= C_1G_1 \cup C_2G_2 \cup \ldots \cup C_nG_m$$
$$= (x_1^1 \ldots x_{k_1}^1)G_1 \cup \ldots \cup (x_1^m \ldots x_{k_m}^m)G_m.$$

Then C is a m-code word of C and $r(x+c) = r_1(x_1 + c_1) \cup \ldots \cup r_m(x_m + c_m) \leq n_1 - k_1 \cup \ldots \cup n_m - k_m$. Hence the proof.

For any $[n_1, k_1] \cup [n_2, k_1] \cup \ldots \cup [n_m, k_1]$ repetition RD m-code generated by the m-matrix $G = G_1 \cup G_2 \cup \ldots \cup G_m = (1\ 1\ \ldots\ 1) \cup (1\ 1\ \ldots\ 1) \cup \ldots \cup (1\ 1\ \ldots\ 1)$ ($G_i \neq G_j$ if $i \neq j$; $1 \leq i, j \leq m$) over $F_{2^N}$. A non zero m-code word of it has m-rank $1 \cup 1 \cup \ldots \cup 1$.

We proceed on to define the notion of covering m-radius.

**THEOREM 3.9:** *The covering m-radius of a $[n_1, 1] \cup [n_2, 1] \cup \ldots \cup [n_m, 1]$ repetition RD m-code over $F_{2^N}$ is $(n_1 - 1) \cup (n_2 - 1) \cup \ldots \cup (n_m - 1)$.*

*Proof:* The Cartesian m-product of 2 linear RD m-codes
$$C = C_1[n_1^1, k_1^1] \cup C_2[n_2^1, k_2^1] \cup \ldots \cup C_m[n_m^1, k_m^1]$$
and
$$D = D_1[n_1^2, k_1^2] \cup D_2[n_2^2, k_2^2] \cup \ldots \cup D_m[n_m^2, k_m^2]$$
over $F_{2^N}$ is given by
$$C \times D = C_1 \times D_1 \cup C_2 \times D_2 \cup \ldots \cup C_m \times D_m$$
$$= \{(a_1^1, b_1^1) / a_1^1 \in C_1 \text{ and } b_1^1 \in D_1\} \cup$$
$$\{(a_1^2, b_1^2) / a_1^2 \in C_2 \text{ and } b_1^2 \in D_2\} \cup \ldots \cup$$
$$\{(a_1^m, b_1^m) / a_1^m \in C_m \text{ and } b_1^m \in D_m\}.$$

$C \times D$ is a $\{(n_1^1 + n_1^2) \cup (n_2^1 + n_2^2) \cup \ldots \cup (n_m^1 + n_m^2), (k_1^1 + k_1^2) \cup (k_2^1 + k_2^2) \cup \ldots \cup (k_m^1 + k_m^2)\}$ linear RD m-code.
(We assume $n_i^1 + n_i^2 \leq N$ for $i = 1, 2, \ldots, m$).

Now the reader is expected to prove the following theorem:

**THEOREM 3.10:** *If $C = C_1 \cup C_2 \cup \ldots \cup C_m$ and $D = D_1 \cup D_2 \cup \ldots \cup D_m$ be two linear RD m-codes then $t(C \times D) \leq (t(C_1) + t(D_1)) \cup (t(C_2) + t(D_2)) \cup \ldots \cup (t(C_m) + t(D_m))$.*

*Hint:* If $C = C_1 \cup C_2 \cup \ldots \cup C_m$ and $D = D_1 \cup D_2 \cup \ldots \cup D_m$ then $C \times D = \{C_1 \times D_1\} \cup \{C_2 \times D_2\} \cup \ldots \cup \{C_m \times D_m\}$ and $t(C \times D) = t(C_1 \times D_1) \cup t(C_2 \times D_2) \cup \ldots \cup t(C_m \times D_m) \leq \{t(C_1) + t(D_1)\} \cup \{t(C_2) + t(D_2)\} \cup \ldots \cup \{t(C_m) + t(D_m)\}$.

Next we proceed on to define the notion of m-divisible linear RD m-codes we have earlier defined the notion of bidivisible linear RD bicodes

**DEFINITION 3.28:** $C = C_1[n_1, k_1, d_1] \cup C_2[n_2, k_2, d_2] \cup \ldots \cup C_m[n_m, k_m, d_m]$ *be a linear RD m-code over* $F_{q^N}$, $n_i \leq N$, $1 \leq i \leq m$ *and* $N > 1$. *If there exists* $(m_1, m_2, \ldots, m_m)$ $(m_i > 1; i = 1, 2, \ldots, m)$ *such that*

$$\frac{m_i}{r_i(c_i; q)};$$

$1 \leq i \leq n_i$; $i = 1, 2, \ldots, m$ *for all* $c_i \in C_i$; *then we say the m-code C is m-divisible.*

**THEOREM 3.11:** *Let* $C = C_1[n_1, 1, n_1] \cup C_2[n_2, 1, n_2] \cup \ldots \cup C_m[n_m, 1, n_m]$ $(n_i \neq n_j, i \neq j; 1 \leq i, j \leq m)$ *be a MRD m-code for all* $n_i \leq N$, $1 \leq i \leq m$. *Then C is a m-divisible m-code.*

*Proof:* Since there cannot exists m-code words of m-rank greater than $(n_1, n_2, \ldots, n_m)$ in an $[n_1, 1, n_1] \cup [n_2, 1, n_2] \cup \ldots \cup [n_m, 1, n_m]$ MRD m-code. C is a m-divisible m-code.

**DEFINITION 3.29:** *Let* $C_i = [n_i, k_i]$ *be a linear RD-code,* $i = 1, 2, 3, \ldots, m_1$ *and* $C_j = [n_j, k_j, d_j]$ *linear divisible RD codes,* $j = 1, 2, 3, \ldots, m_2$ *defined over* $GF(2^N)$. *Let* $m = m_1 + m_2$, *then the RD linear m-code* $C = C_1 \cup C_2 \cup \ldots \cup C_m$ *is defined as quasi divisible RD m-code,* $n_i \leq N$, $1 \leq i \leq m$.

**DEFINITION 3.30:** *Let $C_i = C_i[n_i, k_i, d_i]$ be a MRD code which is not divisible and $C_j = C_j[n_j, k_j, d_j]$ be a divisible MRD code defined over $GF(2^N)$; $i = 1, 2, 3, …, m_1$ and $j = 1, 2, 3, …, m_2$ such that $m = m_1 + m_2$. $C = C_1 \cup C_2 \cup … \cup C_m$ is defined to be a quasi divisible MRD m-code.*

**DEFINITION 3.31:** *Let $C_1, C_2, …, C_{m_1}$ be circulant rank codes and $C_j[n_j, k_j, d_j]$ a divisible RD-code defined over $GF(2^N)$; $j = 1, 2, …, m_2$ such that $m_1 + m_2 = m$. Then $C = C_1 \cup C_2 \cup … \cup C_m$ is defined to be a quasi divisible circulant rank m-code.*

**DEFINITION 3.32:** *Let $C_1, C_2, …, C_{m_1}$ be AMRD codes and $C_j[n_j, k_j, d_j]$ be a divisible RD code defined over $GF(2^N)$, $j = 1, 2, …, m_2$ such that $m_1 + m_2 = m$. Then $C = C_1 \cup C_2 \cup … \cup C_m$ is defined to be the quasi divisible AMRD m-code.*

We see non divisible MRD m-codes exists as there exists non divisible MRD bicodes.

**DEFINITION 3.33:** *Let $C_i = C_i[n_i, k_i, d_i]$, $i = 1, 2, …, m$ be MRD codes defined over $F_{q^N}$, $n_i \leq N$; $i = 1, 2, …, m$ with $n_i \neq n_j$ if $i \neq j$, $1 \leq i, j \leq m$.*
$A_{s_1}[n_1, d_1] \cup A_{s_2}[n_2, d_2] \cup … \cup A_{s_m}[n_m, d_m]$ *be the number of m-code words with rank m-norms $s_i$ in the linear $[n_i, k_i, d_i]$ MRD-code $1 \leq i \leq m$. Then m-spectrum of the MRD m-code $C = C_1 \cup C_2 \cup … \cup C_m$ is described by the formulae*

$A_0(n_1, d_1) \cup A_0(n_2, d_2) \cup … \cup A_0(n_m, d_m) = 1 \cup 1 \cup … \cup 1$

$$A_{d_1+m_1}(n_1, d_1) \cup A_{d_2+m_2}(n_2, d_2) \cup … \cup A_{d_m+m_m}(n_m, d_m)$$

$$= \begin{bmatrix} n_1 \\ d_1 + m_1 \end{bmatrix} \sum_{j_1=0}^{m_1} (-1)^{j_1+m_1} \begin{bmatrix} d_1 + m_1 \\ d_1 + j_1 \end{bmatrix} q^{\frac{(m_1-j_1)(m_1-j_1-1)(Q^{j_1+1}-1)}{2}} \cup$$

$$\begin{bmatrix} n_2 \\ d_2+m_2 \end{bmatrix} \sum_{j_2=0}^{m_2} (-1)^{j_2+m_2} \begin{bmatrix} d_2+m_2 \\ d_2+j_2 \end{bmatrix} q^{\frac{(m_2-j_2)(m_2-j_2-1)(Q^{j_2+1}-1)}{2}} \cup ... \cup$$

$$\begin{bmatrix} n_m \\ d_m+m_m \end{bmatrix} \sum_{j_m=0}^{m_m} (-1)^{j_m+m_m} \begin{bmatrix} d_m+m_m \\ d_m+j_m \end{bmatrix} q^{\frac{(m_m-j_m)(m_m-j_m-1)(Q^{j_m+1}-1)}{2}}$$

where $Q = q^N$;

$$\begin{bmatrix} n_i \\ m_i \end{bmatrix} = \frac{(q^{n_i}-1)(q^{n_i}-2)...(q^{n_i}-q^{m_i-1})}{(q^{m_i}-1)(q^{m_i}-2)...(q^{m_i}-q^{m_i-1})}; i = 1, 2, ..., m.$$

Using the m-spectrum of a MRD m-code we prove the following theorem:

**THEOREM 3.12:** *All MRD m-codes $C_1[n_1, k_1, d_1] \cup C_2[n_2, k_2, d_2] \cup ... \cup C_m[n_m, k_m, d_m]$ with $d_i < n_i$ (i.e., with $k_i \geq 2$) $1 \leq i \leq m$ are non m-divisible.*

*Proof:* This is proved by making use of the m-spectrum of the MRD m-code. Clearly
$$A_{d_1}(n_1,d_1) \cup A_{d_2}(n_2,d_2) \cup ... \cup A_{d_m}(n_m,d_m) \neq 0 \cup 0 \cup ... \cup 0.$$
If the existence of a m-code word with m-rank $(d_1+1) \cup (d_2+1) \cup ... \cup (d_m+1)$ is established then the proof is complete as the m-gcd $\{(d_1, d_1+1) \cup (d_2, d_2+1) \cup ... \cup (d_m, d_m+1)\} = 1 \cup 1 \cup ... \cup 1$. So the proof is to show that,

$$A_{d_1+1}(n_1,d_1) \cup A_{d_2+1}(n_2,d_2) \cup ... \cup A_{d_m+1}(n_m,d_m)$$

is non zero (i.e., $A_{d_i+1}(n_i,d_i) \neq 0$; $i = 1, 2, ..., m$).

Now
$$A_{d_1+1}(n_1,d_1) \cup A_{d_2+1}(n_2,d_2) \cup ... \cup A_{d_m+1}(n_m,d_m)$$

$$= \begin{bmatrix} n_1 \\ d_1+1 \end{bmatrix} \left( -\begin{bmatrix} d_1+1 \\ d_1 \end{bmatrix} [(Q-1)+(Q^2-1)] \right) \cup$$

$$\begin{bmatrix} n_2 \\ d_2+1 \end{bmatrix}\left(-\begin{bmatrix} d_2+1 \\ d_2 \end{bmatrix}\left[(Q-1)+(Q^2-1)\right]\right)\cup\ldots\cup$$

$$\begin{bmatrix} n_m \\ d_m+1 \end{bmatrix}\left(-\begin{bmatrix} d_m+1 \\ d_m \end{bmatrix}\left[(Q-1)+(Q^2-1)\right]\right)$$

$$=\begin{bmatrix} n_1 \\ d_1+1 \end{bmatrix}(Q-1)\left(Q+1-\begin{bmatrix} d_1+1 \\ d_1 \end{bmatrix}\right)\cup$$

$$\begin{bmatrix} n_2 \\ d_2+1 \end{bmatrix}(Q-1)\left(Q+1-\begin{bmatrix} d_2+1 \\ d_2 \end{bmatrix}\right)\cup\ldots\cup$$

$$\begin{bmatrix} n_m \\ d_m+1 \end{bmatrix}(Q-1)\left(Q+1-\begin{bmatrix} d_m+1 \\ d_m \end{bmatrix}\right).$$

Suppose

$$Q+1-\begin{bmatrix} d_1+1 \\ d_1 \end{bmatrix}\cup Q+1-\begin{bmatrix} d_2+1 \\ d_2 \end{bmatrix}\cup\ldots\cup Q+1-\begin{bmatrix} d_m+1 \\ d_m \end{bmatrix}$$

$$=0\cup 0\cup\ldots\cup 0;$$

i.e.,

$$q^N+1=\frac{q^{d+1}-1}{q-1}$$

i.e.,

$$q-1=\frac{q^d-1}{q^{N-1}}.$$

Clearly,

$$\frac{q^d-1}{q^{N-1}}<1$$

For,

$$\text{if } \frac{q^d-1}{q^{N-1}}\geq 1 \text{ then } q^{N-1}<q^d-1$$

which is impossible as d < n ≤ N. Thus q – 1 < 1 which implies q < 2 a contradiction.

Hence, $A_{d_1+1}(n_1,d_1) \cup A_{d_2+1}(n_2,d_2) \cup ... \cup A_{d_m+1}(n_m,d_m)$ is non zero. Thus expect $C_1(n_1, 1, n_1) \cup C_2(n_2, 1, n_2) \cup ... \cup C_m(n_m, 1, n_m)$ MRD m-codes all $C_1[n_1, k_1, d_1] \cup C_2[n_2, k_2, d_2] \cup ... \cup C_m[n_m, k_m, d_m]$ MRD m-codes with $d_i \le n_i$; i = 1, 2, …, m are non divisible.

Now finally we define the fuzzy rank distance m-codes (m ≥ 3).

Recall Von Kaenel introduced the idea of fuzzy codes with hamming metric and we have defined fuzzy RD codes with Rank metric, we have in the earlier chapter defined the notion of fuzzy RD bicodes. We now proceed onto define the new notion of fuzzy RD m-codes (m ≥ 3) when (m = 3), we call the fuzzy RD m-code to be a fuzzy RD tricode. In the chapter two we have recalled the notion of three types of errors namely asymmetric, symmetric and unidirectional.

We proceed onto define the notion of fuzzy RD m-codes m ≥ 3.

**DEFINITION 3.34:** *Let $V^{n_1} \cup V^{n_2} \cup ... \cup V^{n_m}$ denote the $(n_1, n_2, ..., n_m)$ dimensional vector m-space of $(n_1, n_2, ..., n_m)$-tuples over $F_{2^N}$; $n_i \le N$ and $N > 1$; $1 \le i \le m$.*

*Let $u_i, v_i \in V^{n_i}$; i = 1, 2, ..., m, where,*
$$u_i = (u_1^i, u_2^i, ..., u_{n_i}^i)$$

*and*
$$v_i = (v_1^i, v_2^i, ..., v_{n_i}^i)$$

*with*
$$u_j^i, v_j^i \in F_{2^N};\ 1 \le j \le n_i;\ 1 \le i \le m.$$

*A fuzzy RD m-code word $f_{u_1 \cup u_2 \cup ... \cup u_m} = f_{u_1}^1 \cup f_{u_2}^2 \cup ... \cup f_{u_m}^m$ is a fuzzy m-subset of $V^{n_1} \cup V^{n_2} \cup ... \cup V^{n_m}$ defined by,*

$$f_{u_1 \cup u_2 \cup ... \cup u_m} = f_{u_1}^1 \cup f_{u_2}^2 \cup ... \cup f_{u_m}^m$$
$$= \{(v_1, f_{u_1}^1(v_1)) / v_1 \in V^{n_1}\} \cup \{(v_2, f_{u_2}^2(v_2)) / v_2 \in V^{n_2}\} \cup ... \cup$$

$$\{(v_m, f_{u_m}^m(v_m))/v_m \in V^{n_m}\}$$

where $f_{u_1}^1(v_1) \cup f_{u_2}^2(v_2) \cup \ldots \cup f_{u_m}^m(v_m)$ is the membership m-function.

**DEFINITION 3.35:** *For the symmetric error m-model assume $p_1 \cup p_2 \cup \ldots \cup p_m$ to represent the m-probability that no transition (i.e., error) is made and $q_1 \cup q_2 \cup \ldots \cup q_m$ to represent the m-probability that a m-rank error occurs so that, $p_1 + q_1 \cup p_2 + q_2 \cup \ldots \cup p_m + q_m = 1 \cup 1 \cup \ldots \cup 1$, then*

$$f_{u_1}^1(v_1) \cup f_{u_2}^2(v_2) \cup \ldots \cup f_{u_m}^m(v_m)$$
$$= p_1^{n_1 - r_1} q_1^{r_1} \cup p_2^{n_2 - r_2} q_2^{r_2} \cup \ldots \cup p_m^{n_m - r_m} q_m^{r_m}$$

*where $r_i = r_i(u_i - v_i, 2) = \|u_i - v_i\|$, $i = 1, 2, \ldots, m$.*

**DEFINITION 3.36:** *For unidirectional and asymmetric error m-models assume $q_1 \cup q_2 \cup \ldots \cup q_m$ to represent the probability that $(1 \to 0) \cup (1 \to 0) \cup \ldots \cup (1 \to 0)$ m-transition or $(0 \to 1) \cup (0 \to 1) \cup \ldots \cup (0 \to 1)$ m-transition occurs. Then*

$$f_{u_1}^1(v_1) \cup f_{u_2}^2(v_2) \cup \ldots \cup f_{u_m}^m(v_m)$$
$$= \prod_{i=1}^{n_1} f_{u_i^1}^1(v_i^1) \cup \prod_{i=1}^{n_2} f_{u_i^2}^2(v_i^2) \cup \ldots \cup \prod_{i=1}^{n_m} f_{u_i^m}^m(v_i^m)$$

*where $f_{u_i^1}^1(v_i^1) \cup f_{u_i^2}^2(v_i^2) \cup \ldots \cup f_{u_i^m}^m(v_i^m)$ inherits its definition from the unidirectional and asymmetric m-models respectively, since each $u_i^1 \cup u_i^2 \cup \ldots \cup u_i^m$ or $v_i^1 \cup v_i^2 \cup \ldots \cup v_i^m$ itself is an N-m-tuple over $F_2$. That is since $u_i^j, v_i^j \in F_{2^N}$; $j = 1, 2, \ldots, m$; each $u_i^1 \cup u_i^2 \cup \ldots \cup u_i^m$ or $v_i^1 \cup v_i^2 \cup \ldots \cup v_i^m$ itself is an N, m-tuple from $F_2$.*

$$u_i^j = (u_{i1}^j \cup u_{i2}^j \cup \ldots \cup u_{iN}^j),$$
$$v_k^j = (v_{k1}^j \cup v_{k2}^j \cup \ldots \cup v_{kN}^j)$$

*where, $u_{ip}^j, v_{kl}^j \in F_2$, $1 \leq p, l \leq N$ and $1 \leq i \leq m$.*

*Then for unidirectional error m-model*

$$f^1_{u^1_i}(v^1_i) \cup f^2_{u^2_i}(v^2_i) \cup \ldots \cup f^m_{u^m_i}(v^m_i) =$$

$$\begin{cases} 0 \cup 0 \cup \ldots \cup 0 \\ \quad min(k^1_{i1},k^1_{i2}) \cup \ldots \cup min(k^m_{i1},k^m_{i2}) \neq 0 \cup \ldots \cup 0 \\ p_1^{m^1_i-d^1_i} q_1^{d^1_i} \cup p_2^{m^2_i-d^2_i} q_2^{d^2_i} \cup \ldots \cup p_m^{m^m_i-d^m_i} q_m^{d^m_i} \quad otherwise \end{cases}$$

where $k^j_{it} = \sum_{s=1}^{N} max(0, u^j_{is} - v^j_{is})$ where $j = 1, 2, \ldots, m$ and $i = 1, 2, \ldots, m$.

$$d^j_i = \begin{cases} k^j_{i1} & if \quad k^j_{i2} = 0 \\ k^j_{i2} & if \quad k^j_{i1} = 0 \end{cases}$$

$j = 1, 2, \ldots, m$.

$$m^j_i = \begin{cases} \sum_{s=1}^{N} u^j_{is} & if \quad k^j_{i2} = 0 \\ N - \sum_{s=1}^{N} u^j_{is} & if \quad k^j_{i1} = 0 \\ max\left(\sum u^j_{is}, N - \sum u^j_{is}\right) & if \quad k^j_{i1} = k^d_{i2} = 0 \end{cases}$$

$j = 1, 2, \ldots, m$.

*For the asymmetric error m-model*

$$f^1_{u^1_i}(v^1_i) \cup f^2_{u^2_i}(v^2_i) \cup \ldots \cup f^m_{u^m_i}(v^m_i) =$$

$$\begin{cases} 0 \cup \ldots \cup 0 \\ \quad if \quad min(k^1_{i1},k^1_{i2}) \cup \ldots \cup min(k^m_{i1},k^m_{i2}) \neq 0 \cup \ldots \cup 0 \\ p_1^{m^1_i-d^1_i} q_1^{d^1_i} \cup p_2^{m^2_i-d^2_i} q_2^{d^2_i} \cup \ldots \cup p_m^{m^m_i-d^m_i} q_m^{d^m_i} \quad otherwise \end{cases}$$

where $d^j_i = k^j_{i1}$ and $j = 1, 2, \ldots, m$.

$$m^j_i = \sum_{s=1}^{N} u^j_{is}, j = 1, 2, \ldots, m$$

for asymmetric $(1 \to 0) \cup (1 \to 0)$ error m-model and $d^j_i = k^j_{i2}$, $j = 1, 2, \ldots, m$ and

$$m_i^j = N - \sum_{s=1}^{N} u_{is}^j,$$

$j = 1, 2, \ldots, m$ for the asymmetric $0 \to 1 \cup 0 \to 1$ error m-model.

The results for minimum m-distance of a fuzzy RD-m-code can be derived as in case of minimum bidistance of a fuzzy RD bicode. The notions related to m-covering radius of RD bicodes can be analogously transformed to RD-p-codes ($p \geq 3$).

**PROPOSITION 3.1:** *If*
$$C_1^1 \cup C_2^1 \cup \ldots \cup C_m^1 \text{ and } C_1^2 \cup C_2^2 \cup \ldots \cup C_m^2$$
*are RD m-codes with*
$$C_1^1 \cup C_2^1 \cup \ldots \cup C_m^1 \subseteq C_1^2 \cup C_2^2 \cup \ldots \cup C_m^2$$
$$(i.e., \ C_j^1 \subseteq C_j^2;\ j = 1, 2, \ldots, m)$$

*then*
$$t_{m_1}(C_1^1) \cup t_{m_2}(C_2^1) \cup \ldots \cup t_{m_m}(C_m^1) \geq$$
$$t_{m_1}(C_1^2) \cup t_{m_2}(C_2^2) \cup \ldots \cup t_{m_m}(C_m^2)$$
*(i.e., $t_{m_i}(C_i^1) \geq t_{m_i}(C_i^2)$; $i = 1, 2, \ldots, m$).*

*Proof:* Let $S_i \subseteq V^{n_i}$ with $1 \leq j \leq m$; $|S_i| = m_i$; $i = 1, 2, \ldots, m$

$$\text{cov}(C_1^2, S_1) \cup \text{cov}(C_2^2, S_2) \cup \ldots \cup \text{cov}(C_m^2, S_m)$$
$$= \min\{\text{cov}(x_1, S_1); x_1 \in C_1^2\} \cup$$
$$\min\{\text{cov}(x_2, S_2); x_2 \in C_2^2\} \cup \ldots \cup \min\{\text{cov}(x_m, S_m); x_m \in C_m^2\}$$
$$\leq \min\{\text{cov}(x_1, S_1)/x_1 \in C_1^1\} \cup$$
$$\min\{\text{cov}(x_2, S_2)/x_2 \in C_2^1\} \cup \ldots \cup \min\{\text{cov}(x_m, S_m)/x_m \in C_m^1\}$$
$$= \text{cov}(C_1^1, S_1) \cup \text{cov}(C_2^1, S_2) \cup \ldots \cup \text{cov}(C_m^1, S_m).$$

Thus
$$t_{m_1}(C_1^2) \cup t_{m_2}(C_2^2) \cup \ldots \cup t_{m_m}(C_m^2)$$
$$\leq t_{m_1}(C_1^1) \cup t_{m_2}(C_2^1) \cup \ldots \cup t_{m_m}(C_m^1).$$

**PROPOSITION 3.2:** *For any RD m-code $C = C_1 \cup C_2 \cup \ldots \cup C_m$ and a m-tuple of positive integers $(m_1, m_2, \ldots, m_m)$*

$$t_{m_1}(C_1) \cup t_{m_2}(C_2) \cup \ldots \cup t_{m_m}(C_m) \leq$$
$$t_{m_1+1}(C_1) \cup t_{m_2+1}(C_2) \cup \ldots \cup t_{m_m+1}(C_m).$$

*Proof:* Since $S_i \subseteq V^{n_i}$; $i = 1, 2, \ldots, m$; $S_1 \cup S_2 \cup \ldots \cup S_m$ is a m-subset of $V^{n_1} \cup V^{n_2} \cup \ldots \cup V^{n_m}$.

Now

$$t_{m_1}(C_1) \cup t_{m_2}(C_2) \cup \ldots \cup t_{m_m}(C_m)$$

$$= \max\{cov(C_1, S_1) / S_1 \subseteq V^{n_1}, |S_1| = m_1\} \cup$$
$$\max\{cov(C_2, S_2) / S_2 \subseteq V^{n_2}, |S_2| = m_2\} \cup \ldots \cup$$
$$\max\{cov(C_m, S_m) / S_m \subseteq V^{n_m}, |S_m| = m_m\}$$

$$\leq \max\{cov(C_1, S_1) / S_1 \subseteq V^{n_1}, |S_1| = m_1 + 1\} \cup$$
$$\max\{cov(C_2, S_2) / S_2 \subseteq V^{n_2}, |S_2| = m_2 + 1\} \cup \ldots \cup$$
$$\max\{cov(C_m, S_m) / S_m \subseteq V^{n_m}, |S_m| = m_m + 1\}$$

$$= t_{m_1+1}(C_1) \cup t_{m_2+1}(C_2) \cup \ldots \cup t_{m_m+1}(C_m).$$

**PROPOSITION 3.3:** *For any m-set of positive integers*
$\{n_1, m_1, k_1, K_1\} \cup \{n_2, m_2, k_2, K_2\} \cup \ldots \cup \{n_m, m_m, k_m, K_m\}$;

$$t_{m_1}[n_1, k_1] \cup t_{m_2}[n_2, k_2] \cup \ldots \cup t_{m_m}[n_m, k_m]$$
$$\leq t_{m_1+1}[n_1, k_1] \cup t_{m_2+1}[n_2, k_2] \cup \ldots \cup t_{m_m+1}[n_m, k_m].$$

*Proof:* Given $C_1[n_1, k_1] \cup C_2[n_2, k_2] \cup \ldots \cup C_m[n_m, k_m]$ to be a RD m-code with $C_i \subseteq V^{n_i}$; $i = 1, 2, \ldots, m$.

Now

$$t_{m_1}[n_1, k_1] \cup t_{m_2}[n_2, k_2] \cup \ldots \cup t_{m_m}[n_m, k_m]$$
$$= \min\{t_{m_1}(C_1) / C_1 \subseteq V^{n_1}, \dim C_1 = k_1\} \cup$$
$$\min\{t_{m_2}(C_2) / C_2 \subseteq V^{n_2}, \dim C_2 = k_2\} \cup \ldots \cup$$

$$\min\{t_{m_m}(C_m)/C_m \subseteq V^{n_m}, \dim C_m = k_m\}$$

$$\leq \min\{t_{m_1+1}(C_1)/C_1 \subseteq V^{n_1}, \dim C_1 = k_1\} \cup$$
$$\min\{t_{m_2+1}(C_2)/C_2 \subseteq V^{n_2}, \dim C_2 = k_2\} \cup \ldots \cup$$
$$\min\{t_{m_m+1}(C_m)/C_m \subseteq V^{n_m}, \dim C_m = k_m\}$$

$$= t_{m_1+1}[n_1, k_1] \cup t_{m_2+1}[n_2, k_2] \cup \ldots \cup t_{m_m+1}[n_m, k_m] \ .$$

Similarly we have

$$t_{m_1}[n_1, K_1] \cup t_{m_2}[n_2, K_2] \cup \ldots \cup t_{m_m}[n_m, K_m]$$
$$\leq t_{m_1+1}[n_1, K_1] \cup t_{m_2+1}[n_2, K_2] \cup \ldots \cup t_{m_m+1}[n_m, K_m] \ .$$

That is

$$t_{m_1}[n_1, K_1] \cup t_{m_2}[n_2, K_2] \cup \ldots \cup t_{m_m}[n_m, K_m]$$

$$= \min\{t_{m_1}(C_1)/C_1 \subseteq V^{n_1}, |C_1| = K_1\} \cup$$
$$\min\{t_{m_2}(C_2)/C_2 \subseteq V^{n_2}, |C_2| = K_2\} \cup \ldots \cup$$
$$\min\{t_{m_m}(C_m)/C_m \subseteq V^{n_m}, |C_m| = K_m\}$$

$$\leq \min\{t_{m_1+1}(C_1)/C_1 \subseteq V^{n_1}, |C_1| = K_1\} \cup$$
$$\min\{t_{m_2+1}(C_2)/C_2 \subseteq V^{n_2}, |C_2| = K_2\} \cup \ldots \cup$$
$$\min\{t_{m_m+1}(C_m)/C_m \subseteq V^{n_m}, |C_m| = K_m\}$$

$$\leq t_{m_1+1}[n_1, K_1] \cup t_{m_2+1}[n_2, K_2] \cup \ldots \cup t_{m_m+1}[n_m, K_m] \ .$$

**PROPOSITION 3.4:** *For any m-set of positive integers*
*$\{n_1, m_1, k_1, K_1\} \cup \{n_2, m_2, k_2, K_2\} \cup \ldots \cup \{n_m, m_m, k_m, K_m\}$;*
$$t_{m_1}[n_1, k_1] \cup t_{m_2}[n_2, k_2] \cup \ldots \cup t_{m_m}[n_m, k_m]$$
$$\geq t_{m_1}[n_1, k_1+1] \cup t_{m_2}[n_2, k_2+1] \cup \ldots \cup t_{m_m}[n_m, k_m+1] \ .$$

*Proof:* Given $C = C_1 \cup C_2 \cup \ldots \cup C_m$ is a RD m-code, hence a m-subspace of $V^{n_1} \cup V^{n_2} \cup \ldots \cup V^{n_m}$.
Consider

$$t_{m_1}[n_1, k_1+1] \cup t_{m_2}[n_2, k_2+1] \cup \ldots \cup t_{m_m}[n_m, k_m+1]$$

$$= \min\{t_{m_1}(C_1)/C_1 \subseteq V^{n_1}, \dim C_1 = k_1+1\} \cup$$
$$\min\{t_{m_2}(C_2)/C_2 \subseteq V^{n_2}, \dim C_2 = k_2+1\} \cup \ldots \cup$$
$$\min\{t_{m_m}(C_m)/C_m \subseteq V^{n_m}, \dim C_m = k_m+1\}$$

$$\leq \min\{t_{m_1}(C_1)/C_1 \subseteq V^{n_1}, \dim C_1 = k_1\} \cup$$
$$\min\{t_{m_2}(C_2)/C_2 \subseteq V^{n_2}, \dim C_2 = k_2\} \cup \ldots \cup$$
$$\min\{t_{m_m}(C_m)/C_m \subseteq V^{n_m}, \dim C_m = k_m\}$$

since for each $C_1 \cup C_2 \cup \ldots \cup C_m \subseteq C_{12} \cup C_{22} \cup \ldots \cup C_{m2}$

$$t_{m_1}(C_{12}) \cup t_{m_2}(C_{22}) \cup \ldots \cup t_{m_m}(C_{m2})$$
$$\leq t_{m_1}(C_1) \cup t_{m_2}(C_2) \cup \ldots \cup t_{m_m}(C_m)$$
$$= t_{m_1}[n_1, k_1] \cup t_{m_2}[n_2, k_2] \cup \ldots \cup t_{m_m}[n_m, k_m]..$$

Similarly
$$t_{m_1}(n_1, K_1+1) \cup t_{m_2}(n_2, K_2+1) \cup \ldots \cup t_{m_m}(n_m, K_m+1)$$
$$\leq t_{m_1}(n_1, K_1) \cup t_{m_2}(n_2, K_2) \cup \ldots \cup t_{m_m}(n_m, K_m).$$

Using these results and the fact $k_{im_i}[n_i, t_i]$ denotes the smallest dimension of a linear RD code of length $n_i$ and $m_i$ covering radius $t_i$ and $K_{im_i}[n_i, t_i]$ denotes the least cardinality of the RD codes of length $n_i$ and $m_i$-covering radius $t_i$ the following results can be easily proved.

*Result 1:* For any m-set of positive integers
$$\{n_1, m_1, t_1\} \cup \{n_2, m_2, t_2\} \cup \ldots \cup \{n_m, m_m, t_m\};$$
and

$$k_{m_1}[n_1,t_1] \cup k_{m_2}[n_2,t_2] \cup ... \cup k_{m_m}[n_m,t_m]$$
$$\leq k_{m_1+1}[n_1,t_1] \cup k_{m_2+1}[n_2,t_2] \cup ... \cup k_{m_m+1}[n_m,t_m]$$

and

$$K_{m_1}(n_1,t_1) \cup K_{m_2}(n_2,t_2) \cup ... \cup K_{m_m}(n_m,t_m)$$
$$\leq K_{m_1+1}(n_1,t_1) \cup K_{m_2+1}(n_2,t_2) \cup ... \cup K_{m_m+1}(n_m,t_m).$$

*Result 2:* For any m-set of positive integers
$$\{n_1, m_1, t_1\} \cup \{n_2, m_2, t_2\} \cup ... \cup \{n_m, m_m, t_m\}.$$
we have,
$$k_{m_1}[n_1,t_1] \cup k_{m_2}[n_2,t_2] \cup ... \cup k_{m_m}[n_m,t_m]$$
$$\geq k_{m_1}[n_1,t_1+1] \cup k_{m_2}[n_2,t_2+1] \cup ... \cup k_{m_m}[n_m,t_m+1]$$

and

$$K_{m_1}(n_1,t_1) \cup K_{m_2}(n_2,t_2) \cup ... \cup K_{m_m}(n_m,t_m)$$
$$\geq K_{m_1}(n_1,t_1+1) \cup K_{m_2}(n_2,t_2+1) \cup ... \cup K_{m_m}(n_m,t_m+1).$$

We say a m-function $f_1 \cup f_2 \cup ... \cup f_m$ is a non decreasing m-function in some m-variable say $x_1 \cup x_2 \cup ... \cup x_m$ if each $f_i$ happen to be a non-decreasing function in the variable $x_i$; $i = 1, 2, ..., m$.

With this understanding we have for $(m_1, m_2, ..., m_m)$-covering m-radius of a fixed RD m-code $C_1 \cup C_2 \cup ... \cup C_m$,

$$t_{m_1}[n_1,k_1] \cup t_{m_2}[n_2,k_2] \cup ... \cup t_{m_m}[n_m,k_m],$$
$$k_{m_1}[n_1,t_1] \cup k_{m_2}[n_2,t_2] \cup ... \cup k_{m_m}[n_m,t_m],$$
$$t_{m_1}(n_1,K_1) \cup t_{m_2}(n_2,K_2) \cup ... \cup t_{m_m}(n_m,K_m)$$

and

$$K_{m_1}(n_1,t_1) \cup K_{m_2}(n_2,t_2) \cup ... \cup K_{m_m}(n_m,t_m)$$

are non decreasing m-functions of $(m_1, m_2, ..., m_m)$.
The relationships between the multi covering m-radii of two RD m-codes that are built using them are described.
Let
$$C^i = C_1^i \cup C_2^i \cup .... \cup C_m^i$$
for i = 1, 2 be a

$[n_1^1, k_1^1, d_1^1] \cup [n_1^2, k_1^2, d_1^2] \cup \ldots \cup [n_1^m, k_1^m, d_1^m]$

and

$$[n_2^1, k_2^1, d_2^1] \cup [n_2^2, k_2^2, d_2^2] \cup \ldots \cup [n_2^m, k_2^m, d_2^m]$$

RD m-codes over $F_{2^N}$ with $n_1^i, n_2^i, n_1^i + n_2^i \leq N$ for i = 1, 2, …, m.

**PROPOSITION 3.5:** *Let*
$$C^1 = C_1^1 \cup C_2^1 \cup \ldots \cup C_m^1 \text{ and } C^2 = C_1^2 \cup C_2^2 \cup \ldots \cup C_m^2$$
*be RD m-codes described above*
$$C = C^1 \times C^2 = (C_1^1 \times C_1^2) \cup (C_2^1 \times C_2^2) \cup \ldots \cup (C_m^1 \times C_m^2)$$
$$= \{(x_1/y_1) / x_1 \in C_1^1, y_1 \in C_1^2\} \cup$$
$$\{(x_2/y_2) / x_2 \in C_2^1, y_2 \in C_2^2\} \cup \ldots \cup$$
$$\{(x_m/y_m) / x_m \in C_m^1, y_m \in C_m^2\}.$$
*Then $C^1 \times C^2$ is a*
$$[\, n_1^1 + n_2^1 \cup n_1^2 + n_2^2 \cup \ldots \cup n_1^m + n_2^m,$$
$$k_1^1 + k_2^1 \cup k_1^2 + k_2^2 \cup \ldots \cup k_1^m + k_2^m,$$
$$min\{d_1^1, d_2^1\} \cup min\{d_1^2, d_2^2\} \cup \ldots \cup min\{d_1^m, d_2^m\}\,]$$
*rank distance m-code over $F_{2^N}$ and*
$$t_{m_1}(C_1^1 \times C_1^2) \cup t_{m_2}(C_2^1 \times C_2^2) \cup \ldots \cup t_{m_m}(C_m^1 \times C_m^2)$$
$$\leq t_{m_1}(C_1^1) + t_{m_1}(C_1^2) \cup t_{m_2}(C_2^1)$$
$$+ t_{m_2}(C_2^2) \cup \ldots \cup t_{m_m}(C_m^1) + t_{m_m}(C_m^2).$$

*Proof:* Let $S_i \subseteq V^{n_1^i + n_2^i}$; for i = 1, 2, …, m. and $S_i = \{s_1^i, \ldots, s_{m_i}^i\}$ for i = 1, 2, …, m with $s_i^j(x_{ji}/y_{ji})$ for i = 1, 2, …, m; $x_{ji} \in V^{n_1^j}$ and $y_{ji} \in V^{n_2^j}$, $1 \leq i \leq m_i$; $1 \leq i \leq m$. Let

$$S_1^1 = \{x_{11}, \ldots, x_{1m_1}\}, \quad S_1^2 = \{y_{11}, \ldots, y_{1m_1}\},$$
$$S_2^1 = \{x_{21}, \ldots, x_{2m_2}\}, \quad S_2^2 = \{y_{21}, \ldots, y_{2m_2}\}, \ldots,$$
$$S_m^1 = \{x_{m1}, \ldots, x_{mm_m}\} \text{ and } S_m^2 = \{y_{m1}, \ldots, y_{mm_m}\}.$$

Now $t_{m_i}(C_1^i)$ being the $m_i$ covering radius of $C_1^i$; $i = 1, 2, \ldots, m$, there exists $c_1^i \in C_1^i$ such that $S_1^i \subseteq B_{t_{m_i}}^i(C_1^i)$; $i = 1, 2, \ldots, m$.

This implies as in case of RD bicodes
$$r_j(s_{ji} + C^j) = r_j((x_{ji}/y_{ji}) + (C_1^j/C_2^j)) = r_j(x_{ji} + C_1^j/y_{ji} + C_2^j)$$
$$\leq r_j(x_{ji} + C_1^j) + r_j(y_{ji} + C_2^j) \leq t_{m_j}(C_1^j) + t_{m_j}(C_2^j);$$

$j = 1, 2, \ldots, m$.
Thus
$$t_m(C) = t_{m_1}(C_1^1 \times C_1^2) \cup t_{m_2}(C_2^1 \times C_2^2) \cup \ldots \cup t_{m_m}(C_m^1 \times C_m^2) \leq$$
$$t_{m_1}(C_1^1) + t_{m_1}(C_1^2) \cup t_{m_2}(C_2^1) + t_{m_2}(C_2^2) \cup \ldots \cup t_{m_m}(C_m^1) + t_{m_m}(C_m^2).$$

For any $\underbrace{(r,r,\ldots,r)}_{m-\text{times}}$ (r a positive integer (m ≥ 3) the (r, r, …, r) fold repetition RD m-code $C_1 \cup C_2 \cup \ldots \cup C_m$ is the m-code

$$C = \{(c_1 | c_1 | \ldots | c_1)/c_1 \in C_1\} \cup$$
$$\{(c_2 | c_2 | \ldots | c_2)/c_2 \in C_2\} \cup \ldots \cup$$
$$\{(c_m | c_m | \ldots | c_m)/c_m \in C_m\}$$

where the m-code word $C_1 \cup C_2 \cup \ldots \cup C_m$ is a concatenation of (r, r, …, r) times, this is a $[rn_1, k_1, d_1] \cup [rn_2, k_2, d_2] \cup \ldots \cup [rn_m, k_m, d_m]$ rank distance m-code with $n_i \leq N$ and $rn_i \leq N$; $i = 1, 2, \ldots, m$. Thus any m-code word in $C_1 \cup C_2 \cup \ldots \cup C_m$ would be of the form
$$\{(c_1 | c_1 | \ldots | c_1)\} \cup \{(c_2 | c_2 | \ldots | c_2)\} \cup \ldots \cup \{(c_m | c_m | \ldots | c_m)\}$$
such that $x_i \in C_i$ for $i = 1, 2, \ldots, m$.

We can also define $(r_1, r_2, \ldots, r_m)$ fold repetition m-code ($r_i \neq r_j$ if $i \neq j$; $1 \leq i, j \leq m$).

**DEFINITION 3.37:** *Let $C_1 \cup C_2 \cup \ldots \cup C_m$ be a $[n_1, k_1, d_1] \cup [n_2, k_2, d_2] \cup \ldots \cup [n_m, k_m, d_m]$ RD m-code.*

*Let* $c^i = \left\{ \underbrace{(c_i / c_i / ... / c_i)}_{r_i-\text{times}} \middle/ c_i \in C_i \right\}$ *be a $r_i$-fold repetition RD-code* $C_i$, $i = 1, 2, ..., m$. *Then* $C^1 \cup C^2 \cup ... \cup C^m$ *is defined as the* $(r_1, r_2, ..., r_m)$-*fold repetition m-code each* $r_i n_i < N$ *for* $i = 1, 2, ..., m$.

We prove the following interesting theorem.

**THEOREM 3.13:** *For an $(r, r, ..., r)$ fold repetition m-code* $C_1 \cup C_2 \cup ... \cup C_m$

$$t_{m_1}(C_1) \cup t_{m_2}(C_2) \cup ... \cup t_{m_m}(C_m)$$
$$= t_{m_1}(C^1) \cup t_{m_2}(C^2) \cup ... \cup t_{m_m}(C^m).$$

*Proof:* Let $S_i = \{v_{i1}, v_{i2}, ..., v_{im_i}\} \subseteq V^{n_i}$ for $i = 1, 2, ..., m$; such that $\text{cov}(C_i, S_i) = t_{m_i}(C_i)$; $i = 1, 2, ..., m$.
Let

$$v'_{i1} = (v_{i1} | v_{i1} | ... | v_{i1})$$
$$v'_{i2} = (v_{i2} | v_{i2} | ... | v_{i2})$$

and so on

$$v'_{im_i} = (v_{im_i} | v_{im_i} | ... | v_{im_i}); \ 1 \leq i \leq m.$$

Let $S_i = \{v'_{i1}, v'_{i2}, ..., v'_{im_i}\}$ be the set of $m_i$-vectors of length $m_i$ for, $i = 1, 2, ..., m$. An r fold repetition of any RD code word retains the same rank weight. Hence $(C^i, S'_i) = t_{m_i}(C_i)$ true for i = 1, 2, ..., m. Since

$$t_{m_i}(C^i) = \text{cov}(C^i, S'_i)$$

it follows that

$$t_{m_i}(C^i) \geq t_{m_i}(C_i)$$

for $i = 1, 2, ..., m$; i.e.,

$$t_{m_1}(C^1) \cup t_{m_2}(C^2) \cup ... \cup t_{m_m}(C^m) \geq$$
$$t_{m_1}(C_1) \cup t_{m_2}(C_2) \cup ... \cup t_{m_m}(C_m) \qquad ----- \qquad I$$

Conversely let $S_i = \{v_{i1}, v_{i2}, ..., v_{im_i}\}$ be a set of $m_i$ vectors $i = 1, 2, ..., m$ of length $rn_i$ with $v_{ij} = (v'_{ij} | ... | v'_{ij})$; $j = 1, 2, ..., m_i$, $i = 1, 2, ..., m$. $v'_i \in V^{n_i}$, $1 \leq i \leq m_i$. Then there exists $c_i \in C_i$ such that $d_{R_i}(c_i, v'_i) \leq t_{m_i}(C_i)$; $i = 1, 2, ..., m_i$; $1 \leq i \leq m$.

This implies $d_{R_i}((c_i | c_i | ... | c_i), v_{ij}) \leq t_{m_i}(C_i)$ for every $i$ ($1 \leq i \leq m$). Thus $t_{m_i}(C^i) \leq t_{m_i}(C_i)$, $i = 1, 2, ..., m$;
i.e.,
$$t_{m_1}(C^1) \cup t_{m_2}(C^2) \cup ... \cup t_{m_m}(C^m) \leq$$
$$t_{m_1}(C_1) \cup t_{m_2}(C_2) \cup ... \cup t_{m_m}(C_m) \qquad \text{------- II}$$

From I and II
$$t_{m_1}(C^1) \cup t_{m_2}(C^2) \cup ... \cup t_{m_m}(C^m) =$$
$$t_{m_1}(C_1) \cup t_{m_2}(C_2) \cup ... \cup t_{m_m}(C_m).$$

Now we proceed on to analyse the notion of multi covering m-bounds for RD m-codes. The $(m_1, m_2, ..., m_m)$ covering m-radius $t_{m_1}(C_1) \cup t_{m_2}(C_2) \cup ... \cup t_{m_m}(C_m)$ of a RD m-code $C = C_1 \cup C_2 \cup ... \cup C_m$ is a non-decreasing m-function of $m_1 \cup m_2 \cup ... \cup m_m$. Thus a lower m-bound for
$$t_{m_1}(C_1) \cup t_{m_2}(C_2) \cup ... \cup t_{m_m}(C_m)$$
implies a m-bound for
$$t_{m_1+1}(C_1) \cup t_{m_2+1}(C_2) \cup ... \cup t_{m_m+1}(C_m).$$

First m-bound exhibits that for $m_1 \cup m_2 \cup ... \cup m_m \geq 2 \cup 2 \cup ... \cup 2$ the situation of $(m_1, m_2, ..., m_m)$ covering m-radii is quite different for ordinary covering radii.

**PREPOSITION 3.6:** *If $m_1 \cup m_2 \cup ... \cup m_m > 2 \cup 2 \cup ... \cup 2$ then the $(m_1, m_2, ..., m_m)$ covering m-radii of a RD m-code $C = C_1 \cup C_2 \cup ... \cup C_m$ of m-length $(n_1, n_2, ..., n_m)$ is at least*
$$\left\lceil \frac{n_1}{2} \right\rceil \cup \left\lceil \frac{n_2}{2} \right\rceil \cup ... \cup \left\lceil \frac{n_m}{2} \right\rceil.$$

*Proof:* Let $C = C_1 \cup C_2 \cup \ldots \cup C_m$ be a RD m-code of m-length $(n_1, n_2, \ldots, n_m)$ over $GF(2^N)$. Let $m_1 \cup m_2 \cup \ldots \cup m_m = 2 \cup 2 \cup \ldots \cup 2$. Let $t_1, t_2, \ldots, t_m$ be the 2-covering m-radii of the RD m-code $C = C_1 \cup C_2 \cup \ldots \cup C_m$. Let $x = x_1 \cup x_2 \cup \ldots \cup x_m \in V^{n_1} \cup V^{n_2} \cup \ldots \cup V^{n_m}$. Choose $y = y_1 \cup y_2 \cup \ldots \cup y_m \in V^{n_1} \cup V^{n_2} \cup \ldots \cup V^{n_m}$ such that all the $(n_1, n_2, \ldots, n_m)$ coordinate of

$$x - y = (x_1 - y_1) \cup (x_2 - y_2) \cup \ldots \cup (x_m - y_m)$$

are linearly independent that is

$$d_R(x, y) = d_R(x_1 \cup x_2 \cup \ldots \cup x_m, y_1 \cup y_2 \cup \ldots \cup y_m)$$
$$= d_{R_1}(x_1, y_1) \cup d_{R_2}(x_2, y_2) \cup \ldots \cup d_{R_m}(x_m, y_m)$$

($R = R_1 \cup R_2 \cup \ldots \cup R_m$ and $d_R = d_{R_1} \cup d_{R_2} \cup \ldots \cup d_{R_m}) = n_1 \cup n_2 \cup \ldots \cup n_m$. Then for any $c = c_1 \cup c_2 \cup \ldots \cup c_m \in C_1 \cup C_2 \cup \ldots \cup C_m$.

$$d_R(x+c) + d_R(c+y)$$
$$= d_{R_1}(x_1, c_1) + d_{R_1}(c_1, y_1) \cup d_{R_2}(x_2, c_2) + d_{R_2}(c_2, y_2) \cup \ldots \cup$$
$$d_{R_m}(x_m, c_m) + d_{R_m}(c_m, y_m)$$
$$\geq d_{R_1}(x_1, y_1) \cup d_{R_2}(x_2, y_2) \cup \ldots \cup d_{R_m}(x_m, y_m)$$
$$= n_1 \cup n_2 \cup \ldots \cup n_m;$$

this implies that one of

$$d_{R_1}(x_1, c_1) \cup d_{R_2}(x_2, c_2) \cup \ldots \cup d_{R_m}(x_m, c_m)$$

and

$$d_{R_1}(c_1, y_1) \cup d_{R_2}(c_2, y_2) \cup \ldots \cup d_{R_m}(c_m, y_m)$$

is at least

$$\frac{n_1}{2} \cup \frac{n_2}{2} \cup \ldots \cup \frac{n_m}{2}.$$

(That is one of $d_{R_i}(x_i, c_i)$ and $d_{R_i}(c_i, y_i)$ is atleast $\frac{n_i}{2}$; $i = 1, 2, \ldots, m$) and hence

$$t = t_1 \cup t_2 \cup \ldots \cup t_m \geq \left\lceil \frac{n_1}{2} \right\rceil \cup \left\lceil \frac{n_2}{2} \right\rceil \cup \ldots \cup \left\lceil \frac{n_m}{2} \right\rceil.$$

Since t is a non-decreasing m-function of $m_1 \cup m_2 \cup \ldots \cup m_m$ it follows that

$$t_m(C) = t_{m_1}(C_1) \cup t_{m_2}(C_2) \cup \ldots \cup t_{m_m}(C_m)$$
$$\geq \left\lceil \frac{n_1}{2} \right\rceil \cup \left\lceil \frac{n_2}{2} \right\rceil \cup \ldots \cup \left\lceil \frac{n_m}{2} \right\rceil$$

for $m_1 \cup m_2 \cup \ldots \cup m_m \geq 2 \cup 2 \cup \ldots \cup 2$. m-bounds of the multi covering m-radius of $V^{n_1} \cup V^{n_2} \cup \ldots \cup V^{n_m}$ can be used to obtain m-bounds on the multi covering m-radii of arbitrary m-codes. Thus a relationship between $(m_1, m_2, \ldots, m_m)$ covering m-radii of an RD m-code and that of its ambient m-space $V^{n_1} \cup V^{n_2} \cup \ldots \cup V^{n_m}$ is established.

**THEOREM 3.14:** *Let $C = C_1 \cup C_2 \cup \ldots \cup C_m$ be RD m-code of m-length $n_1 \cup n_2 \cup \ldots \cup n_m$ over $F_{2^N} \cup F_{2^N} \cup \ldots \cup F_{2^N}$. Then for any positive m-integer tuple $(m_1, m_2, \ldots, m_m)$*

$$t_{m_1}^1(C_1) \cup t_{m_2}^2(C_2) \cup \ldots \cup t_{m_m}^m(C_m)$$
$$\leq t_1^1(C_1) + t_{m_1}^1(V^{n_1}) \cup t_1^2(C_2) + t_{m_2}^2(V^{n_2}) \cup \ldots \cup$$
$$t_1^m(C_m) + t_{m_m}^m(V^{n_m}).$$

*Proof:* Let $S = S_1 \cup S_2 \cup \ldots \cup S_m \subseteq V^{n_1} \cup V^{n_2} \cup \ldots \cup V^{n_m}$ (i.e., $S_i \subseteq V^{n_i}$; $i = 1, 2, \ldots, m$) with $|S| = |S_1| \cup |S_2| \cup \ldots \cup |S_m| = m_1 \cup m_2 \cup \ldots \cup m_m$. Then there exists $u = u_1 \cup u_2 \cup \ldots \cup u_m \in V^{n_1} \cup V^{n_2} \cup \ldots \cup V^{n_m}$ such that

$$\text{cov}(u, S) = \text{cov}(u_1, S_1) \cup \text{cov}(u_2, S_2) \cup \ldots \cup \text{cov}(u_m, S_m)$$
$$\leq t_{m_1}^1(V^{n_1}) \cup t_{m_2}^2(V^{n_2}) \cup \ldots \cup t_{m_m}^m(V^{n_m}).$$

Also there is a $c = c_1 \cup c_2 \cup \ldots \cup c_m \in C_1 \cup C_2 \cup \ldots \cup C_m$ such that
$$d_R(c, u) = d_{R_1}(c_1, u_1) \cup d_{R_2}(c_2, u_2) \cup \ldots \cup d_{R_m}(c_m, u_m)$$
$$\leq t_1^1(C_1) \cup t_1^2(C_2) \cup \ldots \cup t_1^m(C_m).$$

Now,

$$\mathrm{cov}(c,S) = \mathrm{cov}(c_1,S_1) \cup \mathrm{cov}(c_2,S_2) \cup ... \cup \mathrm{cov}(c_m,S_m)$$

$$= \max\{d_{R_1}(c_1,y_1)/y_1 \in S_1\} \cup$$
$$\max\{d_{R_2}(c_2,y_2)/y_2 \in S_2\} \cup ... \cup$$
$$\max\{d_{R_m}(c_m,y_m)/y_m \in S_m\}$$

$$\leq \max\{d_{R_1}(c_1,u_1) + d_{R_1}(u_1,y_1)/y_1 \in S_1\} \cup$$
$$\max\{d_{R_2}(c_2,u_2) + d_{R_2}(u_2,y_2)/y_2 \in S_2\} \cup ... \cup$$
$$\max\{d_{R_m}(c_m,u_m) + d_{R_m}(u_m,y_m)/y_m \in S_m\}$$

$$= d_{R_1}(c_1,u_1) + \mathrm{cov}(u_1,S_1) \cup$$
$$d_{R_2}(c_2,u_2) + \mathrm{cov}(u_2,S_2) \cup ... \cup$$
$$d_{R_m}(c_m,u_m) + \mathrm{cov}(u_m,S_m)$$

$$\leq t_1^1(C_1) + t_{m_1}^1(V^{n_1}) \cup t_1^2(C_2) + t_{m_2}^2(V^{n_2}) \cup ... \cup$$
$$t_1^m(C_m) + t_{m_m}^m(V^{n_m}) .$$

Thus for every $S = S_1 \cup S_2 \cup ... \cup S_m \subseteq V^{n_1} \cup V^{n_2} \cup ... \cup V^{n_m}$ with $|S| = |S_1| \cup |S_2| \cup ... \cup |S_m| = m = m_1 \cup m_2 \cup ... \cup m_m$ one can find $c = c_1 \cup c_2 \cup ... \cup c_m \in C_1 \cup C_2 \cup ... \cup C_m$ such that
$$\mathrm{cov}(c,S) = \mathrm{cov}(c_1,S_1) \cup \mathrm{cov}(c_2,S_2) \cup ... \cup \mathrm{cov}(c_m,S_m)$$
$$\leq t_1^1(C_1) + t_{m_1}^1(V^{n_1}) \cup t_1^2(C_2) + t_{m_2}^2(V^{n_2}) \cup ... \cup$$
$$t_1^m(C_m) + t_{m_m}^m(V^{n_m}) .$$

Since
$$\mathrm{cov}(c,S) = \mathrm{cov}(c_1,S_1) \cup \mathrm{cov}(c_2,S_2) \cup ... \cup \mathrm{cov}(c_m,S_m)$$

$$= \min\{\mathrm{cov}(a_1,S_1)/a_1 \in C_1\} \cup$$
$$\min\{\mathrm{cov}(a_2,S_2)/a_2 \in C_2\} \cup ... \cup \min\{\mathrm{cov}(a_m,S_m)/a_m \in C_m\}$$

$$\leq \{t_1^1(C_1) + t_{m_1}^1(V^{n_1})\} \cup \{t_1^2(C_2) + t_{m_2}^2(V^{n_2})\} \cup ... \cup$$
$$\{t_1^m(C_m) + t_{m_m}^m(V^{n_m})\}$$

for all $S = S_1 \cup S_2 \cup ... \cup S_m \subseteq V^{n_1} \cup V^{n_2} \cup ... \cup V^{n_m}$ with $|S| = |S_1| \cup |S_2| \cup ... \cup |S_m| = m_1 \cup m_2 \cup ... \cup m_m$, it follows that

$$t_{m_1}^1(C_1) \cup t_{m_2}^2(C_2) \cup ... \cup t_{m_m}^m(C_m)$$

$$= \max\{cov(C_1, S_1)/S_1 \subseteq V^{n_1}, |S_1| = m_1\} \cup$$
$$\max\{cov(C_2, S_2)/S_2 \subseteq V^{n_2}, |S_2| = m_2\} \cup ... \cup$$
$$\max\{cov(C_m, S_m)/S_m \subseteq V^{n_m}, |S_m| = m_m\}$$

$$\leq t_1^1(C_1) + t_{m_1}^1(V^{n_1}) \cup t_1^2(C_2) + t_{m_2}^2(V^{n_2}) \cup ... \cup$$
$$t_1^m(C_m) + t_{m_m}^m(V^{n_m}) \ .$$

**PROPOSITION 3.7:** *For any m-tuple of integer $(n_1, n_2, ..., n_m)$; $n_1 \cup n_2 \cup ... \cup n_m \geq 2 \cup 2 \cup ... \cup 2$;*
$$t_2^1(V^{n_1}) \cup t_2^2(V^{n_2}) \cup ... \cup t_2^m(V^{n_m})$$
$$\leq n_1 - 1 \cup n_2 - 1 \cup ... \cup n_m - 1$$
*where $V^{n_i} = F_{2^N}^{n_i}$; $i = 1, 2, ..., n_i$; $n_i \leq N$, $i = 1, 2, ..., m$.*

*Proof:* Let
$$x_i = (x_1^i, x_2^i, ..., x_{n_i}^i) \text{ and } y_i = (y_1^i, y_2^i, ..., y_{n_i}^i) \in V^{n_i};$$
$i = 1, 2, ..., m$.

Let $u = u_1 \cup u_2 \cup ... \cup u_m \in V^{n_1} \cup V^{n_2} \cup ... \cup V^{n_m}$ where $u_i = (x_1^i, u_2^i, u_3^i, ..., u_{n_i-1}^i, y_{n_i}^i)$; $i = 1, 2, ..., m$. Thus $u = u_1 \cup u_2 \cup ... \cup u_m$ m-covers $x_1 \cup x_2 \cup ... \cup x_m$ and $y_1 \cup y_2 \cup ... \cup y_m \in V^{n_1} \cup V^{n_2} \cup ... \cup V^{n_m}$ within a m-radius $n_1-1 \cup n_2-1 \cup ... \cup n_m-1$ as $d_{R_i}(u_i, x_i) \leq n_i - 1$; $i = 1, 2, ..., m$. Thus for any pair of m-vectors $x_1 \cup x_2 \cup ... \cup x_m$, $y_1 \cup y_2 \cup ... \cup y_m$ in $V^{n_1} \cup V^{n_2} \cup ... \cup V^{n_m}$ there always exists a m-vector namely $u = u_1 \cup u_2 \cup ... \cup u_m$ which m-covers $x_1 \cup x_2 \cup ... \cup x_m$ and $y_1 \cup y_2 \cup ... \cup y_m$ within a m-radius $n_1-1 \cup n_2-1 \cup ... \cup n_m-1$.

Hence

$$t_2^1(V^{n_1}) \cup t_2^2(V^{n_2}) \cup ... \cup t_2^m(V^{n_m})$$
$$\leq n_1 - 1 \cup n_2 - 1 \cup ... \cup n_m - 1.$$

Now we proceed on to describe the notion of generalized sphere m-covering m-bounds for RD – m-codes. A natural question is for a given $t^1 \cup t^2 \cup ... \cup t^m$, $m_1 \cup m_2 \cup ... \cup m_m$ and $n_1 \cup n_2 \cup ... \cup n_m$ what is the smallest RD m-code whose $m_1 \cup m_2 \cup ... \cup m_m$, m-covering m-radius is atmost $t^1 \cup t^2 \cup ... \cup t^m$. As it turns out even for $m_1 \cup m_2 \cup ... \cup m_m \geq 2 \cup 2 \cup ... \cup 2$, it is necessary that $t^1 \cup t^2 \cup ... \cup t^m$ be atleast

$$\frac{n_1}{2} \cup \frac{n_2}{2} \cup ... \cup \frac{n_m}{2}.$$

Infact the minimal $t^1 \cup t^2 \cup ... \cup t^m$ for which such a m-code exists is the $(m_1, m_2, ..., m_m)$, m-covering m-radius of $C_1 \cup C_2 \cup ... \cup C_m = F_{2^N}^{n_1} \cup F_{2^N}^{n_2} \cup ... \cup F_{2^N}^{n_m}$. Various external values associated with this notion are $t_{m_1}^1(V^{n_1}) \cup t_{m_2}^2(V^{n_2}) \cup ... \cup t_{m_m}^m(V^{n_m})$, the smallest $(m_1, m_2, ..., m_m)$ covering m-radius among m-length $n_1 \cup n_2 \cup ... \cup n_m$ RD-m-codes

$$t_{m_1}^1(n_1, K_1) \cup t_{m_2}^2(n_2, K_2) \cup ... \cup t_{m_m}^m(n_m, K_m);$$

the smallest $(m_1, m_2, ..., m_m)$ covering m-radius among all $(n_1, K_1) \cup (n_2, K_2) \cup ... \cup (n_m, K_m)$ RD-m-codes.

$$K_{m_1}^1(n_1, t^1) \cup K_{m_2}^2(n_2, t^2) \cup ... \cup K_{m_m}^m(n_m, t^m)$$

is the smallest m-cardinality of a m-length $n_1 \cup n_2 \cup ... \cup n_m$. RD-m-code with $m_1 \cup m_2 \cup ... \cup m_m$ covering m-radius $t^1 \cup t^2 \cup ... \cup t^m$ and so on.

It is the latter quantity that is studied in the book for deriving new lower m-bounds. From the earlier results

$$K_{m_1}^1(n_1, t^1) \cup K_{m_2}^2(n_2, t^2) \cup ... \cup K_{m_m}^m(n_m, t^m)$$

is undefined if

$$t^1 \cup t^2 \cup \ldots \cup t^m < \frac{n_1}{2} \cup \frac{n_2}{2} \cup \ldots \cup \frac{n_m}{2}.$$

When this is the case it is accepted to say

$$K^1_{m_1}(n_1,t^1) \cup K^2_{m_2}(n_2,t^2) \cup \ldots \cup K^m_{m_m}(n_m,t^m)$$
$$= \infty \cup \infty \cup \ldots \cup \infty.$$

There are other circumstances when

$$K^1_{m_1}(n_1,t^1) \cup K^2_{m_2}(n_2,t^2) \cup \ldots \cup K^m_{m_m}(n_m,t^m)$$

is undefined. For instance

$$K^1_{2^{N_{n_1}}}(n_1,n_1-1) \cup K^2_{2^{N_{n_2}}}(n_2,n_2-1) \cup \ldots \cup K^m_{2^{N_{n_m}}}(n_m,n_m-1)$$
$$= \infty \cup \infty \cup \ldots \cup \infty.$$
$$m_1 > V(n_1,t^1),\ m_2 > V(n_2,t^2),\ldots,m_m > V(n_m,t^m);$$

since in this case no m-ball of m-radius $t^1 \cup t^2 \cup \ldots \cup t^m$ m-covers any m-set of $m_1 \cup m_2 \cup \ldots \cup m_m$ distinct m-vectors.

More generally one has the fundamental issue of whether
$$K^1_{m_1}(n_1,t^1) \cup K^2_{m_2}(n_2,t^2) \cup \ldots \cup K^m_{m_m}(n_m,t^m)$$
is m-finite for a given $n_1,m_1,t^1,\ n_2,m_2,t^2,\ldots,n_m,m_m,t^m$. This is the case if and only if
$$t^1_{m_1}(V^{n_1}) \leq t^1, t^2_{m_2}(V^{n_2}) \leq t^2,\ldots, t^m_{m_m}(V^{n_m}) \leq t^m$$

since
$$t^1_{m_1}(V^{n_1}) \cup t^2_{m_2}(V^{n_2}) \cup \ldots \cup t^m_{m_m}(V^{n_m})$$
lower m-bounds the $(m_1, m_2, \ldots, m_m)$ covering m-radius of all other m-codes of m-dimension $n_1 \cup n_2 \cup \ldots \cup n_m$ when $t^1 \cup t^2 \cup \ldots \cup t^m = n_1 \cup n_2 \cup \ldots \cup n_m$ every m-code word m-covers every m-vector, so a m-code of size $1 \cup 1 \cup \ldots \cup 1$ will $(m_1, m_2, \ldots, m_m)$ m-cover $V^{n_1} \cup V^{n_2} \cup \ldots \cup V^{n_m}$ for every $m_1 \cup m_2 \cup \ldots \cup m_m$. Thus
$$K^1_{m_1}(n_1,n_1) \cup K^2_{m_2}(n_2,n_2) \cup \ldots \cup K^m_{m_m}(n_m,n_m)$$
$$= 1 \cup 1 \cup \ldots \cup 1$$

for every $m_1 \cup m_2 \cup \ldots \cup m_m$.
If $t^1 \cup t^2 \cup \ldots \cup t^m = n_1 - 1 \cup n_2 - 1 \cup \ldots \cup n_m - 1$ what happens to $K_{m_1}^1(n_1, t^1) \cup K_{m_2}^2(n_2, t^2) \cup \ldots \cup K_{m_m}^m(n_m, t^m)$?

When $m_1 = m_2 = \ldots = m_m$,

$$K_2^1(n_1, n_1 - 1) \cup K_2^2(n_2, n_2 - 1) \cup \ldots \cup K_2^m(n_m, n_m - 1)$$
$$\leq 1 + L_{n_1}(n_1) \cup 1 + L_{n_2}(n_2) \cup \ldots \cup 1 + L_{n_m}(n_m).$$

For $\overline{0} \cup \overline{0} \cup \ldots \cup \overline{0} = (0, 0, \ldots, 0) \cup (0, 0, \ldots, 0) \cup (0, 0, \ldots, 0)$ will m-cover m-norm less than or equal to $n_1 - 1 \cup n_2 - 1 \cup \ldots \cup n_m - 1$ within m-radius $n_1 - 1 \cup n_2 - 1 \cup \ldots \cup n_m - 1$.
That is $\overline{0} \cup \overline{0} \cup \ldots \cup \overline{0} = (0, 0, \ldots, 0) \cup (0, 0, \ldots, 0) \cup (0, 0, \ldots, 0)$ will m-cover all m-norm $n_1 - 1 \cup n_2 - 1 \cup \ldots \cup n_m - 1$ m-vectors within the m-radius $n_1 - 1 \cup n_2 - 1 \cup \ldots \cup n_m - 1$.
Hence remaining m-vectors are m-rank $n_1 \cup n_2 \cup \ldots \cup n_m$ m-vectors.
Thus $\overline{0} \cup \overline{0} \cup \ldots \cup \overline{0} = (0, 0, \ldots, 0) \cup (0, 0, \ldots, 0) \cup (0, 0, \ldots, 0)$ and these m-rank $(n_1 \cup n_2 \cup \ldots \cup n_m)$ m-vector can m-cover the ambient m-space within the m-radius $n_1 - 1 \cup n_2 - 1 \cup \ldots \cup n_m - 1$. Therefore,

$$K_2^1(n_1, n_1 - 1) \cup K_2^2(n_2, n_2 - 1) \cup \ldots \cup K_2^m(n_m, n_m - 1)$$
$$\leq 1 + L_{n_1}(n_1) \cup 1 + L_{n_2}(n_2) \cup \ldots \cup 1 + L_{n_m}(n_m).$$

**PROPOSITION 3.8:** *For any RD m-code of m-length $n_1 \cup n_2 \cup \ldots \cup n_m$ over $F_{2^N} \cup F_{2^N} \cup \ldots \cup F_{2^N}$,*

$$K_2^1(n_1, n_1 - 1) \cup K_2^2(n_2, n_2 - 1) \cup \ldots \cup K_2^m(n_m, n_m - 1)$$
$$\leq m_1 L_{n_1}(n_1) + 1 \cup m_2 L_{n_2}(n_2) + 1 \cup \ldots \cup m_m L_{n_m}(n_m) + 1$$

*provided $m_1 \cup m_2 \cup \ldots \cup m_m$ is such that*
$$m_1 L_{n_1}(n_1) + 1 \cup m_2 L_{n_2}(n_2) + 1 \cup \ldots \cup m_m L_{n_m}(n_m) + 1$$
$$\leq |V^{n_1}| \cup |V^{n_2}| \cup \ldots \cup |V^{n_m}|.$$

*Proof:* Consider a RD m-code $C = C_1 \cup C_2 \cup \ldots \cup C_m$ such that

$$|C| = |C_1| \cup |C_2| \cup ... \cup |C_m|$$
$$= m_1 L_{n_1}(n_1) + 1 \cup m_2 L_{n_2}(n_2) + 1 \cup ... \cup m_m L_{n_m}(n_m) + 1.$$

Each m-vector in $V^{n_1} \cup V^{n_2} \cup ... \cup V^{n_m}$ has
$$L_{n_1}(n_1) \cup L_{n_2}(n_2) \cup ... \cup L_{n_m}(n_m)$$
rank complements, that is from each m-vector $v_1 \cup v_2 ... \cup v_m \in V^{n_1} \cup V^{n_2} \cup ... \cup V^{n_m}$ there are
$$L_{n_1}(n_1) \cup L_{n_2}(n_2) \cup ... \cup L_{n_m}(n_m)$$
m-vectors at rank m-distance $n_1 \cup n_2 \cup ... \cup n_m$.
This means for any set
$$S_1 \cup S_2 \cup ... \cup S_m \subseteq V^{n_1} \cup V^{n_2} \cup ... \cup V^{n_m}$$
of $(m_1, m_2, ..., m_m)$ m-vectors there always exists a $c_1 \cup c_2 \cup ... \cup c_m \in C_1 \cup C_2 \cup ... \cup C_m$ which m-covers $S_1 \cup S_2 \cup ... \cup S_m$; m-rank distance $n_1 - 1 \cup n_2 - 1 \cup ... \cup n_m - 1$. Thus

$$cov(u_1, S_1) \cup cov(u_2, S_2) \cup ... \cup cov(u_m, S_m)$$
$$\leq n_1 - 1 \cup n_2 - 1 \cup ... \cup n_m - 1$$

which implies
$$cov(u_1, S_1) \cup cov(u_2, S_2) \leq n_1 - 1 \cup n_2 - 1.$$

Hence
$$K^1_{m_1}(n_1, n_1 - 1) \cup K^2_{m_2}(n_2, n_2 - 1) \cup ... \cup K^m_{m_m}(n_m, n_m - 1)$$
$$\leq m_1 L_{n_1}(n_1) + 1 \cup m_2 L_{n_2}(n_2) + 1 \cup ... \cup m_m L_{n_m}(n_m) + 1.$$

By bounding the number of $(m_1, m_2, ..., m_m)$ m-sets that can be covered by a given m-code word, one obtains a straight forward generalization of the classical sphere m-bound.

**THEOREM** *(Generalized sphere bound for RD-m-codes): For any* $(n_1, K_1) \cup (n_2, K_2) \cup ... \cup (n_m, K_m)$ *RD m-code* $C = C_1 \cup C_2 \cup ... \cup C_m$,

$$K^1 \binom{V(n_1, t_{m_1}(C_1))}{m_1} \cup K^2 \binom{V(n_2, t_{m_2}(C_2))}{m_2}$$
$$\cup ... \cup K^m \binom{V(n_m, t_{m_m}(C_m))}{m_m}$$

$$\geq \binom{2^{N_{n_1}}}{m_1} \cup \binom{2^{N_{n_2}}}{m_2} \cup \ldots \cup \binom{2^{N_{n_m}}}{m_m}.$$

*Hence for any $n_i, t_i, m_i$ triple $i = 1, 2, \ldots, m$*

$$K_{m_1}^1(n_1, t_1) \cup K_{m_2}^2(n_2, t_2) \cup \ldots \cup K_{m_m}^m(n_m, t_m)$$

$$\geq \frac{\binom{2^{N_{n_1}}}{m_1}}{\binom{V(n_1,t_1)}{m_1}} \cup \frac{\binom{2^{N_{n_2}}}{m_2}}{\binom{V(n_2,t_2)}{m_2}} \cup \ldots \cup \frac{\binom{2^{N_{n_m}}}{m_m}}{\binom{V(n_m,t_m)}{m_m}}$$

*where*

$$V(n_1, t_1) \cup V(n_2, t_2) \cup \ldots \cup V(n_m, t_m)$$
$$= \sum_{i_1=0}^{t_1} L_{i_1}^1(n_1) \cup \sum_{i_2=0}^{t_2} L_{i_2}^2(n_2) \cup \ldots \cup \sum_{i_m=0}^{t_m} L_{i_m}^m(n_m)$$

*number of m-vectors in a sphere m-radius $t^1 \cup t^2 \cup \ldots \cup t^m$ and $L_{i_1}^1(n_1) \cup L_{i_2}^2(n_2) \cup \ldots \cup L_{i_m}^m(n_m)$ is the number of m-vectors in $V^{n_1} \cup V^{n_2} \cup \ldots \cup V^{n_m}$ whose rank m-norm is $i_1 \cup i_2 \cup \ldots \cup i_m$.*

*Proof:* Each set of $(m_1, m_2, \ldots, m_m)$ m-vectors in $V^{n_1} \cup V^{n_2} \cup \ldots \cup V^{n_m}$ must occur in a sphere of m-radius $t_{m_1}(C_1) \cup t_{m_2}(C_2) \cup \ldots \cup t_{m_m}(C_m)$ around at least one code m-word. Total number of such m-sets is

$$|V^{n_1}| \cup |V^{n_2}| \cup \ldots \cup |V^{n_m}|$$

choose $m_1 \cup m_2 \cup \ldots \cup m_m$ where

$$|V^{n_1}| \cup |V^{n_2}| \cup \ldots \cup |V^{n_m}| = 2^{Nn_1} \cup 2^{Nn_2} \cup \ldots \cup 2^{Nn_m}.$$

The number of m-sets of $(m_1, m_2, \ldots, m_m)$ m-vectors in a neighborhood of m-radius

$$t_{m_1}(C_1) \cup t_{m_2}(C_2) \cup \ldots \cup t_{m_m}(C_m)$$

is

$$V(n_1, t_{m_1}(C_1)) \cup V(n_2, t_{m_2}(C_2)) \cup \ldots \cup V(n_m, t_{m_m}(C_m))$$

choose $m_1 \cup m_2 \cup \ldots \cup m_m$.
There are K-code m-words.
Hence

$$K^1 \binom{V(n_1, t_{m_1}(C_1))}{m_1} \cup K^2 \binom{V(n_2, t_{m_2}(C_2))}{m_2}$$
$$\cup \ldots \cup K^m \binom{V(n_m, t_{m_m}(C_m))}{m_m}$$

$$\geq \binom{2^{N_{n_1}}}{m_1} \cup \binom{2^{N_{n_2}}}{m_2} \cup \ldots \cup \binom{2^{N_{n_m}}}{m_m}.$$

Thus for any $n_1 \cup n_2 \cup \ldots \cup n_m$, $t^1 \cup t^2 \cup \ldots \cup t^m$ and $m_1 \cup m_2 \cup \ldots \cup m_m$

$$K^1_{m_1}(n_1, t_1) \cup K^2_{m_2}(n_2, t_2) \cup \ldots \cup K^m_{m_m}(n_m, t_m)$$

$$\geq \frac{\binom{2^{N_{n_1}}}{m_1}}{\binom{V(n_1, t_1)}{m_1}} \cup \frac{\binom{2^{N_{n_2}}}{m_2}}{\binom{V(n_2, t_2)}{m_2}} \cup \ldots \cup \frac{\binom{2^{N_{n_m}}}{m_m}}{\binom{V(n_m, t_m)}{m_m}}.$$

**COROLLARY 3.2:** *If*

$$\binom{2^{N_{n_1}}}{m_1} \cup \binom{2^{N_{n_2}}}{m_2} \cup \ldots \cup \binom{2^{N_{n_m}}}{m_m}$$
$$> 2^{N_{n_1}} \binom{V(n_1, t_1)}{m_1} \cup 2^{N_{n_2}} \binom{V(n_2, t_2)}{m_2} \cup \ldots \cup 2^{N_{n_m}} \binom{V(n_m, t_m)}{m_m}$$

*then*

$K^1_{m_1}(n_1, t_1) \cup K^2_{m_2}(n_2, t_2) \cup \ldots \cup K^m_{m_m}(n_m, t_m) = \infty \cup \infty \cup \ldots \cup \infty$.

**Chapter Four**

# APPLICATIONS OF RANK DISTANCE m-CODES

In this chapter we proceed onto give some applications of Rank Distance bicodes and those of the new classes of rank distance bicodes and their generalizations.

Rank distance m-codes can be used in multi disk storage systems by constructing or building m- Redundant Array of Inexpensive Disks (m-RAID). These linear MRD m-codes can also be used in m-public key m-cryptosystems.

Circulant rank m-codes can be used in multi communication channels or m-channels having very high m-error m-probability for m-error correction.

The AMRD m-codes is useful for m-error correction in data multi(m-) storage systems. These m-codes will equally be as good as the MRD m-codes and are better than the corresponding m-codes in the Hamming metric.

In data multi storage systems these MRD and AMRD m-codes can be used simultaneously in criss cross error corrections by suitably programming; appropriately the functioning of the implementation.

Using these m-codes one can save time, space and economy. With computerization in every walk of life these m-codes can perform simultaneously bulk m-error correction in bulk data transmission if appropriately programmed.

Interested researcher can develop m-algorithms (algorithms that can work in m-channels simultaneously) of these rank distance m-codes.

# FURTHER READING

1. Adams M.J., *Subcodes and covering radius*, IEEE Trans. Inform. Theory, Vol.IT-32, pp.700-701, 1986.

2. Andrews G.E., *The Theory of Partitions*, Vol. 2 of Encyclopaedia of Mathematics and Its Applications, Addison – Wesley, Reading, Massachusetts, 1976.

3. Assmuss Jr. E.F. and Pless V., *On the Covering radius of extremal self dual codes*, IEEE Trans. Inform. Theory, Vol.IT-29, pp.359-363, 1983.

4. Bannai, E., and Ito, T., *Algebraic Combinatorics I. Association Schemes*, The Benjamin/Cummings Publishing Company, 1983.

5. Berger, T. P., *Isometries for rank distance and permutation group of Gabidulin codes*, IEEE Transactions on Information Theory, Vol. 49, no. 11, pp. 3016–3019, 2003.

6. Berlekamp E.R., *Algebraic Coding Theory*, McGraw Hill, 1968.

7. Best M.R., *Perfect codes hardly exist*, IEEE Trans. Inform. Theory, Vol.IT-29, pp.349-351, 1983.

8. Blahut, R., *Theory and Practice of Error Control Codes*. Addison-Wesley, 1983.

9. Blaum M. and McElice R.J., *Coding protection for magnetic tapes; a generalization of the Patel-Hong code*, IEEE Trans. Inform. Theory, Vol.IT-31, pp.690-693, 1985.

10. Blaum M., *A family of efficient burst-correcting array codes*, IEEE Trans. Inform. Theory, Vol.IT-36, pp.671-674, 1990.

11. Blaum M., Farrel F.G. and Van Tilborg H.C.A., *A class of error correcting array codes*, IEEE Trans. Inform. Theory, Vol.IT-32, pp.836-839, 1989.

12. Brouwer A.E., *A few new constant weight codes*, IEEE Trans. Inform. Theory, Vol.IT-26, p.366, 1980.

13. Carlitz, L., *q-Bernoulli numbers and polynomials*, Duke Mathematical Journal, Vol. 15, no. 4, pp. 987–1000, 1948.

14. Chabaud, F., and Stern, J., *The cryptographic security of the syndrome decoding problem for rank distance codes*, Advances in Cryptology: ASIACRYPT '96, Volume 1163 of LNCS, pp. 368–381, Springer-Verlag, 1996.

15. Chen, K., *On the non-existence of perfect codes with rank distance*, Mathematische Nachrichten, Vol. 182, no. 1, pp. 89–98, 1996.

16. Clark, W.E., and Dunning, L.A., *Tight upper bounds for the domination numbers of graphs with given order and minimum degree*, The Electronic Journal of Combinatorics, Vol. 4, 1997.

17. Cohen G.D., *A non constructive upper bound on covering radius*, IEEE Trans. Inform. Theory, Vol.IT-29, pp.352-353, 1983.

18. Cohen G.D., Karpovski M.R., Mattson Jr. H.F. and Schatz J.R., *Covering Radius – Survey and Recent Results*, IEEE Trans. Inform., Theory, Vol.IT-31, pp.328-343, 1985.

19. Cohen G.D., Lobstein A.C. and Sloane N.J.A., *Further results on the covering radius of codes*, IEEE Trans. Inform. Theory, Vol.IT-32, pp.680-694, 1986.

20. Delsarte P., *Four fundamental parameters of a code and their combinatorial significance*, Inform. and Control, Vol.23, pp.407-438, 1973.

21. Delsarte, P., *Bilinear forms over a finite field, with applications to coding theory*, Journal of Combinatorial Theory A, Vol. 25, no. 3, pp. 226–241, 1978.

22. Delsarte, P., *Properties and applications of the recurrence F(i + 1, k + 1, n + 1) = qk+1F(i, k + 1, n) − qkF(i, k, n)*, SIAM Journal of Applied Mathematics, Vol. 31, no. 2, pp. 262–270, September 1976.

23. Erickson T., *Bounds on the size of a code*, Lecture notes in Control and Inform. Sciences, Vol.128, pp.45-72, Springer – Verlag, 1989.

24. Gabidulin E.M., *Theory of codes with Maximum Rank Distance*, Problemy Peredachi Informatsii, Vol.21, pp.17-27, 1985.

25. Gabidulin, E. M. and Loidreau, P., *On subcodes of codes in the rank metric*, Proc. IEEE Int. Symp. on Information Theory, pp. 121–123, Sept. 2005.

26. Gabidulin, E. M. and Obernikhin, V.A., *Codes in the Vandermonde F-metric and their application*, Problems of Information Transmission, Vol. 39, no. 2, pp. 159–169, 2003.

27. Gabidulin, E. M., *Optimal codes correcting lattice-pattern errors*, Problems of Information Transmission, Vol. 21, no. 2, pp. 3–11, 1985.

28. Gabidulin, E. M., Paramonov, A.V., and Tretjakov, O.V., Ideals over a non-commutative ring and their application in cryptology," in *Proceedings of the Workshop on the Theory and Application of Cryptographic Techniques (EUROCRYPT '91)*, Vol. 547 of *Lecture Notes in Computer Science*, pp. 482–489, Brighton, UK, April 1991.

29. Gabidulin, E. M., *Public-key cryptosystems based on linear codes over large alphabets: efficiency and weakness*, Codes and Cyphers, pp. 17– 31, 1995.

30. Gadouleau, M., and Yan, Z*., Packing and covering properties of rank metric codes,* IEEE Trans. Info. Theory, Vol. 54, no. 9,pp. 3873–3883, September 2008.

31. Gadouleau, M., and Yan, Z., "Properties of codes with the rank metric," in *Proceedings of IEEE Global Telecommunications Conference (GLOBECOM '06)*, pp. 1–5, San Francisco, Calif, USA, November 2006.

32. Gadouleau, M., and Yan, Z., *Construction and covering properties of constant-dimension codes*, submitted to IEEE Trans. Info. Theory, 2009.

33. Gadouleau, M., and Yan, Z., Decoder error probability of MRD codes, in *Proceedings of IEEE Information Theory Workshop (ITW '06)*, pp. 264–268, Chengdu, China, October 2006.

34. Gadouleau, M., and Yan, Z., *On the decoder error probability of bounded rank-distance decoders for maximum rank distance codes*, to appear in IEEE Transactions on Information Theory.

35. Gasper, G., and Rahman, M., *Basic Hypergeometric Series*, Vol. 96 of Encyclopedia of Mathematics and Its Applications, Cambridge University Press, New York, NY, USA, 2nd edition, 2004.

36. Gibson, J.K., *The security of the Gabidulin public-key*

*cryptosystem*, in EUROCRYPT'96, pp. 212-223. 1996

37. Graham R.L. and Sloane N.J.A., *Lower Bounds for constant weight codes*, IEEE Trans. Inform. Theory, Vol.IT-26, pp.37-43, 1980.

38. Graham R.L. and Sloane N.J.A., *On the covering radius of codes*, IEEE Trans. Inform. Theory, Vol.IT-31, pp.385-401, 1985.

39. Grant, D. and Varanasi, M., *Duality theory for space-time codes over finite fields*, to appear in Advance in Mathematics of Communications.

40. Grant, D. and Varanasi, M., Weight enumerators and a MacWilliams-type identity for space-time rank codes over finite fields, in *Proceedings of the 43rd Allerton Conference on Communication, Control, and Computing*, pp. 2137–2146, Monticello, Ill, USA, October 2005.

41. Hall, L.O., and Gur Dial, *On fuzzy codes for Asymmetric and Unidirectional error,* Fuzzy sets and Systems, 36, pp. 365-373, 1990.

42. Harold N.Ward, *Divisible codes*, Archiv der Mathematic, 36, pp. 485-494, 1981.

43. Helgert H.J. and Stinaff R.D., *Minimum distance bounds for binary linear codes*, IEEE Trans. Inform. Theory, Vol.IT-19, pp.344-356, 1973.

44. Helleseth T., Klove T. and Mykkeltveit J., *On the covering radius of binary codes*, IEEE Trans. Inform. Theory, Vol.IT-24, pp.627-628, 1978.

45. Helleseth T., *The weight distribution of the coset leaders for some classes of codes with related parity check matrices*, Discrete Math., Vol.28, pp.161-171, 1979.

46. Herstein I.N., *Topics in Algebra*, 2$^{nd}$, Wiley, New York, 1975.

47. Hill R., *A first course in coding theory*, Clarendon, Oxford, 1986.

48. Hoffman K. and Kunze R., *Linear Algebra*, 2$^{nd}$, Prentice Hall of India, New Delhi, 1990.

49. Hua, L., *A theorem on matrices over a field and its applications*, Chinese Mathematical Society, Vol. 1, no. 2, pp. 109–163, 1951.

50. Janwa H.L., *Some new upper bounds on the covering radius of binary linear codes*, IEEE Trans. Inform. Theory, Vol.IT-35, pp.110-122, 1989.

51. Johnson S.M., *Upper bounds for constant weight error correcting codes*, Discrete Math., Vol.3, pp.109-124, 1972.

52. Koetter, R., and Kschischang, F.R., *Coding for errors and erasures in random network coding*, IEEE Trans. Info. Theory, Vol. 54, no. 8, pp. 3579–3591, August 2008.

53. Kshevetskiy, A. and Gabidulin, E. M., *The new construction of rank codes*, Proc. IEEE Int. Symp. on Information Theory, pp. 2105–2108, Sept. 2005.

54. Lidl R. and Niederreiter H., *Finite fields*, Addison-Wesley, Reading, Massachusetts, 1983.

55. Lidl R. and Pilz G., *Applied Abstract Algebra*, Springer-Verlag, New York, 1984.

56. Loidreau, P., *A Welch-Berlekamp like algorithm for decoding Gabidulin codes*, in Proceedings of the 4th International Workshop on Coding and Cryptography (WCC '05), Vol. 3969, pp. 36–45, Bergen, Norway, March 2005.

57. Loidreau, P., *Etude et optimisation de cryptosyst`emes `a cl´e publique fond´es sur la th´eorie des codes correcteurs*, Ph.D. Dissertation, Ecole Polytechnique, Paris, France, May 2001.

58. Loidreau, P., *Properties of codes in rank metric*, http://arxiv.org/pdf/cs.DM/0610057/.

59. Lusina, P., Gabidulin, E., and Bossert, M., *Maximum Rank Distance codes as space-time codes*, IEEE Trans. Info. Theory, Vol. 49, pp. 2757–2760, Oct. 2003.

60. MacLane S. and Birkhoff G., *Algebra*, Macmillan, New York, 1967.

61. MacWilliams, F., and Sloane, N., *The Theory of Error-Correcting Codes*, North-Holland, Amsterdam, The Netherlands, 1977.

62. Manoj, K.N., and Sundar Rajan, B., *Full Rank Distance codes*, Technical Report, IISc Bangalore, Oct. 2002.

63. Mattson Jr. H.F., *An improved upper bound on the covering radius*, Lecture notes in Computer Science, 228, Springer – Verlag, New York, pp.90-106, 1986.

64. Mattson Jr. H.F., *An upper bound on covering radius*, Ann. Discrete Math., Vol.17, pp.453-458, 1983.

65. Mattson Jr. H.F., *Another upper bound on covering radius*, IEEE Trans. Inform. Theory, Vol.IT-29, pp.356-359, 1983.

66. McEliece R.J., and Swanson, L., *On the decoder error probability for Reed-Solomon codes*, IEEE Trans. Info. Theory, Vol. 32, no. 5, pp. 701–703, Sept. 1986.

67. McEliece, R.J., *A public-key cryptosystem based on algebraic coding theory*, Jet Propulsion Lab. DSN Progress Report, pp. 114–116, 1978.

68. Ourivski, A. V. and Gabidulin, E. M., *Column scrambler for the GPT cryptosystem*, Discrete Applied Mathematics, Vol. 128, pp. 207–221, 2003.

69. Peterson W.W. and Weldon Jr. E.J., *Error correcting codes*, M.I.T. Press, Cambridge, Massachusetts, 1972.

70. Pless, V., *Power moment identities on weight distributions in error correcting codes*, Information and Control, Vol. 6, no. 2, pp. 147–152, 1963.

71. Prunsinkiewicz P. and Budkowski S., *A double track error correction code for magnetic tape*, IEEE Trans. Comput., Vol.C-25, pp.642-645, 1976.

72. Richter, G., and Plass, S., Fast decoding of rank-codes with rank errors and column erasures," in *Proceedings of IEEE International Symposium on Information Theory (ISIT '04)*, p.398, Chicago, Ill, USA, June-July 2004.

73. Roth R.M., *Maximum Rank Array codes and their application to crisscross error correction*, IEEE Trans. Inform. Theory, Vol.IT-37, pp.328-336, 1991.

74. Roth, R. M., *Probabilistic criss cross error correction*, IEEE Transactions on Information Theory, Vol. 43, no. 5, pp. 1425–1438, 1997.

75. Schonheim J., *On linear and non-linear single error correcting q-nary perfect codes*, Inform. and Control, Vol.12, pp.23-26, 1968.

76. Schwartz, M., and Etzion, T., *Two-dimensional cluster correcting codes*, IEEE Transactions on Information Theory, Vol. 51, no. 6, pp. 2121–2132, 2005.

77. Selvaraj, R.S., *Metric properties of rank distance codes*, Ph.D. Thesis, Thesis Guide; Dr. W. B. Vasantha Kandasamy, Department of Mathematics, Indian Institute of Technology, Madras, .2006.

78. Shannon C.E., *A mathematical theory of communication*, Bell. Syst. Tech., Jl.27, pp.379-423 and pp.623-656, 1948.

79. Silva, D., Kschischang, F. R. and Koetter, R., *A rank-metric approach to error control in random network coding*, IEEE Trans. Info. Theory, Vol. 54, no. 9, pp. 3951–3967, September 2008.

80. Sripati, U., and Sundar Rajan, B., *On the rank distance of cyclic codes*, Proc. IEEE Int. Symp. on Information Theory, p. 72, July 2003.

81. Strang G., *Linear Algebra and its applications*, 2$^{nd}$, Academic Press, New York, 1980.

82. Suresh Babu, N., S*tudies on Rank Distance Codes*, Ph.D. Thesis. Thesis Guide; Dr. W.B.Vasantha Kandasamy, Department of Mathematics, Indian Institute of Technology, Madras, 1995.

83. Tarokh, V., Seshadri, N., and Calderbank, A.R., *Space-time codes for high data rate wireless communication: performance criterion and code construction*, IEEE Transactions on Information Theory, Vol. 44, no. 2, pp. 744–765, 1998.

84. Tietavanien A., *On the non-existence of perfect codes over finite fields*, Siam. Jl. Appl. Math., Vol.24, 1973.

85. Van Lint J.H. and Wilson R.M, *A course in combinatories*, Cambridge University Press, 1992.

86. Van Lint J.H., *A Survey of perfect codes*, Rocky Moundain Jl. of Math., Vol.5, pp.199-224, 1975.

87. Van Lint J.H., *Coding Theory*, Springer-Verlag, New York, 1971.

88. Vasantha Kandasamy, W.B., and Raja Durai, R. S., *Maximum rank distance codes with complementary duals*, Mathematics and Information Theory, Recent topics and Applications, pp. 86-90, 2002.

89. Vasantha Kandasamy, W.B., and Selvaraj, R. S., *Multi-covering radii of codes with rank metric*, Proc. Information Theory Workshop, p. 215, Oct. 2002.

90. Vasantha Kandasamy, W.B., and Suresh Babu, N., *On the covering radius of rank-distance codes,* Ganita Sandesh, Vol. 13, pp. 43–48, 1999.

91. Vasantha Kandasamy, W.B., *Bialgebraic structures and Smarandache bialgebraic structures*, American Research Press, Rehoboth, 2003.

92. Vasantha Kandasamy, W.B., Smarandache, Florentin and K. Ilanthenral, *Introduction to bimatrices*, Hexis, Phoenix, 2005.

93. Vasantha Kandasamy, W.B., Smarandache, Florentin and K. Ilanthenral, *Introduction to Linear Bialgebra*, Hexis, Phoenix, 2005.

94. Vasantha, W. B. and R.S.Selvaraj, *Multi-covering Radii of Codes with Rank Metric*, Proceedings of 2002 IEEE Information Theory Workshop, Indian Institute of Science, Bangalore, October 20-25, p. 215, 2002.

95. Vasantha, W.B. and Selvaraj R.S., *Divisible MRD Codes,* Proceedings of National Conference on Challenges of the 21$^{st}$ century in Mathematics and its Allied Topics, University of Mysore, Karnataka, India, pp., 242-246, February 3-4, 2001.

96. Vasantha, W.B. and Selvaraj, R.S. *Fuzzy RD Codes with Rank Metric and Distance Properties*, Proceedings of the International Conference on Recent Advances in Mathematical Sciences, Indian Institute of Technology Kharagpur, India, pp., 341-348, Narosa Publishing House, New Delhi, December 20-22, 2000

97. Vasantha, W.B. and Selvaraj, R.S. *Multicovering radius of Codes for Rank Metric*, Proceedings of 4$^{th}$ Asia Europe Workshop on Information Theory Concepts, Viareggio, Italy, October 6-8, pp., 96-98, 2004.

98. Vasantha Kandasamy, W.B., and Smarandache, Florentin, *Rank distance m-codes and their properties (To appear)*, Maths Tiger.

99. Von Kaenel, P.A, *Fuzzy codes and distance properties*, Fuzzy sets and systems, 8, pp., 199-204, 1982.

100. Zhen Zhang, *Limiting efficiencies of Burst correcting Array codes*, IEEE Trans. Inform. Theory, Vol.IT-37, pp.976-982, 1991.

101. Zimmermann K.H. *The weight distribution of indecomposable cyclic codes over 2-groups*, Jour. of Comb. Theory, Series A, Vol.60, 1992.

102. Zinover V.A. and Leontev V.K., *The non-existence of perfect codes over Galois fields,* Problems of Control and Info. Theory, Vol.2, pp.123-132, 1973.

# INDEX

# C

## F

## G

## H

## L

## M

# Q

# R

# S

# ABOUT THE AUTHORS

**Dr.W.B.Vasantha Kandasamy** is an Associate Professor in the Department of Mathematics, Indian Institute of Technology Madras, Chennai. In the past decade she has guided 13 Ph.D. scholars in the different fields of non-associative algebras, algebraic coding theory, transportation theory, fuzzy groups, and applications of fuzzy theory of the problems faced in chemical industries and cement industries.

She has to her credit 646 research papers. She has guided over 68 M.Sc. and M.Tech. projects. She has worked in collaboration projects with the Indian Space Research Organization and with the Tamil Nadu State AIDS Control Society. She is presently working on a research project funded by the Government of India's Department of Atomic Energy. This is her $46^{th}$ book.

On India's 60th Independence Day, Dr.Vasantha was conferred the Kalpana Chawla Award for Courage and Daring Enterprise by the State Government of Tamil Nadu in recognition of her sustained fight for social justice in the Indian Institute of Technology (IIT) Madras and for her contribution to mathematics. The award, instituted in the memory of Indian-American astronaut Kalpana Chawla who died aboard Space Shuttle Columbia, carried a cash prize of five lakh rupees (the highest prize-money for any Indian award) and a gold medal.
She can be contacted at vasanthakandasamy@gmail.com
Web Site: http://mat.iitm.ac.in/home/wbv/public_html/

---

**Dr. Florentin Smarandache** is a Professor of Mathematics at the University of New Mexico in USA. He published over 75 books and 150 articles and notes in mathematics, physics, philosophy, psychology, rebus, literature.

In mathematics his research is in number theory, non-Euclidean geometry, synthetic geometry, algebraic structures, statistics, neutrosophic logic and set (generalizations of fuzzy logic and set respectively), neutrosophic probability (generalization of classical and imprecise probability). Also, small contributions to nuclear and particle physics, information fusion, neutrosophy (a generalization of dialectics), law of sensations and stimuli, etc. He can be contacted at smarand@unm.edu

---

**Dr. N. Suresh Babu** is a faculty in the College of Engineering, Trivandrum. He has obtained his doctorate degree in coding/ communication theory from department of mathematics, Indian Institute of Technology, Madras. He has worked as a post doctoral research fellow, under Prof. Zimmermann, K. H., a renowned coding theorist, working in the University of Bayreuth, Germany.

---

**Dr. R. S. Selvaraj** is presently working as a faculty in NIIT Warangal. He has obtained his doctor of philosophy in coding/ communication theory from department of mathematics, Indian Institute of Technology, Madras.